\magnification=\magstephalf
 \ifnum\mag=\magstephalf
  \hoffset=0.52truein \voffset=0.125truein  \else
 \ifnum\mag=\magstep0
  \hoffset=0.75truein \voffset=0.30truein   \else
 \ifnum\mag=\magstep1
  \hoffset=0.25truein \voffset=0.375truein \else
 \fi\fi\fi
\input amstex
\documentstyle{amsppt}

\def\next{AMSPPT}\ifx\styname\next \else\input amsppt.sty \relax\fi



\brokenpenalty=10000
\clubpenalty=10000
\widowpenalty=10000

\catcode`\@=11

\def\pagewidth#1{\hsize#1%
   \captionwidth@24pc}

\pagewidth{30pc}

\parindent=18\p@
\normalparindent\parindent

\parskip=0pt

\def\makefootline{\baselineskip18\p@\line{\the\footline}}

\def\foliofont@{\sevenrm}
\def\headlinefont@{\sevenpoint}


\def\specialheadfont@{\elevenpoint\smc}
\def\headfont@{\bf}
\def\subheadfont@{\bf}
\def\refheadfont@{\bf}
\def\abstractfont@{\smc}
\def\proclaimheadfont@{\smc}
\def\xcaheadfont@{\smc}
\def\captionfont@{\smc}
\def\citefont@{\bf}
\def\refsfont@{\eightpoint}

\font\sixsy=cmsy6

\font@\twelverm=cmr10 scaled \magstep1
\font@\twelvebf=cmbx10 scaled \magstep1
\font@\twelveit=cmti10 scaled \magstep1
\font@\twelvesl=cmsl10 scaled \magstep1
\font@\twelvesmc=cmcsc10 scaled \magstep1
\font@\twelvett=cmtt10 scaled \magstep1
\font@\twelvei=cmmi10 scaled \magstep1
\font@\twelvesy=cmsy10 scaled \magstep1
\font@\twelveex=cmex10 scaled \magstep1

\newtoks\twelvepoint@
\def\twelvepoint{\normalbaselineskip14\p@
 \abovedisplayskip12\p@ plus3\p@ minus9\p@
 \belowdisplayskip\abovedisplayskip
 \abovedisplayshortskip\z@ plus3\p@
 \belowdisplayshortskip7\p@ plus3\p@ minus4\p@
 \textonlyfont@\rm\twelverm \textonlyfont@\it\twelveit
 \textonlyfont@\sl\twelvesl \textonlyfont@\bf\twelvebf
 \textonlyfont@\smc\twelvesmc \textonlyfont@\tt\twelvett
 \ifsyntax@ \def\big##1{{\hbox{$\left##1\right.$}}}%
  \let\Big\big \let\bigg\big \let\Bigg\big
 \else
   \textfont\z@\twelverm  \scriptfont\z@\eightrm
       \scriptscriptfont\z@\sixrm
   \textfont\@ne\twelvei  \scriptfont\@ne\eighti
       \scriptscriptfont\@ne\sixi
   \textfont\tw@\twelvesy \scriptfont\tw@\eightsy
       \scriptscriptfont\tw@\sixsy
   \textfont\thr@@\twelveex \scriptfont\thr@@\eightex
        \scriptscriptfont\thr@@\eightex
   \textfont\itfam\twelveit \scriptfont\itfam\eightit
        \scriptscriptfont\itfam\eightit
   \textfont\bffam\twelvebf \scriptfont\bffam\eightbf
        \scriptscriptfont\bffam\sixbf
   \setbox\strutbox\hbox{\vrule height8.5\p@ depth3.5\p@ width\z@}%
   \setbox\strutbox@\hbox{\lower.5\normallineskiplimit\vbox{%
        \kern-\normallineskiplimit\copy\strutbox}}%
   \setbox\z@\vbox{\hbox{$($}\kern\z@}\bigsize@1.2\ht\z@
  \fi
  \normalbaselines\rm\dotsspace@1.5mu\ex@.2326ex\jot3\ex@
  \the\twelvepoint@}

\font@\elevenrm=cmr10 scaled \magstephalf
\font@\elevenbf=cmbx10 scaled \magstephalf
\font@\elevenit=cmti10 scaled \magstephalf
\font@\elevensl=cmsl10 scaled \magstephalf
\font@\elevensmc=cmcsc10 scaled \magstephalf
\font@\eleventt=cmtt10 scaled \magstephalf
\font@\eleveni=cmmi10 scaled \magstephalf
\font@\elevensy=cmsy10 scaled \magstephalf
\font@\elevenex=cmex10 scaled \magstephalf

\newtoks\elevenpoint@
\def\elevenpoint{\normalbaselineskip13\p@
 \abovedisplayskip12\p@ plus3\p@ minus9\p@
 \belowdisplayskip\abovedisplayskip
 \abovedisplayshortskip\z@ plus3\p@
 \belowdisplayshortskip7\p@ plus3\p@ minus4\p@
 \textonlyfont@\rm\elevenrm \textonlyfont@\it\elevenit
 \textonlyfont@\sl\elevensl \textonlyfont@\bf\elevenbf
 \textonlyfont@\smc\elevensmc \textonlyfont@\tt\eleventt
 \ifsyntax@ \def\big##1{{\hbox{$\left##1\right.$}}}%
  \let\Big\big \let\bigg\big \let\Bigg\big
 \else
   \textfont\z@\elevenrm  \scriptfont\z@\eightrm
       \scriptscriptfont\z@\sixrm
   \textfont\@ne\eleveni  \scriptfont\@ne\eighti
       \scriptscriptfont\@ne\sixi
   \textfont\tw@\elevensy \scriptfont\tw@\eightsy
       \scriptscriptfont\tw@\sixsy
   \textfont\thr@@\elevenex \scriptfont\thr@@\eightex
        \scriptscriptfont\thr@@\eightex
   \textfont\itfam\elevenit \scriptfont\itfam\eightit
        \scriptscriptfont\itfam\eightit
   \textfont\bffam\elevenbf \scriptfont\bffam\eightbf
        \scriptscriptfont\bffam\sixbf
   \setbox\strutbox\hbox{\vrule height8.5\p@ depth3.5\p@ width\z@}%
   \setbox\strutbox@\hbox{\lower.5\normallineskiplimit\vbox{%
        \kern-\normallineskiplimit\copy\strutbox}}%
   \setbox\z@\vbox{\hbox{$($}\kern\z@}\bigsize@1.2\ht\z@
  \fi
  \normalbaselines\rm\dotsspace@1.5mu\ex@.2326ex\jot3\ex@
  \the\elevenpoint@}

\addto\tenpoint{\normalbaselineskip12\p@
 \abovedisplayskip6\p@ plus6\p@ minus0\p@
 \belowdisplayskip6\p@ plus6\p@ minus0\p@
 \abovedisplayshortskip0\p@ plus3\p@ minus0\p@
 \belowdisplayshortskip2\p@ plus3\p@ minus0\p@
 \ifsyntax@
 \else
  \setbox\strutbox\hbox{\vrule height9\p@ depth4\p@ width\z@}%
  \setbox\strutbox@\hbox{\vrule height8\p@ depth3\p@ width\z@}%
 \fi
 \normalbaselines\rm}

\newtoks\sevenpoint@
\def\sevenpoint{\normalbaselineskip9\p@
 \textonlyfont@\rm\sevenrm \textonlyfont@\it\sevenit
 \textonlyfont@\sl\sevensl \textonlyfont@\bf\sevenbf
 \textonlyfont@\smc\sevensmc \textonlyfont@\tt\seventt
  \textfont\z@\sevenrm \scriptfont\z@\sixrm
       \scriptscriptfont\z@\fiverm
  \textfont\@ne\seveni \scriptfont\@ne\sixi
       \scriptscriptfont\@ne\fivei
  \textfont\tw@\sevensy \scriptfont\tw@\sixsy
       \scriptscriptfont\tw@\fivesy
  \textfont\thr@@\sevenex \scriptfont\thr@@\sevenex
   \scriptscriptfont\thr@@\sevenex
  \textfont\itfam\sevenit \scriptfont\itfam\sevenit
   \scriptscriptfont\itfam\sevenit
  \textfont\bffam\sevenbf \scriptfont\bffam\sixbf
   \scriptscriptfont\bffam\fivebf
 \setbox\strutbox\hbox{\vrule height7\p@ depth3\p@ width\z@}%
 \setbox\strutbox@\hbox{\raise.5\normallineskiplimit\vbox{%
   \kern-\normallineskiplimit\copy\strutbox}}%
 \setbox\z@\vbox{\hbox{$($}\kern\z@}\bigsize@1.2\ht\z@
 \normalbaselines\sevenrm\dotsspace@1.5mu\ex@.2326ex\jot3\ex@
 \the\sevenpoint@}

\newskip\abovespecheadskip \abovespecheadskip18\p@ plus 4\p@ minus 2\p@
\newdimen\belowspecheadskip      \belowspecheadskip6\p@

\outer\def\specialhead{%
  \add@missing\endroster \add@missing\enddefinition
  \add@missing\enddemo \add@missing\endexample
  \add@missing\endproclaim
  \penaltyandskip@{-200}\abovespecheadskip
  \begingroup\interlinepenalty\@M\rightskip\z@ plus\hsize
  \let\\\linebreak
  \specialheadfont@\raggedcenter@\noindent}
\def\endspecialhead{\endgraf\endgroup\nobreak\vskip\belowspecheadskip}

\let\varindent@\indent

\let\subsubhead\relax
\outer\def\subsubhead{%
  \add@missing\endroster \add@missing\enddefinition
  \add@missing\enddemo
  \add@missing\endexample \add@missing\endproclaim
  \let\savedef@\subsubhead \let\subsubhead\relax
  \def\subsubhead##1\endsubsubhead{\restoredef@\subsubhead
      {\def\usualspace{\/{\it\enspace}}%
    \it##1\unskip\frills@{.\enspace}}}%
  \nofrillscheck\subsubhead}

\sl).

\newskip\abstractindent 	\abstractindent=3pc
\long\def\block #1\endblock{\vskip 6pt
	{\leftskip=\abstractindent \rightskip=\abstractindent
	\noindent #1\endgraf}\vskip 6pt}

\long\def\ext #1\endext{\removelastskip\block #1\endblock}

\outer\def\xca{\let\savedef@\xca \let\xca\relax
  \add@missing\endproclaim \add@missing\endroster
  \add@missing\endxca \envir@stack\endxca
   \def\xca##1{\restoredef@\xca
   \penaltyandskip@{-100}\medskipamount
        \bgroup{\def\usualspace{{\xcaheadfont@\enspace}}%
        \varindent@\xcaheadfont@\ignorespaces##1\unskip
        \frills@{.\xcaheadfont@\enspace}}%
        \ignorespaces}%
  \nofrillscheck\xca}
\def\endxca{\egroup\revert@envir\endxca
  \medskip}

\let\proclaim\relax
\outer\def\proclaim{%
  \let\savedef@\proclaim \let\proclaim\relax
  \add@missing\endroster \add@missing\enddefinition
\add@missing\endproclaim
  \envir@stack\endproclaim
 \def\proclaim##1{\restoredef@\proclaim
   \penaltyandskip@{-100}\medskipamount\varindent@
   \def\usualspace{{\proclaimheadfont@\enspace}}\proclaimheadfont@
   \ignorespaces##1\unskip\frills@{.\enspace}%
  \it\ignorespaces}
\it.
	\nofrillscheck\proclaim}

\def\remark{\let\savedef@\remark \let\remark\relax
  \add@missing\endroster \add@missing\endproclaim
  \envir@stack\endremark
  \def\remark##1{\restoredef@\remark
  {\def\usualspace{{\remarkheadfont@\enspace}}%
   \penaltyandskip@{-100}\medskipamount\varindent@
  \remarkheadfont@\ignorespaces##1\unskip\frills@{.\enspace}}\rm
  \ignorespaces}\nofrillscheck\remark}

\def\remarkheadfont@{\smc}

\def\endremark{\revert@envir\endremark\medskip}

\def\qed{\ifhmode\unskip\nobreak\fi\hfill
  \ifmmode\square\else$\m@th\square$\fi}

\def\definition{\let\savedef@\definition \let\definition\relax
  \add@missing\endproclaim \add@missing\endroster
  \add@missing\enddefinition \envir@stack\enddefinition
   \def\definition##1{\restoredef@\definition
   \penaltyandskip@{-100}\medskipamount\varindent@
        {\def\usualspace{{\proclaimheadfont@\enspace}}%
        \proclaimheadfont@\ignorespaces##1\unskip
        \frills@{.\proclaimheadfont@\enspace}}%
        \rm \ignorespaces}%
  \nofrillscheck\definition}


\newdimen\rosteritemsep
\rosteritemsep=.5pc

\newdimen\rosteritemitemitemwd
\newdimen\rosteritemitemwd

\newbox\setwdbox
\setbox\setwdbox\hbox{0.}\rosteritemwd=\wd\setwdbox
  \setbox\setwdbox\hbox{0.\hskip.5pc(c)}\rosteritemitemwd=\wd\setwdbox
  \setbox\setwdbox\hbox{0.\hskip.5pc(c)\hskip.5pc(iii)}%
	\rosteritemitemitemwd=\wd\setwdbox

\def\roster{%
  \envir@stack\endroster
 \edef\leftskip@{\leftskip\the\leftskip}%
 \relaxnext@
 \rostercount@\z@
 \def\item{\FN@\rosteritem@}%
 \def\itemitem{\FN@\rosteritemitem@}%
 \def\itemitemitem{\FN@\rosteritemitemitem@}%
 \DN@{\ifx\next\runinitem\let\next@\nextii@\else
  \let\next@\nextiii@\fi\next@}%
 \DNii@\runinitem
  {\unskip
   \DN@{\ifx\next[\let\next@\nextii@\else
    \ifx\next"\let\next@\nextiii@\else\let\next@\nextiv@\fi\fi\next@}%
   \DNii@[####1]{\rostercount@####1\relax
    \therosteritem{\number\rostercount@}~\ignorespaces}%
   \def\nextiii@"####1"{{\rm####1}~\ignorespaces}%
   \def\nextiv@{\therosteritem1\rostercount@\@ne~}%
   \par@\firstitem@false
   \FN@\next@}
 \def\nextiii@{\par\par@
  \penalty\@m\vskip-\parskip
  \firstitem@true}%
 \FN@\next@}

\def\rosteritem@{\iffirstitem@\firstitem@false
  \else\par\vskip-\parskip\fi
 \leftskip\rosteritemwd \advance\leftskip\normalparindent
 \advance\leftskip.5pc \noindent
 \DNii@[##1]{\rostercount@##1\relax\itembox@}%
 \def\nextiii@"##1"{\def\therosteritem@{\rm##1}\itembox@}%
 \def\nextiv@{\advance\rostercount@\@ne\itembox@}%
 \def\therosteritem@{\therosteritem{\number\rostercount@}}%
 \ifx\next[\let\next@\nextii@\else\ifx\next"\let\next@\nextiii@\else
  \let\next@\nextiv@\fi\fi\next@}

\def\itembox@{\llap{\hbox to\rosteritemwd{\hss
  \kern\z@ 
  \therosteritem@}\hskip.5pc}\ignorespaces}

\def\therosteritem#1{\rom{\ignorespaces#1.\unskip}}

\def\rosteritemitem@{\iffirstitem@\firstitem@false
  \else\par\vskip-\parskip\fi
 \leftskip\rosteritemitemwd \advance\leftskip\normalparindent
 \advance\leftskip.5pc \noindent
 \DNii@[##1]{\rostercount@##1\relax\itemitembox@}%
 \def\nextiii@"##1"{\def\therosteritemitem@{\rm##1}\itemitembox@}%
 \def\nextiv@{\advance\rostercount@\@ne\itemitembox@}%
 \def\therosteritemitem@{\therosteritemitem{\number\rostercount@}}%
 \ifx\next[\let\next@\nextii@\else\ifx\next"\let\next@\nextiii@\else
  \let\next@\nextiv@\fi\fi\next@}

\def\itemitembox@{\llap{\hbox to\rosteritemitemwd{\hss
  \kern\z@ 
  \therosteritemitem@}\hskip.5pc}\ignorespaces}

\def\therosteritemitem#1{\rom{(\ignorespaces#1\unskip)}}

\def\rosteritemitemitem@{\iffirstitem@\firstitem@false
  \else\par\vskip-\parskip\fi
 \leftskip\rosteritemitemitemwd \advance\leftskip\normalparindent
 \advance\leftskip.5pc \noindent
 \DNii@[##1]{\rostercount@##1\relax\itemitemitembox@}%
 \def\nextiii@"##1"{\def\therosteritemitemitem@{\rm##1}\itemitemitembox@}%
 \def\nextiv@{\advance\rostercount@\@ne\itemitemitembox@}%

\def\therosteritemitemitem@{\therosteritemitemitem{\number\rostercount@}}%
 \ifx\next[\let\next@\nextii@\else\ifx\next"\let\next@\nextiii@\else
  \let\next@\nextiv@\fi\fi\next@}

\def\itemitemitembox@{\llap{\hbox to\rosteritemitemitemwd{\hss
  \kern\z@ 
  \therosteritemitemitem@}\hskip.5pc}\ignorespaces}

\def\therosteritemitemitem#1{\rom{(\ignorespaces#1\unskip)}}

\def\endroster{\relaxnext@
 \revert@envir\endroster 
 \par\leftskip@
 \penalty-50
 \DN@{\ifx\next\Runinitem\let\next@\relax
  \else\nextRunin@false\let\item\plainitem@
   \ifx\next\par
    \DN@\par{\everypar\expandafter{\the\everypartoks@}}%
   \else
    \DN@{\noindent\everypar\expandafter{\the\everypartoks@}}%
  \fi\fi\next@}%
 \FN@\next@}

\def\address#1\endaddress{\global\advance\addresscount@\@ne
  \expandafter\gdef\csname address\number\addresscount@\endcsname
  {\vskip12\p@ minus6\p@\indent\addressfont@{\smc\ignorespaces#1}\par}}


\def\email{\let\savedef@\email
  \def\email##1\endemail{\let\email\savedef@
  \toks@{\def\usualspace{{\it\enspace}}\endgraf\indent\addressfont@}%
  \toks@@{{\tt ##1}\par}%
  \expandafter\xdef\csname email\number\addresscount@\endcsname
  {\the\toks@\frills@{{\noexpand\it E-mail address\noexpand\/}:%
     \noexpand\enspace}\the\toks@@}}%
  \nofrillscheck\email}

\def\curraddr{\let\savedef@\curraddr
  \def\curraddr##1\endcurraddr{\let\curraddr\savedef@
  \toks@\expandafter\expandafter\expandafter{%
       \csname address\number\addresscount@\endcsname}%
  \toks@@{##1}%
  \expandafter\xdef\csname address\number\addresscount@\endcsname
  {\the\toks@\endgraf\noexpand\nobreak
    \indent\noexpand\addressfont@{\noexpand\rm
    \frills@{{\noexpand\it Current address\noexpand\/}:\space}%
    \def\noexpand\usualspace{\space}\the\toks@@\unskip}}}%
  \nofrillscheck\curraddr}

\def\topcaption#1#2\endcaption{%
  {\dimen@\hsize \advance\dimen@-\captionwidth@
   \rm\raggedcenter@ \advance\leftskip.5\dimen@ \rightskip\leftskip
  {\captionfont@#1}%
  \if\notempty{#2}.\enspace\ignorespaces#2\fi
  \endgraf}\nobreak\bigskip}
\def\botcaption#1#2\endcaption{%
  \nobreak\bigskip
  \setboxz@h{\captionfont@#1\if\notempty{#2}.\enspace\rm#2\fi}%
  {\dimen@\hsize \advance\dimen@-\captionwidth@
   \leftskip.5\dimen@ \rightskip\leftskip
   \noindent \ifdim\wdz@>\captionwidth@
   \else\hfil\fi
  {\captionfont@#1}%
  \if\notempty{#2}.\enspace\rm#2\fi\endgraf}}

\def\rom#1{{\rm #1}}

matters.

\def\bysame{\by\hbox to2pc{\hrulefill}\thinspace\kern\z@}

\def\refstyle#1{\uppercase{%
  \if#1A\relax \def\keyformat##1{[##1]\enspace\hfil}%
  \else\if#1B\relax
    \def\keyformat##1{\aftergroup\kern
              \aftergroup-\aftergroup\refindentwd}%
    \refindentwd2pc
 \else\if#1C\relax
    \refindentwd\parindent
   \def\keyformat##1{\hfil##1.\enspace}%
 \fi\fi\fi}
}

\let\Refs\relax
\outer\def\Refs{\add@missing\endroster \add@missing\endproclaim
 \let\savedef@\Refs \let\Refs\relax 
 \def\Refs##1{\restoredef@\Refs
   \if\notempty{##1}\penaltyandskip@{-200}\aboveheadskip
     \begingroup \raggedcenter@\refheadfont@
       \ignorespaces##1\endgraf\endgroup
     \penaltyandskip@\@M\belowheadskip
   \fi
   \begingroup\def\envir@end{\endRefs}\refsfont@\sfcode`\.\@m
   }%
 \nofrillscheck{\csname Refs\expandafter\endcsname
  \frills@{{References}}}}

\catcode`\@=11


\pageheight{586pt}	

\def\keyboarder#1{}%


\def\addressfont@{\tenpoint}

\skip\footins=\bigskipamount 

\font@\titlebf=cmbx10 scaled \magstep2	
\font@\titlei=cmmi10 scaled \magstep2
\font@\titlesy=cmsy10 scaled \magstep2
\def\titlefont{\normalbaselineskip18\p@
 \textonlyfont@\bf\titlebf
 \ifsyntax@\else
  \textfont\z@\titlebf  \scriptfont\z@\tenbf  \scriptscriptfont\z@\sevenbf
  \textfont\@ne\titlei  \scriptfont\@ne\teni  \scriptscriptfont\@ne\seveni
  \textfont\tw@\titlesy \scriptfont\tw@\tensy
\scriptscriptfont\tw@\sevensy
  \textfont\thr@@\tenex \scriptfont\thr@@\tenex
\scriptscriptfont\thr@@\tenex
 \fi
 \normalbaselines\titlebf}

\newif\ifpart	\partfalse

\font\partnofont=cmbx10 scaled \magstep3
\font\partfont=cmbx10 scaled \magstep4

\def\part#1\\#2\endpart{\global\parttrue
  \gdef\thepart{\def\\{\hfil\break}%
  {\raggedcenter@
	{\partnofont Part #1\endgraf}%
	\vskip 20pt
	\partfont\baselineskip=25pt
   #2\endgraf}}}

\def\title#1\endtitle{%
  \global\setbox\titlebox@\vtop{\titlefont
   \raggedcenter@
   #1\endgraf}%
 \ifmonograph@ \edef\next{\the\leftheadtoks}%
    \ifx\next\empty@
    \leftheadtext{#1}\fi
 \fi
 \edef\next{\the\rightheadtoks}\ifx\next\empty@ \rightheadtext{#1}\fi
 }


\def\romanchapternum{\gdef\chapterno@{\uppercase\expandafter{\romannumeral
    \chaptercount@}}}

\def\chapterno@{\the\chaptercount@}

\def\chapter{\let\savedef@\chapter
  \def\chapter##1{\let\chapter\savedef@
  \leavevmode\hskip-\leftskip
   \rlap{\vbox to\z@{\vss\centerline{\tenpoint
   \frills@{CHAPTER\space\afterassignment\chapterno@
       \global\chaptercount@=}%
   ##1\unskip}\baselineskip36pt\null}}\hskip\leftskip}%
 \nofrillscheck\chapter}


\def\author{\let\savedef@\author
  \def\author##1\endauthor{\global\setbox\authorbox@
 \vbox{\tenpoint\raggedcenter@
  \frills@{\ignorespaces##1}\endgraf}
 \edef\next{\the\leftheadtoks}%
 \ifx\next\empty@\expandafter\uppercase{\leftheadtext{##1}}\fi}
\nofrillscheck\author}


\def\abstract{\let\savedef@\abstract
 \def\abstract{\let\abstract\savedef@
  \setbox\abstractbox@\vbox\bgroup\noindent$$\vbox\bgroup\indenti=0pt
  \def\envir@end{\endabstract}\advance\hsize-2\indenti
  \def\usualspace{\enspace}\tenpoint \noindent
  \frills@{{\abstractfont@ Abstract.\enspace}}}%
 \nofrillscheck\abstract}

\def\thanks#1\endthanks{%
  \gdef\thethanks@{\tenpoint\raggedcenter@#1\endgraf}}

\def\keywords{\let\savedef@\keywords
  \def\keywords##1\endkeywords{\let\keywords\savedef@
  \toks@{\def\usualspace{{\it\enspace}}\raggedcenter@\tenpoint}%
  \toks@@{##1\unskip.}%
  \edef\thekeywords@{\the\toks@\frills@{{\noexpand\it
    Key words and phrases.\noexpand\enspace}}\the\toks@@}}%
 \nofrillscheck\keywords}

\def\subjclass{\let\savedef@\subjclass
 \def\subjclass##1\endsubjclass{\let\subjclass\savedef@
   \toks@{\def\usualspace{{\rm\enspace}}\tenpoint\raggedcenter@}%
   \toks@@{##1\unskip.}%
   \edef\thesubjclass@{\the\toks@
     \frills@{{\noexpand\rm1991 {\noexpand\it Mathematics Subject
       Classification}.\noexpand\enspace}}%
     \the\toks@@}}%
  \nofrillscheck\subjclass}

\newskip\xcskip   \newskip\afterxcskip
\xcskip=10pt plus2pt minus0pt
\afterxcskip=0pt
\long\def\xcb#1{\par\ifnum\lastskip<\xcskip
  \removelastskip\penalty-100\vskip\xcskip\fi
  \noindent{\bf#1}%
  \nobreak\bgroup
  \xcbrosterdefs}

\def\xcbrosterdefs{\normalparindent=0pt
\setbox\setwdbox\hbox{0.}\rosteritemwd=\wd\setwdbox\relax

\setbox\setwdbox\hbox{0.\hskip.5pc(c)}\rosteritemitemwd=\wd\setwdbox\relax
  \setbox\setwdbox\hbox{0.\hskip.5pc(c)\hskip.5pc(iii)}%
	\rosteritemitemitemwd=\wd\setwdbox\relax
}

\def\endxcb{\par\egroup}

\newif\ifplain@  \plain@false

\def\output@{\shipout\vbox{%
 \ifplain@
    \pagebody
  \else
 	\iffirstpage@
		\global\firstpage@false
		 \pagebody \logo@ \makefootline%
 	\else
  		\ifrunheads@
     			\makeheadline \pagebody
  		\else
      			\pagebody \makefootline
  		\fi
 	\fi
  \fi}%
 \advancepageno \ifnum\outputpenalty>-\@MM\else\dosupereject\fi}

\outer\def\endtopmatter{\add@missing\endabstract
 \edef\next{\the\leftheadtoks}\ifx\next\empty@
  \expandafter\leftheadtext\expandafter{\the\rightheadtoks}\fi
   \ifpart
	  \global\plain@true
	  \null\vskip9.5pc
	  \thepart
	  \vfill\eject
	  \null\vfill\eject
   \fi
  \global\plain@false
  \global\firstpage@true
  \begingroup 
	  \topskip100pt
  \box\titlebox@
  \endgroup
  \ifvoid\authorbox@\else \vskip2.5pcplus1pc\unvbox\authorbox@\fi
  \ifnum\addresscount@>\z@
	\vfill
	Author addresses:
 \count@\z@ \loop\ifnum\count@<\addresscount@\advance\count@\@ne
 \csname address\number\count@\endcsname
 \csname email\number\count@\endcsname
 \repeat
	\vfill\eject\fi
  \ifvoid\affilbox@\else
	\vskip1pcplus\unvbox\affilbox@\fi
   \ifx\thesubjclass@\empty@\else \vfil\thesubjclass@\fi
   \ifx\thekeywords@\empty@\else \vfil\thekeywords@\fi
   \ifx\thethanks@\empty@\else \vfil\thethanks@\fi
  \ifvoid\abstractbox@\else \vfil\unvbox\abstractbox@\vfil\eject \fi
  \ifvoid\tocbox@\else\vskip1.5pcplus.5pc\unvbox\tocbox@\fi
  \prepaper
  \vskip22pt\relax
}

\def\raggedleft@{\leftskip\z@ plus.4\hsize \rightskip\z@
 \parfillskip\z@ \parindent\z@ \spaceskip.3333em \xspaceskip.5em
 \pretolerance9999\tolerance9999 \exhyphenpenalty\@M
 \hyphenpenalty\@M \let\\\linebreak}

\def\aufm #1\endaufm{\vskip12pt{\raggedleft@ #1\endgraf}}

\widestnumber\key{M} 


\begingroup
\let\head\relax \let\specialhead\relax \let\subhead\relax
\let\subsubhead\relax \let\title\relax \let\chapter\relax
\newbox\tocchapbox@

\gdef\newtocdefs{%
  \def\ptitle##1\endptitle
       {\penalty\z@ \vskip8\p@
        \hangindent\wd\tocheadbox@\noindent{\bf ##1}\endgraf}%
  \def\title##1\endtitle
       {\penalty\z@ \vskip8\p@
        \hangindent\wd\tocchapbox@\noindent{##1}\endgraf}%
  \def\chapter##1{\par
        Chapter ##1.\unskip\enspace}%
  \def\part##1{\par
        {\bf Part ##1.}\unskip\enspace}%
  \def\specialhead##1 ##2\endspecialhead{\par
    \begingroup \hangindent.5em \noindent
    \if\notempty{##1}%
      \leftskip\wd\tocheadbox@
      \llap{\hbox to\wd\tocheadbox@{\hfil##1}}\enspace
    \else
      \leftskip\parindent
    \fi
    ##2\endgraf
    \endgroup}%
  \def\head##1 ##2\endhead{\par
    \begingroup \hangindent.5em \noindent
    \if\notempty{##1}%
      \leftskip\wd\tocheadbox@
      \llap{\hbox to\wd\tocheadbox@{\hfil##1}}\enspace
    \else
      \leftskip\parindent
    \fi
    ##2\endgraf
    \endgroup}%
  \def\subhead##1 ##2\endsubhead{\par
    \begingroup \leftskip4.5pc
    \noindent\llap{##1\enspace}##2\endgraf \endgroup}%
  \def\subsubhead##1 ##2\endsubsubhead{\par
    \begingroup \leftskip5pc
    \noindent\hbox{\ignorespaces##1 }##2\endgraf \endgroup}%
}%
\gdef\toc@#1{\relaxnext@
 \DN@{\ifx\next\nofrills\DN@\nofrills{\nextii@}%
      \else\DN@{\nextii@{{#1}}}\fi
      \next@}%
 \DNii@##1{%
\ifmonograph@\bgroup\else\setbox\tocbox@\vbox\bgroup
   \centerline{\headfont@\ignorespaces##1\unskip}\nobreak
   \vskip\belowheadskip \fi
   \def\page####1%
       {\unskip\penalty\z@\null\hfil
        \rlap{\hbox to\pagenumwd{\quad\hfil\rm ####1}}%
   \global\setbox\tocchapbox@\hbox{Chapter 1.\enspace}%
   \global\setbox\tocheadbox@\hbox{\hskip18pt \S0.0.}%
         \hfilneg\penalty\@M}%
   \leftskip\z@ \rightskip\leftskip
   \setboxz@h{\bf\quad000}\pagenumwd\wd\z@
   \advance\rightskip\pagenumwd
   \newtocdefs
 }%
 \FN@\next@}%
\endgroup

\def\logo@{}

\Monograph

\catcode`\@=13

\refstyle{A} 
\define\dfn#1{{\it #1\/}}
\define\dfnb#1{{\it #1}}
\loadeusm 	               
\define\script{\eusm}

\catcode`\@=11
\def\hookrightarrowfill{$\m@th\mathord\lhook\mkern-3mu%
  \mathord-\mkern-6mu%
  \cleaders\hbox{$\mkern-2mu\mathord-\mkern-2mu$}\hfill
  \mkern-6mu\mathord\rightarrow$}
\def\hookleftarrowfill{$\m@th\mathord\leftarrow\mkern-6mu%
  \cleaders\hbox{$\mkern-2mu\mathord-\mkern-2mu$}\hfill
  \mkern-6mu\mathord-\mkern-3mu\mathord\rhook$}
\atdef@ C#1C#2C{\ampersand@\setbox\z@\hbox{$\ssize
 \;{#1}\;\;$}\setbox\@ne\hbox{$\ssize\;{#2}\;\;$}\setbox\tw@
 \hbox{$#2$}\ifCD@
 \global\bigaw@\minCDaw@\else\global\bigaw@\minaw@\fi
 \ifdim\wd\z@>\bigaw@\global\bigaw@\wd\z@\fi
 \ifdim\wd\@ne>\bigaw@\global\bigaw@\wd\@ne\fi
 \ifCD@\hskip.5em\fi
 \ifdim\wd\tw@>\z@
 \mathrel{\mathop{\hbox to\bigaw@{\hookrightarrowfill}}%
 \limits^{#1}_{#2}}\else
 \mathrel{\mathop{\hbox to\bigaw@{\hookrightarrowfill}}%
 \limits^{#1}}\fi
 \ifCD@\hskip.5em\fi\ampersand@}
\atdef@ D#1D#2D{\ampersand@\setbox\z@\hbox{$\ssize
 \;\;{#1}\;$}\setbox\@ne\hbox{$\ssize\;\;{#2}\;$}\setbox\tw@
 \hbox{$#2$}\ifCD@
 \global\bigaw@\minCDaw@\else\global\bigaw@\minaw@\fi
 \ifdim\wd\z@>\bigaw@\global\bigaw@\wd\z@\fi
 \ifdim\wd\@ne>\bigaw@\global\bigaw@\wd\@ne\fi
 \ifCD@\hskip.5em\fi
 \ifdim\wd\tw@>\z@
 \mathrel{\mathop{\hbox to\bigaw@{\hookleftarrowfill}}%
 \limits^{#1}_{#2}}\else
 \mathrel{\mathop{\hbox to\bigaw@{\hookleftarrowfill}}%
 \limits^{#1}}\fi
 \ifCD@\hskip.5em\fi\ampersand@}
\catcode`\@=13
\define\bs{\backslash}

\define\dual{\spcheck}
\define\hra{\hookrightarrow}

\predefine\isom{\cong}
\redefine\cong{\equiv}
\predefine\imaginary{\Im}
\redefine\Im{\operatorname{Im}}

\define\lra{\longrightarrow}
\define\ol{\overline}

\define\simto{\overset\sim\to\longrightarrow}
\define\tensor{\otimes}


\define\varemptyset{\varnothing}
\define\Qed{\hbox to 0.5em{ }\nobreak\hfill\hbox{$\square$}}

\define\Ad{\operatorname{Ad}}
\define\ad{\operatorname{ad}}
\define\angled#1{\langle #1 \rangle}
\define\Aut{\operatorname{Aut}}
\define\CC{{\Bbb C}}
\define\Div{{\script D}}
\define\End{\operatorname{End}}
\define\Endo{\operatorname{End}^0}
\define\EndoA{\Endo(A)}
\define\GG{{\Bbb G}}
\define\Gal{\operatorname{Gal}}
\define\GHC{\operatorname{GHC}}
\define\GL{\operatorname{GL}}
\define\Gm{{\GG_{\text{m}}}{\vphantom{\GG}}}
\define\GSp{\operatorname{GSp}}
\define\Hdg{{\script H}}
\define\Hg{\operatorname{Hg}}
\define\hg{\frak{hg}}
\define\Hom{\operatorname{Hom}}
\define\Id{\operatorname{Id}}
\define\id{\operatorname{id}}
\define\ind{\operatorname{ind}}
\define\Ker{\operatorname{Ker}}
\define\Lie{\operatorname{Lie}}
\define\Lf{\operatorname{Lf}}
\define\MT{\operatorname{MT}}

\define\mt{\frak{mt}}

\define\NN{{\Bbb N}}
\define\scO{{\script O}}

\redefine\phi{\varphi}
\define\phibar{\ol{\phi}}
\define\QQ{{\Bbb Q}}
\define\Qbar{{\ol{\QQ}}}
\define\RR{{\Bbb R}}
\define\rank{\operatorname{rank}}
\define\rdim{\operatorname{rdim}}
\define\Res{\operatorname{Res}}
\define\SS{{\Bbb S}}
\define\Sbar{\ol{S}}

\define\SL{\operatorname{SL}}

\define\SO{\operatorname{SO}}
\predefine\Sup{\Sp}
\redefine\Sp{\operatorname{Sp}}
\define\spec{\operatorname{spec}}
\define\SU{\operatorname{SU}}
\define\tr{\,\mathstrut^t}
\define\U{\operatorname{U}}
\define\V{{\script V}}
\define\Vr{V_\RR}
\define\twedge{{\tsize \bigwedge}}
\define\Wr{W_\RR}
\define\Weil{{\script W}}
\define\X{\operatorname{X}}
\define\ZZ{{\Bbb Z}}

\pageno=-1
\topmatter
\leftheadtext{appendix b. The Hodge conjecture for abelian varieties}
\endtopmatter
\pageno=306
\document
{\bf Appendix B. A Survey of the Hodge Conjecture for Abelian Varieties}
\smallskip
\indent\hskip1.0in {\bf by B. Brent Gordon}

\head Introduction \endhead
 The goal of this appendix is to review what is known about the Hodge
conjecture for abelian varieties, with an emphasis on how Mumford-Tate
groups have been applied to this problem.  In addition to the book in
which this appears, other survey or general articles that precede this one
are Hodge's original paper \cite{B.53}, Grothendieck's modification of
Hodge's general conjecture \cite{B.43}, Shioda's excellent survey article
\cite{B.117}, Steenbrink's comments on the general Hodge conjecture
\cite{B.120}, and van Geemen's pleasing introduction to the Hodge
conjecture for abelian varieties \cite{B.35}.  Naturally there is some
overlap between this appendix and van Geemen's article, but since his
emphasis is on abelian varieties of Weil type, we hope that this appendix
will be a useful complement.

\medpagebreak

  Since the language of linear algebraic groups and their Lie algebras,
which cannot be avoided in any discussion of Mumford-Tate groups and their
application to the Hodge conjecture for abelian varieties, may not be
familiar to students of complex algebraic geometry and Hodge theory, we
begin by recalling the definitions and facts we need and giving some
examples.   Towards the end of section one we also recall some basic facts
about abelian varieties, including the Albert classification of their
endomorphism algebras (Theorem~1.12.2), and give some example of abelian
varieties to which we refer later.  Most readers will find it more
profitable to begin with section two, where we discuss the definitions and
some general structural properties of the Hodge, Mumford-Tate and
Lefschetz groups associated to an abelian variety, or section three,
and refer to the first section as needed.  

  Starting with section three we have tried to be comprehensive in
summarizing the known results and indicating the main ideas involved in
their proofs, while at the same time selecting some cross-section of
proofs to discuss in more detail.  In section three we follow Murty's
exposition \cite{B.84} of the Hodge $(p,p)$ conjecture for arbitrary
products of elliptic curves.  In section four we summarize Shioda's
results on abelian varieties of Fermat type
\cite{B.116}, but only briefly consider the issues and results related to
abelian varieties of Weil type, since \cite{B.35} treats this topic well. 
In section five we discuss the work of Moonen and Zarhin on
four-dimensional abelian varieties \cite{B.75}, for this provides a nice
illustration of the interplay between the endomorphism algebra and the
Mumford-Tate group of an abelian variety. In section six we look at the
work of Tankeev and Ribet on the Hodge conjecture for simple abelian
varieties that satisfy some conditions on their dimension or endomorphism
algebra \cite{B.124} \cite{B.125} \cite{B.126} \cite{B.93} \cite{B.94};
for example, the Hodge conjecture is true for simple abelian varieties of
prime dimension.  Here we look more closely at Ribet's approach, where he
introduced and used the Lefschetz group of an abelian variety.  Then the
results of Murty and Hazama discussed in section seven build on and go
beyond Ribet's methods to treat abelian varieties not assumed to be
simple, but still assumed to satisfy some conditions on their dimensions
or endomorphism algebras \cite{B.81} \cite{B.82} \cite{B.83} \cite{B.46}
\cite{B.47} \cite{B.49}.

  In section eight we shift directions slightly, for here we have
collected together examples of exceptional Hodge cycles, i.e., Hodge
cycles not accounted for by linear combinations of intersections of
divisor classes, and in this section, not known to be algebraic.  Largely
missing from section eight, but considered in section nine, are the
particular problems posed by abelian varieties of complex multiplication
type.  Dodson \cite{B.28} \cite{B.29} \cite{B.30} and others have
constructed numerous examples of such abelian varieties that support
exceptional Hodge cycles.

 In section ten we examine what is known about the general Hodge
conjecture for abelian varieties.  The majority of the work on this
problem is either a very geometric treatment of special abelian varieties
in low dimension, for example \cite{B.12} or \cite{B.104}, or requires
special assumptions about the endomorphism algebra, dimension or Hodge
group, as in \cite{B.127}, \cite{B.128}, \cite{B.50} or \cite{B.4}.

 In the final section eleven we briefly mention three alternative
approaches to proving the (usual) Hodge conjecture for arbitrary abelian
varieties:  First, a method involving the Weil intermediate Jacobian
\cite{B.98}; then that the Tate conjecture for abelian varieties implies
the Hodge conjecture for abelian varieties \cite{B.88} \cite{B.87} and
\cite{B.27}; and thirdly, that the Hodge conjecture for abelian varieties
would follow from knowing Grothendieck's invariant cycles conjecture
(\cite{B.42}) for certain general families of abelian varieties, and
moreover, that for these families, the invariant cycles conjecture would
follow from the $L_2$-cohomology analogue of Grothendieck's standard
conjecture~(A) that the Hodge $*$-operator is algebraic (\cite{B.44})
\cite{B.4}.  The present state of our knowledge about the Hodge
conjecture for abelian varieties is such that any or none of these
approaches might ultimately work, or a counterexample might exist.

  Preceding the bibliography is a rough chronological table
of the work that directly address some aspect of the Hodge
conjecture for abelian varieties.  I have tried to make sure that this
table and this appendix as a whole mention all the relevant references
through the end of 1996; if I have omitted something or otherwise not done
it justice, that was quite unintentional.


\head 1. Abelian varieties and linear algebraic groups
\endhead
  The purpose of this section is to establish the language we use
throughout the rest of this appendix to discuss abelian varieties and
certain linear algebraic groups and Lie algebras associated with them.
Although abelian varieties and linear algebraic groups are both algebraic
groups, the issues surrounding them tend to be of a very different
nature.  It turns out to be most convenient to begin by recalling some of
the definitions and basic properties of linear algebraic groups and their
Lie algebras, and introducing some of the examples of these to which we
will later refer, and then in the second half of the section review some
of the definitions and basic properties of abelian varieties, and
introduce some of the examples we will investigate later.

\subhead 1.0. Notational conventions
\endsubhead
\nopagebreak
\remark{1.0.1. Field of definition}
  Let $F$ be a field and $V$ and algebraic variety.  Then we will write
$V_F$ to signify or emphasize that $V$ is defined over $F$.  When $V$ is
an algebraic variety defined over $F$ and $K$ is a field containing $F$,
then $V_K = V_F \times_{\spec F} \spec K$ is the base change to $K$,
i.e., $V$ as a variety defined over $K$.  We will generally try to
distinguish the abstract variety $V_F$ defined over $F$ from its concrete
set of $F$-points $V(F)$, and then $V(K) = V_F(K)$ is the set of
$K$-points.
\endremark

\definition{1.0.2. Definition}
  Suppose $K$ is a separable algebraic extension of $F$ of finite
degree~$d$, and $V$ is an algebraic variety defined over the larger
field~$K$.  Let $\{ \sigma_1, \ldots, \sigma_d\}$ be the set of distinct
embeddings of $K$ into the algebraic closure $F^{\text{alg}}$ of $F$.
Then the \dfn{restriction of scalars functor} $\Res_{K/F}$ from
varieties over $K$ to varieties over $F$ is defined as follows:  First
let $V_{\sigma_i} = V_K \times_{\spec K,\,\sigma_i} F^{\text{alg}}$.
Then for any variety $W$ defined over $F$ and a morphism $\phi:W\to V$
defined over $K$ there are morphisms $\phi_{\sigma_i} : W\to
V_{\sigma_i}$.  Then if
$$
 (\phi_{\sigma_1},\ldots , \phi_{\sigma_d}) : W \to V_{\sigma_1} \times
\cdots \times V_{\sigma_d}
$$
 is an isomorphism, then \dfn{$W = \Res_{K/F} V$ is the variety obtained
from $V$ by restriction of the field of definition from $K$ to $F$.}  Its
uniqueness is a consequence of the universal property that whenever $X$ is
any variety defined over $F$ and $\psi: X \to V$ is a morphism defined
over $K$, then there exists a unique $\Psi: X\to W$ defined over $F$ such
that $\psi = \phi\circ\psi$.  In practice it is often easiest to look at
the $K$-points, then
$$
 \Res_{K/F}V(K) \simeq \prod_{\sigma\in \Hom_F(K,F^{\text{alg}})}
V_{K,\sigma}(K)
$$
 together with the action of $\Gal(F^{\text{alg}}/F)$ permuting the
factors according to its action on $\{\sigma_1 ,\ldots , \sigma_d\}$.
For further details see \cite{B.136}~1.3.
\enddefinition

\definition{1.1. Definition}
  An \dfn{algebraic group} over $F$ is an algebraic variety $G$ defined
over $F$ together with morphisms
$$
 \operatorname{mult} : G\times G \to G \qquad \text{ and } \qquad
\operatorname{inv} : G\to G,
$$
 both defined over $F$, and an element $e \in G(F)$ such that $G$ is a
group with identity~$e$, multiplication given by $\operatorname{mult}$,
and inverses given by $\operatorname{inv}$.  A \dfn{morphism of algebraic
groups} is a morphism of algebraic varieties which is also a group
homomorphism.
\enddefinition

 As a variety an algebraic group is smooth, since it contains an open
subvariety of smooth points and the group of translations $h\mapsto gh$
acts transitively.

\definition{1.1.1. Definition}
 An \dfn{abelian variety} is a complete connected algebraic group.  It
follows from this definition that an abelian variety is a smooth
projective variety and that its group law is commutative, see for example
\cite{B.96}, \cite{B.121}, \cite{B.68},
\cite{B.95}, \cite{B.69}, \cite{B.74}, \cite{B.79}
or \cite{B.134}.  It also follows that every morphism of
abelian varieties as varieties can be expressed as a composition of a
homomorphism with a translation, though of course only homomorphisms are
morphisms of abelian varieties as algebraic groups.
  In this appendix we will only be dealing with abelian varieties defined
over $\CC$, that is, complex abelian varieties.  When $A$ is a complex
abelian variety then the manifold underlying $A(\CC)$ is a complex torus.
\enddefinition

\definition{1.1.2. Definition}
 An affine algebraic group is also called a \dfn{linear algebraic group.}
This is justified by the fact that an affine algebraic group is
isomorphic, over its field of definition, to a closed subgroup of
$\GL(n)$ for some ~$n$, see \cite{B.14}, \cite{B.55},
\cite{B.52}, \cite{B.119}, \cite{B.133}.
\enddefinition

\definition{1.2. Definition}
  Any affine algebraic group that is isomorphic (as an algebraic group) to
the diagonal subgroup of $\GL(n)$ for some $n$ is called an \dfn{algebraic
torus.}  For additional basic exposition on algebraic tori see
\cite{B.14}~III.8 or \cite{B.55}~\S16.
\enddefinition

\example{1.2.1. Example}
  Our most basic and important example of an algebraic torus is $\Gm :=
\GL(1)$.  {\it A priori\/} $\Gm = \Gm_{/\QQ}$ is defined over $\QQ$, and
thus $\Gm(F) = F^\times$ for any field $F$ containing~$\QQ$.  Similarly,
with the conventions of ~1.0, $\Gm_{/F}(K) = K^\times$ when $K$ is a field
containing $F$.
\endexample

\example{1.2.2. Example}
  We may also apply the restriction of scalars functor to an algebraic
torus.  For the purposes of this appendix, one of the most important
examples that we will use later is
$$
 \SS := \Res_{\CC/\RR}\Gm_{/\CC} .
$$
 Then $\SS(\RR) = \CC^\times$ and $\SS(\CC) \simeq \CC^\times \times
\CC^\times$, where these last two factors are interchanged by complex
conjugation.  In particular $\SS(\RR)$ embeds as the diagonal in
$\SS(\CC)$.
\endexample

\definition{1.3. Definition}
  A connected linear algebraic group of positive dimension is said to be
\dfn{semisimple} if it has no closed connected commutative normal
subroups except the identity.  A (Zariski-connected) linear algebraic
group $G$ is said to be \dfn{reductive} if is the product of two
(Zariski-connected) normal subgroups $G_{\text{ab}}$ and
$G_{\text{ss}}$, where $G_{\text{ab}}$ is an algebraic torus and
$G_{\text{ss}}$ is semisimple, and $G_{\text{ab}} \cap G_{\text{ss}}$ is
finite.  {\it A~fortiori\/} any semisimple group is reductive.
\enddefinition

\definition{1.4. Definition}
  Recall that a \dfn{representation} of a group $G$ is a homomorphism
$\rho : G \to \GL(V)$ from $G$ to the automorphism group of a vector
space~$V$.  Such a representation may be referred to as $(\rho, V)$ or
simply by $\rho$ or by~$V$.  If $(\sigma, W)$ is another representation of
$G$, a map $\psi: V\to W$ such that $\sigma(g)\circ\psi = \psi \circ
\rho(g)$ for all $g\in G$ is said to be \dfn{$G$-linear} or
\dfn{$G$-equivariant.}  In this case, if $\psi$ is an isomorphism the
representations $(\rho,V)$ and $(\sigma, W)$ are said to be
\dfn{equivalent.}  We frequently identify equivalent representations.

  A subrepresentation of a representation is defined in the natural way,
and a representation is said to be \dfn{irreducible} if it contains no
nontrivial subrepresentations.  Further, given representations $(\rho,V)$
and $(\sigma,W)$ of $G$ we may form their direct sum or their tensor
product.  Thus the $r^{\text{th}}$ exterior power $(\twedge^r \rho,
\twedge^r V)$ of a representation arises naturally as a subrepresentation
of the $r$-fold tensor product of the representation $(\rho, V)$ with
itself.

  Let $V\dual = \Hom(V,F)$ denote the dual space to $V$ (if $V$ is a
vector space over $F$), and let $\angled{\ ,\ } :V\times V \to F$ be the
natural pairing.  Then $\rho$ induces a representation $\rho\dual: G\to
\GL(V\dual)$, called the \dfn{dual,} or \dfn{contragredient
representation,} of $G$.  It is defined by requiring
$$
  \angled{\rho\dual(g)v\dual,\rho(g)v} = \angled{v\dual,v};
$$
  concretely this means that $\rho\dual(g) = \tr\rho(g)^{-1}$.

  When $G$ is a subgroup of $\GL(V)$, for some vector space $V$, in
particular when $G$ comes as a subroup of $\GL(n)$, the group
of invertible $n\times n$ matrices, then by the \dfn{standard
representation} of $G$ we mean the natural inclusion $G\hra \GL(V)$.
\enddefinition

\subhead 1.5. Examples of semisimple and reductive groups
\endsubhead
  The examples that will be of interest to us are all classical
groups, defined from the outset as subgroups of $\GL(n)$.

\example{1.5.1.  Example}
  The first basic example is $\SL(n)$, the subgroup of $\GL(n)$ of
matrices of determinant~$1$.  For $n\ge 2$, $\SL(n)$ is semisimple.
It follows that $\GL(n)$ is reductive, as it is the product of its
subgroup of diagonal matrices and $\SL(n)$.
\endexample

\example{1.5.2. Example}
  Let $F$ be a subfield of the real numbers, in particular $\RR$ itself,
let $K$ be an imaginary quadratic extension of~$F$, and let $V$ be a
vector space over $K$.  Then a \dfn{Hermitian form} on $V$ is an
$F$-bilinear form $H:V\times V \to K$ such that $H(v,u) = \sigma(H(u,v))$,
where $\sigma$ is the nontrivial automorphism of $K$ over $F$, the
restriction of complex conjugation.  Then the \dfn{unitary group}
$\U(V,H)$ is the subgroup of $g\in \GL(V)$ such that $H(gu,gv) = H(u,v)$,
and the \dfn{special unitary group} $\SU(V,H)$ is the subgroup of $\U(V)$
of elements of determinant~$1$.  Note that $\U(V,H)$ and $\SU(V,H)$ are
algebraic groups defined over~$F$.  When $F=\RR$ and $H$ can be
represented by a diagonal matrix with $p$~$1$'s and $q$~$(-1)$'s then we
may write $\U(p,q)$ or $\SU(p,q)$ for $\U(V,H)$ or $\SU(V,H)$,
respectively; when $q=0$, that is when $H$ is equivalent to the standard
form $(u,v) \mapsto \tr\bar u \cdot v$, we write $\U(n)$ or $\SU(n)$.   As
a particular special case, note that $U(1)$ is defined over $\RR$, and
$\U(1,\RR)$ is the group of complex numbers of absolute value~$1$.
Similarly as in the previous example, $\SU(n)$ is semisimple and $\U(n)$
is reductive, for $n\ge 2$.
\endexample

\example{1.5.3. Example}
  When $E$ is a skew-symmetric bilinear form on a vector space $V$, that
is, $E(v,u) = -E(u,v)$, then the \dfn{symplectic group} $\Sp(V,E)$ is the
subgroup of $g\in \GL(V)$ such that $E(gu,gv)= E(u,v)$.  The
\dfn{symplectic similitude group} $\GSp(V,E)$ is the group of
$g\in\GL(V)$ such that there is a scalar $\nu(g)$ such that $E(gu,gv) =
\nu(g)E(u,v)$.  Thus $\GSp(V,E)$ contains $\Sp(E,V)$ as the subgroup of
\dfn{similtude norm} $\nu(g) =1$.  When $E$ can be represented by a
matrix of the form
$$
 \pmatrix 0 & I_n \\ -I_n & 0 \endpmatrix ,
$$
 then we may write $\Sp(2n)$ for $\Sp(V,E)$ and $\GSp(2n)$ for
$\GSp(V,E)$.  Moreover, $\Sp(2n)$ is semisimple.
\endexample

\example{1.5.4. Example}
  When $S$ is a symmetric bilinear form on a vector space $V$, then the
\dfn{special orthogonal group} $\SO(V,S)$ is the subgroup of $g \in
\SL(V)$ such that $S(gu,gv) = S(u,v)$.  In particular, if $S$ can be
represented by an identity matrix $I_n$ then we may write $\SO(n)$
instead of $\SO(V,S)$.  Also $\SO(n)$ is semisimple.
\endexample

\subhead{1.6. Lie algebras}
\endsubhead
  It will be useful later to have available some of the language of Lie
algebras, so we briefly recall some of the definitions here.  Our major
references for this paragraph (and the next two) are \cite{B.34},
\cite{B.54}, \cite{B.14}, \cite{B.55} and \cite{B.103}.

\definition{1.6.1. Definition}
  A \dfn{Lie algebra} is a vector space $\frak g$ together with a
skew-symmetric bilinear map, the \dfn{bracket} operation,
$$
 [\ ,\ ]: \frak g \times \frak g \to \frak g
$$
 satisfying the \dfn{Jacobi identity}
$$
 [X,[Y,Z]] + [Y,[Z,X]] + [Z,[X,Y]] =0.
$$
\enddefinition

\example{Example}
  Let $V$ be a vector space over a field $F$.  The fundamental example of
a Lie algebra is $\frak{gl}(V)$, which is $\End_F(V)$ as a vector space
with the bracket given by $[X,Y] := XY-YX$.
\endexample

\subsubhead{1.6.2. The Lie algebra of an algebraic group}
\endsubsubhead
  Let $G$ be a linear algebraic group over a field $F$.  In order to
review how a Lie algebra $\Lie(G)$ is associated to $G$, first recall
that when $A$ is any $F$-algebra then the Lie algebra of
\dfn{$F$-derivations} from $A$ to $A$ can be described as
$$
 \operatorname{Der}_F(A,A) := \{ X\in \frak{gl}(A) : X(f\cdot g) =
(Xf)\cdot g +f\cdot(Xg), \text{ for } f,g\in A\}
$$
 with the induced bracket $[X,Y] := XY-YX$.  Now let $A = F[G]$, the
coordinate ring of $G$ (as algebraic variety).  Then $G$ acts on $F[G]$ by
left translations:  $(\lambda_g f)(x) := f(g^{-1}x)$, for $f\in F[G]$ and
$g,x \in G$.  Then the set of \dfn{left invariant} derivations
$$
 L(G) := \{X\in \operatorname{Der}_F(F[G],F[G]) : \lambda_g \circ X =
X\circ \lambda_g \text{ for all }g\in G\}
$$
 is a Lie subalgebra of $\operatorname{Der}_F(A,A)$, and some authors
take this as the definition of $\Lie(G)$.  Next, recall that when
$\scO_e$ is the local ring at $e$ and $\frak m_e$ its maximal ideal, the
the \dfn{tangent space} to $G$ at the identity is
$$
 T(G)_e = \operatorname{Der}_F(\scO_e, \scO_e/\frak m_e) \isom
\Hom_{F-\text{mod}}(\frak m_e/\frak m_e^2,F) .
$$
 Of course $\scO_e/\frak m_e$ is the residue field of the local ring at
the identity~$e$.
 Then it turns out that evaluation at the identity $e$ of $G$ gives an
isomorphism from $L(G)$ to $T(G)_e$ to $G$ at the identity (see references
on linear algebraic groups
cited above).  Thus $\Lie(G)$ can also be defined as $T(G)_e$ with the
bracket operation induced by the isomorphism with $L(G)$.

\example{Example}
 As the notation suggests, $\frak{gl}(V) = \Lie(\GL(V))$.
\endexample

\remark{Remark}
 One motivation for working with Lie algebras is that for a connected
linear group $G$ a homomorphism $\phi:G \to H$ to another group $H$ is
determined by its differential at the identity.  In this way Lie algebras
linearize some of the problems of representation theory.  More generally,
when some property of $G$ is determined by an open neighborhood of the
identity, it is often more effective work with $\Lie(G)$.
\endremark

\subsubhead{1.6.3. The adjoint representations}
\endsubsubhead
  In general the \dfn{differential} of a morphism of (irreducible)
algebraic varieties $\phi: X\to Y$ at $x\in X$ is the linear map on
tangent spaces $d\phi_x:T(X)_x \to T(Y)_{\phi(x)}$ induced by $\phi^*:
\scO_{\phi(x)} \to \scO_x$.  In particular, $G$ acts on itself by inner
automorphisms
$$
  \operatorname{Int}_g : h\mapsto g h g^{-1},
$$
 and this action fixes the identity.  Then the differential of this map
$$
  \Ad(g) := d(\operatorname{Int}_g)_e : T(G)_e \to T(G)_e
$$
 defines the \dfn{adjoint representation} of $G$
$$
  \Ad: G\to \Aut(T(G)_e) : g\mapsto \Ad(g).
$$
 If we go one step further and take the differential of the adjoint
representation, we get a Lie algebra morphism
$$
  \ad := d\Ad : T(G)_e \to \End(T(G)_e) ,
$$
 with $\ad(X)(Y) = [X,Y]$.

\definition{1.6.4. Definition}
  Now let $G$ be (the real points of) a connected algebraic group over
{}~$\RR$, and let $K$ be a maximal compact subgroup.  Then a \dfn{Cartan
involution} of $G$ with respect to $K$ is an involutive automorphism of
$G$ whose fixed point set is precisely $K$.  The differential of a Cartan
involution is a Cartan involution of $\Lie(G)$, and the decomposition
$$
 \Lie(G) = \frak k + \frak p
$$
 where $\frak k$ is the fixed point set and $\frak p$ is the
$(-1)$-eigenspace, is called a \dfn{Cartan decomposition.}  It follows
that
$$
 [\frak k, \frak k] \subseteq \frak k, \qquad [\frak k, \frak p]
\subseteq \frak p, \qquad [\frak p, \frak p] \subseteq \frak k.
$$
\enddefinition

\definition{1.6.5. Definition}
  A semisimple real Lie algebra $\frak g$ with a Cartan decomposition
$\frak g = \frak k + \frak p$ is of \dfn{Hermitian type} if there exists
an element $H_0$ in the center of $\frak k$ such that $(\ad(H_0))^2 =
-1$ as endomorphisms of ~$\frak p$.  Recall that a reductive group is an
extension of a semisimple group by an algebraic torus.  We will say that a
real reductive algebraic group $G$ is of \dfn{Hermitian type} if the
abelian part of $G$ is compact and the semisimple part of its Lie algebra,
$\Lie(G)_{\text{ss}}$, is of Hermitian type in the previous sense..
\enddefinition

\example{1.7. Examples}
  Let $V$ be a vector space over a field $F$.  We have already noted that
$\frak{gl}(V) = \Lie(\GL(V))$ is $\End(V)$ as a vector space with the
bracket $[X,Y] = XY -YX$.

 The subgroup of $\frak{gl}(V)$ of endomorphisms with trace~$0$ is
$\frak{sl}(V) = \Lie(\SL(V))$.

 Suppose $V$ is a symplectic space of dimension ~$2n$ whose skew-symmetric
form is represented by
$$
  \pmatrix 0 & I_n \\ -I_n & 0 \endpmatrix
$$
 Then $\frak{sp}(2n) = \Lie(\Sp(2n))$ consists of matrices of the form
$$
 X = \pmatrix M & N \\ P & Q \endpmatrix
$$
 such that $N$ and $P$ are symmetric and $\tr M =-Q$.

  When $V$ has dimension ~$n$ and comes with a symmetric bilinear form
represented by an identity matrix, then $\frak{so}(n) = \Lie(\SO(n))$
consists of $n\times n$ skew-symmetric matrices.
\endexample

\subhead{1.8. The spin representations of $\frak{so}(n)$}
\endsubhead
  As a subgroup of $\frak{gl}(V) = \End(V)$ any of the Lie algebras above
naturally acts on $V$, and we may think of this as the standard
representation of the Lie algebra, cf.~1.4.  Moreover, up to equivalence,
all the representations of $\frak{sl}(n)$, respectively $\frak{sp}(2n)$,
occur in some tensor power of the standard representation.  However, that
is not the case for $\frak{so}(n)$, so here we briefly recall the complex
representation(s) that do not.

  Let $V$ be an $n$-dimensional vector space with a nondegenerate
symmetric bilinear form ~$Q$.  Then the quotient of the tensor algebra of
$V$ by the ideal generated by all elements of the form $v\tensor v -
Q(v,v)$ for $v\in V$ is called the \dfn{Clifford algebra} $C(V)$ of ~$V$.
Since this ideal preserves the property that an element of the tensor
algebra is the product of an even number of vectors, such products
generate a subalgebra $C^+(V)$ of $C(V)$ called the \dfn{even Clifford
algebra.}  For more details about Clifford algebras some good sources are
\cite{B.20} \cite{B.31} or \cite{B.19}.

  Now following \cite{B.34}~Chapter~20 we first
observe that since $C(V)$ and $C^+(V)$ are associative algebras they
determine Lie algebras with $[a,b] = a\cdot b - b \cdot a$.  Moreover,
$\frak{so}(V)$ embeds in $C^+(V)$ as a Lie subalgebra.  Roughly speaking,
on the one hand there is an embedding $\psi: \twedge^2V \to C^+(V)$, given
by $\psi(a\wedge b) = a \cdot b - Q(a,b)$, while on the other hand there
is an isomorphism $\phi: \twedge^2 V \simto \frak{so}(V) \subset
\frak{gl}(V)$ given by
$$
 \phi(a\wedge b)(v) = 2(Q(b,v)a - Q(a,v)b) .
$$

  Now suppose $n=2m$ is even (and the underlying field $F =\CC$).
Then $V$ can be written as the sum of two
$m$-dimensional isotropic subspaces, $V= W\oplus W'$ (meaning that the
restriction of $Q$ to $W$, respectively $W'$, is zero).  Then the key
lemma is that $C(V) \simeq \End(\twedge^*W)$, where $\twedge^*W$ signifies
the exterior algebra of ~$W$.  From this it follows that, if we write
$\twedge^*W = \twedge^+W \oplus \twedge^-W$ corresponding to even and odd
exterior powers, then
$$
 C^+(V) \simeq \End(\twedge^+W) \oplus \End(\twedge^-W).
$$
 Therefore the embedding of $\frak{so}(V)$ into $C^+(V)$ determines two
(inequivalent) representations of $\frak{so}(V)$, namely its actions on
$\twedge^+W$ and $\twedge^-W$ respectively.  These are referred to as
the \dfn{half-spin} representations of $\frak{so}(V)$, and their sum is
the \dfn{spin} representation.

 When $n=2m+1$ is odd, then we may write $V = W \oplus W' \oplus U$,
where $W$ and $W'$ are $m$-dimensional isotropic subspaces, as before,
and $U$ is $1$-dimensional and orthogonal to both $W$ and $W'$.  In this
case
$$
  C(V) \simeq \End(\twedge^*W) \oplus \End(\twedge^*W')
$$
 and $C^+(V) \simeq \End(\twedge^*W)$.  Thus the embedding of
$\frak{so}(V)$ into $C^+(V)$ determines a single \dfn{spin}
representation.

\subhead 1.9. Quaternion algebras
\endsubhead
  A \dfn{quaternion algebra} over a field $F$ (of characteristic not~$2$)
is a simple $F$-algebra of rank~$4$ whose center is precisely~$F$.  Over
the complex numbers, or any algebraically closed field of characteristic
not equal to ~$2$, there is up to isomorphism only one quaternion algebra,
the $2\times 2$ matrix algebra.  Over the real numbers, aside from the
$2\times 2$ matrix algebra there is up to isomorphism only one other
quaternion algebra, the \dfn{Hamiltonian} quaternion algebra ~$\Bbb H$,
generated over $\RR$ by $1$ and elements $i$ and $j$ such that
$$
 i^2 = j^2 =-1, \qquad\quad ij =-ji .
$$
 Note that $\Bbb H$ is a division algebra.  Over $\QQ$ there are
infinitely many non-isomorphic quaternion algebras, all of which except
the $2\times 2$ matrix algebra are (noncommutative) division algebras; see
for example \cite{B.132}.  A quaternion algebra over $\QQ$, or
more generally a quaternion algebra over a subfield of $\RR$, is said to
be \dfn{definite} or \dfn{indefinite} according as its tensor product with
$\RR$ is isomorphic to $\Bbb H$ or to $M_2(\RR)$.  In particular, an
indefinite quaternion algebra can be embedded into $M_2(\RR)$.

  In general, a quaternion algebra over $F$ has a basis consisting of $1$
and elements $\alpha$, $\beta$ and $\alpha\beta$ such that $\alpha^2$ and
$\beta^2$ are nonzero elements of $F$ and $\beta\alpha=-\alpha\beta$.
There is also a \dfn{canonical involution} on a quaternion algebra given
by
$$
 (a+b\alpha +c\beta +d\alpha\beta)' = a -b\alpha -c\beta +d\alpha\beta
\quad \text{ or } \quad \pmatrix a& b\\c&d\endpmatrix' = \pmatrix d& -b
\\ -c & a\endpmatrix .
$$
 Then the \dfn{reduced trace} and \dfn{reduced norm} of an element
$\gamma$ are $\gamma + \gamma'$ and $\gamma\cdot \gamma'$ respectively.

  A subring of a quaternion algebra that is also a lattice, i.e., a free
$\ZZ$-module of rank ~$4$, is called an \dfn{order} in the quaternion
algebra.  Very roughly, maximal orders play a similar role for quaternion
algebras as rings of integers do for number fields, except that maximal
orders in quaternion algebras need not be unique.

\definition{1.10. Definition}
 A \dfn{complex structure} on a real vector space $\Wr$ of dimension~$2g$
is given by any of the following equivalent data:
\roster
\item"(i)" A scalar multiplication by $\CC$ with which $\Wr$ is a
$g$-dimensional complex vector space.
\item"(ii)" An endomorphism $J\in \End(\Wr)$ such that $J^2 =-\Id$.
\item"(iii)" A homomorphism $h_1: \U(1) \to \GL(\Wr)$ of algebraic groups
over ~$\RR$ such that for $u \in \U(1)$ the action of $h(u)$ on $W_\CC$
has only $u^{\pm1}$-eigenspaces, each occuring with equal
multiplicity~$g$.
\item"(iv)" A homomorphism $h: \SS \to \GL(\Wr)$ of algebraic groups over
{}~$\RR$ such that for $(z,w)\in \SS(\CC)$ the action of $h(z,w)$ on
$W_\CC$ has only $z^1w^0$- and $z^0w^1$-eigenspaces, each occuring with
equal multiplicity~$g$.
\endroster
\enddefinition

\demo\nofrills{}
 To see the equivalence of these conditions, if a complex vector space
structure is given, let $J$ be the action of multiplication by
$i=\sqrt{-1}$.  Then $h_1(i) = J$ determines either $h_1$ or $J$ in terms
of the other.  Since $\SS = \Gm_{/\RR} \cdot \U(1)$, then $h$ is
determined by $h_1$ and $\RR$-linearity, or by $h(a+bi) = a\Id +bJ$.  And
$h$ in turn defines a scalar multiplication by ~$\CC$.  Also $h_1$ is the
restriction of $h$ to $\U(1) \subset \SS$, on account of which we
sometimes simply write $h$ instead of ~$h_1$.
\enddemo

\example{1.10.1. Example}
 In the notation of 1.6.5, $\ad(H_0)$ defines a complex structure on
{}~$\frak p$, with which $\frak p$ becomes a complex vector space.
\endexample

\example{1.10.2.  Example}
 To give a complex torus $V/L$, where $V$ is a $g$-dimensional
complex vector space and $L \subset V$ is a lattice, is the same
as giving a real $2g$-dimensional vector space $W$ together with a
complex structure, say $J$, and a lattice $L \subset W$.  We can go back
and forth between these two points of view by thinking of $W$ as the
real vector space underlying $V$ and $J$ as the induced complex
structure, or by thinking of $V$ as the complex vector space defined by
the pair $(W,J)$.  In particular it will sometimes be convenient below to
present a complex torus as a triple $(W, J, L)$ instead of in the form
$V/L$.
\endexample

\subhead{1.11. Complex abelian varieties}
\endsubhead
  After the definition given in 1.1.1, a complex abelian variety $A$ is a
complete, connected algebraic group over~$\CC$ whose group law is
necessarily commutative.  A morphism of abelian varieties will always be
taken to mean a morphism in the sense of algebraic groups (see~1.1).

\definition{1.11.1. Definition}
 A \dfn{Riemann form} on a complex torus $V/L$ is a nondegenerate,
skew-symmetric, real-valued, $\RR$-bilinear form $E:V\times V\to\RR$
such that
\roster
\item"(i)" $E(iv,iw) = E(v,w)$,
\item"(ii)" $(v,w) \mapsto E(v,iw)$ is symmetric and positive definite,
and
\item"(iii)" $E(v,w)\in\ZZ$ whenever $v,w\in L$.
\endroster
\enddefinition

 For a proof of the following proposition, see the references cited
in~1.1.1.

\proclaim{1.11.2. Proposition}
 A complex torus is the underlying manifold of a complex abelian variety
if and only if it admits a Riemann form.
\Qed
\endproclaim

\definition{1.11.3. Definition}
 A morphism of complex abelian varieties is called an \dfn{isogeny} if it
is surjective and has a finite kernel.  Given an isogeny
$\phi:A\to A'$ there exists a dual isogeny $\phi\dual: A'\to A$ such that
$\phi\spcheck\circ\phi = m \Id_{A}$ and $\phi\circ\phi\spcheck = m
\Id_{A'}$ for some positive integer~$m$ called the degree of~$\phi$.
Thus two complex abelian varieties are said to be \dfn{isogenous} iff
there exists an isogeny from one to the other, and being isogenous is an
equivalence relation.  An abelian variety is said to be \dfn{simple} iff
it is not isogenous to a product of (positive dimensional) abelian
varieties.
\enddefinition

 The following proposition is proved in Lecture~12, 12.25, or see the
references on abelian varieties cited above.

\proclaim{1.11.4. Proposition {\rm (Poincar\'e Reducibility Theorem)}}
 If $A$ is an abelian variety and $A' \subset A$ is an abelian subvariety,
then there exists an abelian subvariety $A'' \subset A$ such that $A' \cap
A''$ is finite and $A$ is isogenous to $A'\times A''$.  In particular, any
abelian variety is isogenous to a product of simple abelian varieties.
\endproclaim

\definition{1.11.5. Definition}
  Two Riemann forms $E$ and $E'$ are said to be \dfn{equivalent} iff there
exist positive integers $n$ and $n'$ such that $nE = n'E'$.  A
\dfn{polarization} of a complex torus $T$ is an equivalence class, say
{}~$[E]$, of Riemann forms on~$T$.  By a \dfn{polarized abelian variety}
we mean an abelian variety together with a choice of polarization.  In
light of 1.10.2 and 1.11.2, a polarized abelian variety is determined by
data $(W, J, L, E)$, where $W$ is an even-dimensional real vector space,
$J$ is a complex structure on ~$W$, $L$ is a lattice in ~$W$, and $E$ is
a Riemann form on the complex torus $(W,J,L)$.
\enddefinition

\subhead 1.12.. The endomorphism algebra of an abelian variety
\endsubhead
  For a complex abelian variety $A$ let $\End(A)$ denote its endomorphism
ring, and let
$$
 \EndoA := \End(A)\tensor_\ZZ \QQ .
$$

\proclaim{1.12.1. Lemma}
\roster
\item When $A$ and $A'$ are isogenous abelian varieties, $\EndoA \simeq
\Endo(A')$.
\item $\EndoA$ is a semisimple $\QQ$-algebra with a positive involution.
\endroster
\endproclaim

 Recall that an involution $\iota$ of $\EndoA$ is said to be positive if
for nonzero $\phi\in \EndoA$ the trace $\operatorname{Tr}(\phi\cdot
\phi^\iota) > 0$.

\demo{Proof}
  If $\phi:A\to A'$ is an isogeny and $\phi\dual:A'\to A$ is the dual
isogeny, then $\phi^*:\Endo(A') \to \EndoA$ and $\frac 1 m (\phi\dual)^*
: \EndoA \to \Endo(A')$ are mutually inverse ring homomorphisms.

 To prove part~2, first observe that in general the image of
a homomorphism of complex tori is a subtorus, and the kernel is a closed
subgroup whose connected component of the identity is a subtorus of finite
index in the full kernel.  Thus Schur's Lemma implies that $\EndoA$ is a
division algebra when $A$ is a simple abelian variety, and then the
semisimplicity of $\EndoA$ for general $A$ follows from the Poincar\'e
Reducibility Theorem.

 To see that $\EndoA$ has a positive involution $\iota$, let $E$
be a Riemann form on $A$.  Then $H(u,v) = E(u,iv) + iE(u,v)$ is a
Hermitian form on $\Wr$, and if we take $\iota$ to be the antiautomorphism
that takes $\phi\in\EndoA$ to its adjoint with respect to~$H$ then the
conditions making $E$ a Riemann form imply that $\iota$ is a positive
involution.
\Qed
\enddemo

\definition{Definition}
 The involution $\iota$ of $\EndoA$ described in the proof above is called
the \dfn{Rosati involution.}
\enddefinition

 Thus when $A$ is simple, $\EndoA$ is a division algebra over~$\QQ$ which
admits a positive involution.  Such algebras were classified by Albert
\cite{B.6}, \cite{B.7}, \cite{B.8}, see also \cite{B.109} and
\cite{B.79}.  The result is the following.

\proclaim{1.12.2.  Theorem {\rm (Albert classification)}}
 Let $A$ be a simple complex abelian variety.  Let $K$ be the center of
$\EndoA$ and let $K_0$ be the subfield of elements of $K$ fixed by the
Rosati involution.  Then $\EndoA$ is one of the following types:
\roster
\item"(I)" $\EndoA= K =K_0$ is a totally real algebraic number field, and
the Rosati involution acts as the identity.
\item"(II)" $K=K_0$ is a totally real number algebraic field, and $\EndoA$
is a division quaternion algebra over~$K$ such that every simple component
of $\EndoA\tensor_\QQ \RR$ is isomorphic to $M_2(\RR)$; there is an
element $\beta\in \EndoA$ such that
$\vphantom{\beta}^{\operatorname{t}}\beta = -\beta$, and $\beta^2 \in K$
is totally negative; and the Rosati involution is given by $\alpha^\iota =
\beta^{-1}\cdot\vphantom{\alpha}^{\operatorname{t}}\!\alpha\cdot\beta$.
\item"(III)" $K=K_0$ is a totally real number algebraic field, and
$\EndoA$ is a division quaternion algebra over~$K$ such that every
simple component of $\EndoA\tensor_\QQ \RR$ is isomorphic to the
Hamiltonian quaternion algebra $\Bbb H$ over~$\RR$; and $\alpha^\iota =
\tr\alpha$.
\item"(IV)" $K_0$ is a totally real number field, and $K$ is a totally
imaginary quadratic extension of $K$, and $\EndoA$ is division algebra
with center $K$, and the restriction of the Rosati involution to $K$ acts
as the restriction of complex conjugation to ~$K$.
\endroster
\endproclaim

 Thus we will say that an abelian variety or an abelian manifold ~$A$ is
\dfn{of type} (I), (II), (III) or (IV) if $A$ is simple and $\Endo(A)$ is
of that type in the classification above, or when $A$ is not simple, if it
is isogenous to a product of simple abelian varieties of that type.

\subhead 1.13. Examples of abelian varieties
\endsubhead
 We now introduce some basic constructions of complex abelian varieties
with various endomorphism algebras.  In the next sections we will more
fully analyze their Hodge structures.

\example{1.13.1. Elliptic curves}
  Let $E$ be a $1$-dimensional complex abelian variety, an \dfn{elliptic
curve.}  Then $E(\CC) \simeq \CC/(\tau\ZZ+\ZZ)$ for some $\tau\in\CC$ with
$\Im\tau >0$.  A Riemann form is given by the pairing
$$
 ((a\tau+b),(c\tau+d)) \mapsto \frac 1{\Im \tau}
\Im((a\tau+b)(c\bar\tau+d)) = ad-bc
$$
  It is an elementary exercise to show that $\Endo(E)$ can only be
isomorphic to $\QQ$ or to an imaginary quadratic field, say~$K$; in the
latter case $E$ is said to have \dfn{complex multiplication} by $K$.
Moreover $E$ has complex multiplication if and only if $\tau$ is
quadratic over~$\QQ$, in which case $\QQ(\tau) = K$.  The general elliptic
curve, whose period $\tau$ has algebraically independent transcendental
real and imaginary parts, has $\Endo(E)\simeq \QQ$.
\endexample

\example{1.13.2. Abelian varieties with multiplication by an imaginary
quadratic field}
  We will say that an abelian variety has multiplication by an imaginary
quadratic field $K$ if there is an embedding $K\hra \EndoA$.  Following
\cite{B.135}, to construct a simple complex abelian variety $A$ with
multiplication $K$, let $W$ be an $n$-dimensional vector space over ~$K$.
Then $\Wr = W\tensor_\QQ \RR$ may be identified with $\CC^n$, but we may
also twist the complex structure as follows.  Write $\Wr = \Wr' \oplus
\Wr''$ as the direct sum of two subspaces over~$\CC$, and define a complex
structure $J\in\End_\RR(\Wr)$ by $Jw' = i w'$ for $w'\in \Wr'$ and $J w''
= -i w''$ for $w''\in\Wr''$.  Then with this complex structure, $\alpha
\in K$ acts on $(\Wr,J)$ by $w'\mapsto \alpha w'$ for $w'\in \Wr'$ and
$w'' \mapsto \bar\alpha w''$ for $w''\in \Wr''$.  In particular, when
$L\subset W$ is a lattice such that $L \tensor_\ZZ\QQ = W$, then the set
of $\alpha \in K$ such that $\alpha\cdot L \subseteq L$ is a subring of
$K$
commensurable with the ring of integers of~$K$.  Thus $A(\CC)=(\Wr/L, J)$
is a complex torus with an embedding $K\hra \Endo(A)$.  To exhibit
a Riemann form for $A$, let $H$ be any $\QQ$-valued Hermitian form on
$W\times W$ which as a $\CC$-valued form on $\Wr\times\Wr$ is positive
definite on $\Wr'$ and negative definite on $\Wr''$, meaning in particular
that these two subspaces are orthogonal with respect to ~$H$.  Then $H = S
+ iE$ for $\RR$-valued forms $S$ and $E$, and the imaginary part~$E$ of
{}~$H$ is a Riemann form for ~$A$.

  Let $n'= \dim_\CC \Wr'$ and $n'' = \dim_\CC \Wr''$.  Then $n' +n'' =n$,
and by \cite{B.109} Thm.5, when $n\ge 3$ then both $n'$ and $n''$
are positive, and $n=2$ does not occur.   If $n' = n''$ the pair $(A,K)$
is said to be an abelian variety \dfn{of Weil type.}  Equivalently, an
abelian variety of Weil type is a pair $(A,K)$ consisting of an abelian
variety $A$ and an imaginary quadratic field $K$ with an embedding $K\hra
\Endo(A)$ such that for $\alpha\in K$ the corresponding endomorphism has
the eigenvalues ~$\alpha$ and ~$\bar\alpha$ with equal multiplicity.  The
example of Mumford in \cite{B.88}, see Lecture~7, 7.23--7.28, is
one example of an abelian variety of Weil type.

  Another case that will arise below is when $n'$ and $n''$ are relatively
prime.  We will refer to a pair $(A,K)$ satisfying this condition as an
abelian variety \dfn{of Ribet type,} see \cite{B.94}~Thm.3.
\endexample

\example{1.13.3. Simple abelian varieties of odd prime dimension}
  Let $A$ be a simple complex abelian variety of odd prime dimension~$g$.
Then reading off from the tables in \cite{B.85}, the cases that occur
are:  $\EndoA \simeq \QQ$, which is the general case; or
$\EndoA$ is a totally real number field of degree~$g$ over ~$\QQ$; or
$\EndoA$ is an imaginary quadratic field, in which case $A$ is of Ribet
type; or $A$ is of CM-type, that is, a totally imaginary
quadratic extension of a totally real field of degree~$g$ over~$\QQ$.
\endexample

\example{1.13.4. Simple abelian fourfolds}
 By reading the tables in \cite{B.85}, the endomorphism algebras that
can occur for a simple abelian fourfold are, by Albert type:
\roster
\item"(I)" $\EndoA$ is $\QQ$, or a real quadratic field, or a totally real
quartic field;
\item"(II)" $\EndoA$ is an indefinite division quaternion algebra
over $\QQ$ or a totally indefinite division quaternion algebra over a real
quadratic field;
\item"(III)" $\EndoA$ is a definite division quaternion over $\QQ$;
\item"(IV)" $\EndoA$ is an imaginary quadratic field, in which case it can
only be of Ribet type with $\{n',\,n''\} = \{1,\,3\}$ or of Weil type with
$n'=n''=2$, or else $\EndoA$ is a CM-field of degree~$4$ or $8$
over~$\QQ$.
\endroster
\endexample

\example{1.13.5. Abelian varieties with real multiplication}
  An abelian variety of type~(I) is sometimes said to have \dfn{real
multiplication.}  To construct an example of a simple abelian variety with
real multiplication, let $K$ be a totally real number field with
$[K:\QQ]=g$, and let $\scO$ be the ring of integers of ~$K$.  Then there
are $g$ distinct embeddings $\alpha\mapsto \alpha^{(j)}$ of $K$ into
$\RR$.  Let $\tau_j \in\CC$ with $\Im\tau_j>0$, for $1\le j\le g$.  Then
the image of $\scO \oplus \scO$ under the map
$$
 (\alpha,\beta)\mapsto (\alpha^{(1)}\tau_1 + \beta^{(1)}, \ldots ,
\alpha^{(g)}\tau_g + \beta^{(g)})
$$
 is a lattice $L \subset \CC^g$, and $A = \CC^g/L$ is a complex abelian
variety.  A Riemann form is given by
$$
 E(z,w) = \sum_{j=1}^g (\Im\tau_j)^{-1} \Im(z_j\bar w_j),
$$
 where $z,w\in\CC^g$.  Then $K\hra \Endo(A)$, and this is an isomorphism
for general $(\tau_1,\dots,\tau_g)$.  It can be shown that any simple
abelian variety $A$ for which $\EndoA$ is a totally real number field is
isogenous to one which can be constructed as we have here \cite{B.27}.
\endexample

\example{1.13.6. Abelian varieties of CM-type}
   Recall that an algebraic number field $K$ is said to be a
\dfn{CM-field} iff it is a totally imaginary quadratic extension of a
totally real number field~$K_0$.  The embeddings of a CM-field $K$ into
$\CC$ come in complex conjugate pairs.  Then \dfn{CM-type} for $K$ is a
subset $S\subset\Hom(K,\CC)$ containing exactly one from each pair of
conjugate embeddings, so that $\Hom(K,\CC) = S \cup \ol S$.

 A simple abelian variety $A$ is said to be of \dfn{CM-type,} or to have
\dfn{complex multiplication} by ~$K$, iff there exists a field $K\hra
\Endo(A)$ such that $[K:\QQ] \ge 2\dim A$, in which case equality holds,
$K\simeq \Endo(A)$ and $K$ is a CM-field, see \cite{B.115}
or \cite{B.70}.  More generally, an abelian variety may be
said to be of \dfn{CM-type} if it is isogenous to a product of simple
abelian varieties of CM-type, or equivalently if $\EndoA$ contains a
commutative semisimple $\QQ$-algebra $R$ with $[R:\QQ] = 2\dim A$.

 To construct a simple abelian variety of CM-type, let $K$ be a CM-field
with totally real subfield $K_0$ such that $[K:\QQ]=2g$, and let $S$ be a
CM-type for ~$K$, and let $\scO$ be the ring of integers of ~$K$.  Then
$K\tensor_\QQ \RR \simeq \CC^g$.  If we embed $\scO \hra
\CC^g$ by $\alpha \mapsto (\sigma \alpha)_{\sigma\in S}$ and let $L$ be
the image of this map, then $A=\CC^g/L$ is an abelian variety.  To
construct a Riemann form, choose an element $\beta\in \scO$ such that $K=
K_0(\beta)$, and $-\beta^2$ is a totally positive element of $K_0$, and
$\Im(\sigma\beta)>0$ for $\sigma \in S$.  Then
$$
 E(z,w) = \sum_{j=1}^g \sigma_j(\beta) (z_j\ol w_j - \ol z_j w_j)
$$
  is a Riemann form, where $z,w\in\CC^g$.  When $z= (\sigma
\alpha_1)_{\sigma\in S}$ and $w= (\sigma \alpha_2)_{\sigma\in S}$ with
$\alpha_1, \alpha_2 \in K$, then $E(z,w) =
\operatorname{Tr}_{K/\QQ}(\beta \alpha_1 \ol \alpha_2)$.  Moreover, $A$
has complex multiplication by $K$, with $\alpha \in \scO$ acting by $z_j
\mapsto \sigma_j(\alpha) z_j$ for $1\le j\le g$; for more detail see
\cite{B.70}~\S1.4.  In addition, it can be shown that any simple
abelian variety of CM-type is isogenous to one such as we constructed
above \cite{B.27}.
\endexample

\example{1.13.7. Abelian surfaces with quaternionic multiplication}
 A simple abelian variety of type~(II) may be said to have
\dfn{quaternionic multiplication,} or sometimes, to be of \dfn{QM-type.}
The simplest example is an abelian surface $A$ whose endomorphism algebra
is an indefinite division quaternion algebra $D$ over $\QQ$.  To construct
such an abelian surface, let $\scO$ be an order in ~$D$, fix an
embedding $j: D \hra M_2(\RR)$, and let $\tau\in\CC$ with $\Im\tau >0$.
Then the image of $\scO$ under the map $\psi: \alpha \mapsto
j(\alpha)\left(\smallmatrix \tau \\ 1\endsmallmatrix
\right)$ is a lattice $L \subset \CC^2$, and $A=\CC^2/L$ is a
$2$-dimensional abelian variety. In the special case that $D\simeq
M_2(\QQ)$, then $\scO$ is commensurable with $M_2(\ZZ)$ and $A$ is
isogenous to the product of two isogenous elliptic curves.  Otherwise $D$
is a division algebra and $A$ is a simple abelian surface.  In this latter
case there is an element $\beta\in \scO$ such that
$\beta' = -\beta$ and $\beta^2<0$ in
$\QQ$ (recall that $\beta\mapsto \beta'$ is the canonical involution
on~$D$).  Then a Riemann form on $A$ is given by
$$
E((z_1,z_2),(w_1,w_2)) = \frac 1{\Im \tau} \Im \operatorname{Tr}\left(
j(\beta)\cdot \pmatrix z_1\ol w_2 & \ol z_1 w_1 \\ z_2\ol w_2 & \ol z_2
w_1 \endpmatrix \right) =  \operatorname{Tr} (\beta \cdot \alpha_1' \cdot
\alpha_2 )
$$
  when $\left(\smallmatrix z_1 \\ z_2 \endsmallmatrix \right) =
\psi(\alpha_1)$ and $\left(\smallmatrix w_1 \\ w_2 \endsmallmatrix \right)
= \psi(\alpha_2)$ for $\alpha_1,\alpha_2\in D\tensor_\QQ \RR$.
Furthermore, for $\gamma\in \scO$, multiplication by $j(\gamma)$ on
$\left(\smallmatrix z_1 \\ z_2\endsmallmatrix \right) \in\CC^2$ preserves
$L$.  Then by tensoring with $\QQ$ we get an inclusion $D\hra \Endo(A)$,
which is an isomorphism for general~$\tau$.

  We leave it as an exercise for the reader to combine the construction of
this example with that of 1.13.5 to obtain an arbitrary simple abelian
variety of type~(II).
\endexample

\example{1.13.8. General abelian varieties}
  When $A$ is a $g$-dimensional complex abelian variety then $A(\CC)
\simeq \CC^g/(T\ZZ^g + \ZZ^g)$ for some $T$ in the \dfn{Siegel upper
half-space of genus~$g$,} consisting of symmetric complex $g\times g$
matrices with positive-definite imaginary part.  Then generalizing the
$1$-dimensional case, a Riemann form is given by
$$
 E(\bold z, \bold w) =
 \Im(\tr\bold z (\Im T)^{-1} \ol{\bold w}) = \tr \bold a \cdot \bold d -
\tr \bold b \cdot \bold c
$$
 when $\bold z = T\bold a + \bold b$ and $\bold w = T \bold c + \bold d$,
with $\bold a, \bold b, \bold c, \bold d \in \RR^g$.  Then for general
$T$, that is, when all the real and imaginary parts of the distinct
entries of $T$ are algebraically independent real transcendental numbers,
$\EndoA \simeq \QQ$.
\endexample

\head 2. The Hodge, Mumford-Tate and Lefschetz groups of an abelian
variety
\endhead
\rightheadtext{The Hodge and Mumford-Tate groups}
\nopagebreak
\subhead 2.1. Rational Hodge structures
\endsubhead
   A real Hodge structure is a natural generalization of a complex
structure, and a rational Hodge structure is a real Hodge structure with
an underlying $\QQ$-structure.  For a rational vector space $V=V_\QQ$ we
write $\Vr =V\tensor_\QQ \RR$ and $V_\CC = V\tensor_\QQ \CC$.

\definition{2.1.1. Definition}
 A \dfn{rational Hodge structure of weight~$n$} consists of a
finite-dimensional $\QQ$-vector space $V$ together with any of the
following equivalent data:
\roster
\item"(i)" A decomposition $V_\CC  = \bigoplus_{p+q=n} V^{p,q}$, called
the {\rm Hodge decomposition,} such that $\ol{V^{p,q}} = V^{q,p}$.
\item"(ii)" A decreasing filtration $F_H^rV_\CC$ of $V_\CC$, called the
{\rm Hodge filtration,} such that $F_H^rV_\CC \oplus \ol{F_H^{n-r+1}V_\CC}
=V_\CC$.
\item"(iii)" A homomorphism $h_1:\U(1) \to \GL(\Vr)$ of real algebraic
groups, and also specifying that the weight of the Hodge structure is
{}~$n$.
\item"(iv)" A homomorphism $h:\SS\to\GL(\Vr)$ of real algebraic groups
such that via the composition $\Gm_{/\RR} \hra \SS \to \GL(\Vr)$ an
element $t\in \Gm_{/\RR}$ acts as $t^{-n}\cdot\Id$.
\endroster

\demo\nofrills{}
 To see that the data (i)--(iv) are equivalent, the Hodge decomposition
and the Hodge filtration are related by $F_H^rV_\CC = \bigoplus_{p\ge
r}V^{p,n-p}$ and $V^{p,n-p} = F_H^pV_\CC \cap \ol{F_H^{n-p}V_\CC}$.  The
homomorphism $h$ and the Hodge decomposition are related by $h(z)\cdot v =
z^{-p}\bar z^{-q} v$ for $v\in V^{p,q}$.  The homomorphism $h_1$ can be
obtained as the restriction of $h$; conversely, $V^{p,q}$ can be recovered
as the subspace of $V_\CC$ on which $h_1(u)$ acts as $u^{q-p}$, provided
$n=p+q$ is specified.
\enddemo
\enddefinition

\example{2.1.2. Examples}
 1. Let $V$ be a $\QQ$-vector space of even dimension, and suppose
$h:\SS \to \GL(V_\RR)$ defines a complex structure on ~$V_\RR$.  Then
$V$ is a rational Hodge structure of weight~$-1$:  criterion 1.10(iv)
in the definition of complex structure immediately implies criterion
2.1.1(iv).  Alternatively, the $z^1w^0$- and $z^0w^1$-eigenspaces of
$h(z,w)$ acting on $V_\CC$ are $V^{-1,0}$ and $V^{0,-1}$ respectively,
and these are complex conjugate as require by 2.1.1(i) because
$h$ is defined over $\RR$.

  2.  When $X$ is (the analytic space underlying) a complex projective
variety, or more generally, a compact K\"ahler manifold, then
$H^n(X,\QQ)$ is a rational Hodge structure of weight~$n$.

  3.  In general, when $V$ is a rational Hodge structure of weight~$n$,
its dual $V\dual$ is a rational Hodge structure of weight~$-n$, and
$V^{\tensor r}\tensor_\QQ (V\dual)^{\tensor s}$ is a rational Hodge
structure of weight $(r-s)n$.
\endexample

\remark{Remark}
 The convention that $(z,w)\in \SS(\CC)$ acts as $z^{-p}w^{-q}$ on
$V^{p,q}$ follows \cite{B.25}~\S1 and \cite{B.27}.  But
as we will see below it really is the more natural choice in the context
of Hodge structures associated to complex abelian varieties.
\endremark

\definition{2.1.3. Definition}
 The \dfn{type} of a rational Hodge structure $V$ is the set of pairs
$(p,q)$ such that $V^{p,q}\ne 0$.
\enddefinition

\definition{2.1.4. Definition}
  When $V$ is a rational Hodge structure of even weight~$2p$, the subspace
of \dfn{Hodge vectors} in $V$ is the subspace of $1$-dimensional rational
sub-Hodge structures of~$V$,
$$
 \Hdg(V) : = V_\QQ \cap V^{p,p} .
$$
\enddefinition

\definition{2.1.5. Definition}
  A \dfn{morphism} of rational Hodge structures $\phi:V_1\to V_2$ is a
$\QQ$-vector space map on the underlying vector spaces such that over
$\CC$
\roster
\item"(i)" $\phi(V_1^{p,q}) \subseteq V_2^{p,q}$, for all $p,q$; or
\item"(ii)" $\phi(F_H^rV_{1\CC})\subseteq F_H^r V_{2\CC}$, for all ~$r$;
or
\item"(iii)" $\phi$ commutes with the action of $U(1)$, and preserves the
common weight of $V_1$ and $V_2$; or
\item"(iv)" $\phi$ commutes with the action of $\SS$, preserving the
action of $\Gm_{/\QQ} \hra \SS$.
\endroster
\enddefinition

\definition{2.1.6. Definition}
  A \dfn{polarization} of a rational Hodge structure $V$ is a morphism of
rational Hodge structures $\psi : V\tensor V \to \QQ(-n)$ such that the
real-valued form $(u,v) \mapsto \psi(u,h(i)v)$ on $\Vr$ is symmetric and
positive-definite.  Here $\QQ(m)$ is the is the vector space $\QQ$ as a
one-dimensional rational Hodge structure of type~$(-m,-m)$, for $m\in\ZZ$.
\enddefinition

\subsubhead 2.1.7. The rational Hodge structure associated to an abelian
variety
\endsubsubhead
  When we speak of the rational Hodge structure associated to an abelian
variety $A$ we always mean the rational Hodge structure $H^1(A,\QQ)$.
Moreover, any morphism $\phi : A\to A'$ of abelian varieties induces a
morphism of rational Hodge structures $\phi^*: H^1(A',\QQ) \to
H^1(A,\QQ)$.  In particular, it is easy to check that when $\phi$ is an
isogeny it induces an isomorphism on the associated rational Hodge
structures.  Thus, up to isomorphism, the rational Hodge structure
associated to an abelian variety depends only on its isogeny class.
Similarly, an element of $\EndoA$ induces an endomorphism on the rational
Hodge structure $H^1(A,\QQ)$.

\remark{2.1.8.  Notation}
 For a complex abelian variety $A$ we will regularly let $W= H_1(A,\QQ)$
and $V=H^1(A,\QQ) = W\dual$.  Further, we denote the Hodge classes of $A$
by
$$
 \Hdg^p(A) := H^{2p}(A,\QQ) \cap H^{p,p}(A) \qquad\quad \Hdg(A) :=
\bigoplus_p \Hdg^p(A) ,
$$
 and let $\Div^p(A)\subset H^{2p}(A,\QQ)$ be the $\QQ$-linear span of
$p$-fold intersections of divisors on $A$, and $\Div(A) := \bigoplus_p
\Div^p(A)$.
\endremark

\definition{2.2. Definition}
 Let $(V, h:\SS\to\GL(\Vr))$ be a rational Hodge structure.  The
\dfn{Hodge group} $\Hg(V)$ of~$V$, also called the \dfn{special
Mumford-Tate group} of~$V$, is the smallest algebraic subgroup of $\GL(V)$
defined over $\QQ$ such that $h(\U(1)) \subset \Hg(\Vr)$.  The
\dfn{Mumford-Tate group} $\MT(V)$ of $V$ is the smallest algebraic
subgroup of $\GL(V)$ defined over $\QQ$ such that $h(\SS)\subset
\MT(\Vr)$.  As a matter of notation we let $\hg(V) := \Lie(\Hg(V))$ and
$\mt(V) := \Lie(\MT(V))$ denote their respective Lie algebras, as
subalgebras of $\End_\QQ(V)$.

 When $A$ is a complex abelian variety then by the Hodge or the
Mumford-Tate group of $A$ we mean $\Hg(A) := \Hg(H^1(A,\QQ))$ and $\MT(A)
:= \MT(H^1(A,\QQ))$ respectively.
\enddefinition

\remark{Remark}
  The Hodge group of an abelian variety was introduced in \cite{B.77},
see also \cite{B.78}.  The general notion of the
Mumford-Tate group of a rational Hodge structure seems to appear first in
\cite{B.24}~\S7, and from a rather abstract point of view, in
\cite{B.97}.  A thorough analysis of the Mumford-Tate groups of
Hodge structures that can be generated by the Hodge structures of abelian
varieties can be found in \cite{B.25}~\S1, while the best place
to find proofs of the basic properties of Mumford-Tate groups in general
is \cite{B.27}~\S3.
\endremark

 The first properties to follow directly from the definition are the
following.

\proclaim{2.3. Lemma}
 Let $V$ be a rational Hodge structure.  Then
\roster
\item"(i)" $\Hg(V)$ and $\MT(V)$ are connected linear algebraic groups;
\item"(ii)" $\Hg(V)\subseteq \SL(V)$;
\item"(iii)" $\MT(V) = \Gm \cdot \Hg(V)$.
\Qed
\endroster
\endproclaim

 The vital role of Mumford-Tate groups in analyzing Hodge structures comes
from the following fact.

\proclaim{2.4. Proposition}
 Let $V$ be a rational Hodge structure, and $r,s\in\NN$.  Then $\MT(V)$
acts on the rational Hodge structure $V^{\tensor r}\tensor
(V\dual)^{\tensor s}$, and the rational $\MT(V)$-subrepresentations in
$V^{\tensor r}\tensor (V\dual)^{\tensor s}$ are precisely the rational
sub-Hodge structures of $V^{\tensor r}\tensor (V\dual)^{\tensor s}$.
\endproclaim

 The action of $\MT(V)$ on $V$ extends ``diagonally'' to an action on
$V^{\tensor r}$, and the action of $\MT(V)$ on $V\dual$ is the
contragredient of its action on $V$, as in~1.4.

\demo{Proof}
 To simplify the notation, let $T:= V^{\tensor r}\tensor (V\dual)^{\tensor
s}$.  Then action of $\GL(V)$ on $V$ induces an action of $\GL(V)$, and
thus of $\MT(V)$ on $T$.  Now suppose first that $W\subset T$ is a
$\QQ$-rational subspace preserved by the $\MT(V)$-action.  Then over $\RR$
the composition $h:\SS \hra \MT(\Vr) \to \GL(\Wr)$ describes the sub-Hodge
structure on $\Wr \subset T\tensor \RR$.  Conversely, if $W\subset T$ is a
rational sub-Hodge structure, then $W$ is a rational subspace of $T$ such
that $\Wr \subset T_\RR$ is preserved by the action of $h(\SS)$.
Therefore $W$ is preserved by the action of $\MT(V)$ on $T$.
\Qed
\enddemo

\proclaim{Corollary}
 Let $V$ be a rational Hodge structure of weight~$n$.  Then for any
$r,s\in \NN$ such that $(r-s)n =2p$,
$$
 \Hdg^p(V^{\tensor r}\tensor (V\dual)^{\tensor s}) = (V^{\tensor r}\tensor
(V\dual)^{\tensor s})^{\Hg(V)}
$$
\endproclaim

\remark{Remark}
  Hodge and Mumford-Tate groups have proved to be a powerful for studying
the Hodge conjecture for abelian varieties.  Typically the starting point
is this corollary to Proposition~2.4, which implies that the Hodge cycles
in $H^*(A^n,\QQ)$, say, for a complex abelian variety $A$, are precisely
the invariants under the action of $\Hg(A)$.  Then, the
problem is to determine enough about $\Hg(A)$ to be able to describe its
invariants, or at least determine their dimension.  In many cases it is
possible to show this way that the space of Hodge cycles is generated by
those of degree~$2$, in other words, by divisors on $A$, and then in these
cases the Hodge conjecture is verified.  In other cases it is possible to
show that the space of Hodge cycles is {\sl not\/} generated by divisors,
but still something can be said about the dimension of the space of Hodge
cycles.  However, the actual computations are sometimes quite technical. 
There are many variously narrow results, and a few general results, as we
will try to show in the following sections.
\endremark

 Since a polarization of a complex abelian variety induces a polarization
on its associated Hodge structure, the following proposition insures that
all the Hodge and Mumford-Tate groups with which we will work in this
appendix are reductive (definition~1.4); compare
\cite{B.25}~Principe~1.1.9.

\proclaim{2.5. Proposition}
 Let $V$ be a polarizable rational Hodge structure.  Then $\MT(V)$ and
$\Hg(V)$ are reductive.
\endproclaim

\demo{Proof {\rm (after \cite{B.27}~Prop.3.6)}}
 Suppose $h:\U(1) \to \GL(\Vr)$ defines the Hodge structure on $V$, and
let $\psi : V\tensor V \to \QQ(-n)$ be a polarization.  Then for $u,v \in
V_\CC$ and $g \in \Hg(V,\CC)$ (the complexification, or complex points
of $\Hg(V)$)
$$
 \psi(u,h(i) \bar v) = \psi(gu,gh(i)\bar v) = \psi(gu, h(i) (h(i)^{-1} g
h(i)) \bar v) = \psi(gu, h(i)\ol{g^*v}),
$$
 where the first equality holds because $\psi$ is a Hodge structure
morphism, and $g^* := h(i)^{-1} \bar g h(i)$.  Therefore the
positive-definite form on $\Vr$ given by $(u,v)\mapsto \psi(u,h(i)v)$ is
invariant under the real form $\Hg^*(\Vr)$ of $\Hg(V,\CC)$ fixed by the
involution $g\mapsto g^*$.  It follows that $\Hg^*(\Vr)$ is a compact real
form of $\Hg(V)$, and thus all its finite-dimensional representations are
semisimple.  This is equivalent to $\Hg^*(\Vr)$ being reductive, see
for example \cite{B.103}~I.3.  Then since the linear algebraic
$\Hg(V)$ possesses a reductive real form, it is reductive (as an
algebraic group), and then $\MT(V)$ is reductive as well.
\Qed
\enddemo

\proclaim{2.5.1. Corollary}
 When $V$ is a polarizable rational Hodge structure, $\Hg(V)$ is the
largest subgroup of $\GL(V)$ fixing all Hodge vectors in all $V^{\tensor
r}\tensor (V\dual)^{\tensor s}$, for $r,s\in\NN$.
\endproclaim

\demo{Proof}
 A reductive subgroup of $\GL(V)$ is characterized by its invariants in
the extended tensor algebra of ~$V$ (compare \cite{B.27}~Prop.3.1).
\Qed
\enddemo

 Once we know that $\Hg(V)$ and $\MT(V)$ are reductive, the next corollary
follows from Lemma~2.3(iii).

\proclaim{2.5.2. Corollary}
 Let $V$ be  polarizable rational Hodge structure.  Then $\Hg(V)$ is
semisimple if and only if the center of $\MT(V)$ is $\Gm$, i.e., consists
only of scalars.
\endproclaim

 Now we turn specifically to the Hodge structures of complex abelian
varieties.  When combined with Proposition~2.4, the following lemma says
that the endomorphism algebra of a complex abelian variety may be
identified with the rational Hodge structure endomorphisms of its
associated rational Hodge structure.

\proclaim{2.6. Lemma}
 Let $A$ be a complex abelian variety.  Then
$$
 \EndoA \isom \End_{\MT(A)}(H_1(A,\QQ)) = \End_{\Hg(A)}(H_1(A,\QQ)).
$$
\endproclaim

\demo{Proof}
 Let $\dim A =g$.  Then $A(\CC) = \CC^g/L$ for some lattice $L$, and let
$W=W_\QQ := L\tensor\QQ$.  Then we can identify $W \simeq
H_1(A,\QQ)$ and identify the universal covering space $\CC^g$ of $A(\CC)$
as the real $2g$-dimensional space $\Wr = W\tensor \RR$ together with the
induced complex structure represented as a homomorphism
$h:\SS\to\GL(\Wr)$.  Now an element of $\EndoA$ is characterized by
firstly being a $\QQ$-linear endomorphism of $W$ and secondly being a
complex-linear endomorphism of $\Wr$, which precisely means that it
commutes with $h(\SS(\RR))$.  But a $\QQ$-linear endomorphism of $W$ that
commutes with $h(\SS(\RR))$ must commute with all of $\MT(A)$ acting on
{}~$W$.
\Qed
\enddemo

\proclaim{2.7. Proposition {\rm (\cite{B.17})}}
 Let $A$ be a complex abelian variety.  If $\EndoA$ is a simple
$\QQ$-algebra with center $\QQ$, then $\Hg(A)$ is simple.
\endproclaim

\demo\nofrills
 Recall that $\EndoA$ is a simple $\QQ$-algebra precisely when $A$ is
simple.  Then following \cite{B.142}~p.66, the idea of the proof
is that the center, say $\frak c$, of $\mt(A)$ contains $\QQ\cdot\Id$ and
is contained in the center of $\EndoA$.  So if the center of $\EndoA$ is
$\QQ\cdot \Id$ then $\frak c =\QQ\cdot \Id$, and $\hg(A)$ is semisimple,
and for simple $A$, it is simple.
\enddemo

 As we have already begun to see, most of what can be said about the
structure and classification of the Hodge and Mumford-Tate groups of
abelian varieties is a consequence of the presence and properties of a
polarization.  The most fundamental fact is the following.

\proclaim{2.8. Lemma}
 Let $A$ be a complex abelian variety, and let $[E]$ be a polarization of
$A$ represented by the Riemann form~$E$.  Further, let $W=W_\QQ
=H_1(A,\QQ)$.  Then $E$ is a skew-symmetric bilinear form on $W$, and
there are natural representations
$$
 \Hg(A) \hra \Sp(W,E) \qquad \text{and}\qquad \MT(A) \hra \GSp(W,E) .
$$
\endproclaim

\demo{proof}
 This follows from the observation that the Riemann form $E$ is a
polarization on the rational Hodge structure $W$.  First, $E(h(i)u, h(i)v)
= E(u,v)$, where $u,v\in\Wr$ and $h:\SS\to\GL(\Wr)$ represents the complex
structure.  Thus if we write $h(s) = a\Id + b h(i)$, then $E(h(s)u,h(s)v)
= |a +bi|^2 E(u,v)$.  Therefore $h(\SS) \hra \GSp(\Wr,E)$, and then by
taking the Zariski-closure over~$\QQ$ we find $\MT(A) \hra \GSp(W,E)$.
\Qed
\enddemo

 The following criterion for the semisimplicity of the Hodge group is
linked to whether there are any simple components of Weil type,
see~1.13.2; cf.~also the discussions in sections four and five, below.

\proclaim{2.9. Proposition \rm{(\cite{B.118})}}
 Suppose $A$ is an abelian variety defined over $\CC$.  Then the Hodge
group of $A$ is not semisimple if and only if for some simple component
$B$ of $A$ the center of $\Endo(B)$ is a CM-field ~$K$ such that $(B,K)$,
with $K$ embedded in $\Endo(B)$ by the identity map, is not of Weil type.
\endproclaim

 The next proposition is that the real Lie groups $\Hg(A,\RR)$ and
$\MT(A,\RR)$ are of Hermitian type (see~1.6.5).

\proclaim{2.10. Proposition \rm{(\cite{B.77})}}
 Let $A$ be a complex abelian variety.  Then $\Hg(A,\RR)$ is of Hermitian
type.  Further, letting $K = K_\RR$ denote the centralizer of $h(i)$, or
equivalently of $h(\U(1))$, then the topologically connected component
$K^+$ of $K$ is a maximal compact subgroup of the topologically connected
component $\Hg(A,\RR)^+$ of $\Hg(A,\RR)$, and the quotient
$\Hg(A,\RR)^+/K^+$ is a Hermitian symmetric space of noncompact type,
i.e., a bounded symmetric domain.
\endproclaim

\demo{Proof {\rm (after \cite{B.25}, but see also \cite{B.103} and
\cite{B.51}~Ch.VIII)}}
 First, the connected center of $\Hg(A,\RR)$ is compact, since the group
itself is generated by compact subgroups, namely the
$\Aut_\QQ(\CC)$-conjugates of $h(\U(1))$.  From the proof that $\Hg(A)$ is
reductive it follows that $\Ad(h(i))$ defines a Cartan involution on
$\Hg(A,\RR)$, and thus $\ad(h(i))$ is a Cartan involution on $\frak g =
\hg(A,\RR)_{\text{ss}}$, the semisimple part of the reductive Lie
algebra $\hg(A,\RR)$.  Let $\frak k + \frak p$ be the corresponding
Cartan decomposition.  Then the restriction $\tilde J$ of $h(i)$
to $\frak p$ is a derivation of $\frak g$.  Since $\frak g$ is
semisimple, there therefore exists $H_0 \in \frak g$ such that $\tilde J
= \ad(H_0)$..  And since $\tilde J$ commutes with the Cartan decomposition
$\ad(h(i))$, we have that $H_0$ is in $\frak k$, and in the center of
$\frak k$, as required.
\Qed
\enddemo

  The following result points to how the complexified Lie algebras of the
Hodge and Mumford-Tate groups of a complex abelian variety fit into the
general classification of complex semisimple Lie algebras.  This
formulation of the result follows \cite{B.142}.

\proclaim{2.11. Theorem {\rm (\cite{B.25}~\S1, but see also
\cite{B.108}~\S3 and Appendix and \cite{B.141})}}
 Let $A$ be a complex abelian variety, let $\frak g$ be a simple factor of
the complex semisimple Lie algebra $\mt(A,\CC)_{\text{ss}}$, and let $r$
denote the rank of~$\frak g$. Further, let $W = H_1(A,\QQ)$, and let $V
\subset W_\CC$ be an irreducible subrepresentation for the action of
$\frak g$ on $W_\CC$.  Then $\frak g$ and $V$ must be one of the
following:
\roster
\item"(A)" $\frak g \simeq \frak{sl}(r+1)$ and $V$ is equivalent to the
$s$-th exterior power of the standard representation of dimension $r+1$,
for some $1\le s \le r$.
\item"(B)" $\frak g \simeq \frak{so}(2r+1)$, and $V$ is equivalent to the
spin representation of dimension ~$2^r$.
\item"(C)" $\frak g \simeq \frak{sp}(2r)$, and $V$ is equivalent to the
standard representation of dimension~$2r$.
\item"(D)" $\frak g \simeq \frak{so}(2r)$, and $V$ is equivalent to the
standard representation of dimension~$2r$, or to one of the two
half-spin representations of dimension~$2^{r-1}$.
\endroster
\endproclaim

 We do not give the proof here, but we just mention that the proof depends
essentially on the symplectic representation in ~2.8.

\proclaim{2.12. Proposition {\rm(\cite{B.78})}}
 A complex abelian variety is of CM-type if and only if $\Hg(A)$ is an
algebraic torus.
\endproclaim

\demo{Proof}
 Suppose first that $A$ is of CM-type.  Then $\EndoA$ contains a
commutative semisimple $\QQ$-algebra of dimension $2\dim A$ over $\QQ$.
{}From Lemma~2.6 it follows that $\Hg(A)$ commutes with a maximal
commutative semisimple subalgebra $R'\subset \End(W)$, where
$W=H_1(A,\QQ)$.  Therefore $\Hg(A)$ is contained in the units of $R'$ and
thus must be an algebraic torus.  Conversely, if $\Hg(A)$ is an algebraic
torus, then it is diagonalizable over $\CC$.  Therefore its centralizer in
$\End(W) \tensor \CC$, and thus its centralizer in $\End(W)$, contains a
maximal commutative semisimple subalgebra $R'\subset \End(W)$.  But then
$[R':\QQ] = \dim W = 2\dim A$ and $R'\subset \EndoA$, so $A$ is of
CM-type.
\Qed
\enddemo

\remark{Remark}
 If $A$ is an abelian variety with complex multiplication by $K$, and
$K_0$
is the maximal totally real subfield of $K$, then a more precise statement
is that
$$
 \Hg(A) \subseteq \Ker\{\Res_{K/\QQ} \Gm_{/K} \lra \Res_{K_0/\QQ}
\Gm_{/K_0} \} ,
$$
 where the arrow is induced by the norm map from $K$ to $K_0$.  The
arguement above shows that $\Hg(A) \subseteq \Res_{K/\QQ} \Gm_{/K}$.  Then
observe that $h_1(\U(1))$ is contained in the real points of the indicated
kernel, and recall that $\Hg(A)$ is the smallest algebraic subgroup
defined over $\QQ$ which over $\RR$ contains $h_1(\U(1))$.
\endremark

\definition{2.13. Definition}
 Let $K$ be a CM field, $S$ a CM-type for $K$, and let $A$ be the
corresponding abelian variety (up to isogeny), as described in~1.5.4.
Then the CM-type $(K,S)$ or the abelian variety $A$ with that CM-type is
said to be \dfn{nondegenerate} if $\dim\Hg(A) = \dim A = \frac12 [K:\QQ]$.
\enddefinition

\subhead  The Lefschetz group of an abelian variety
\endsubhead
 The Lefschetz group of an abelian variety was first studied by Ribet
\cite{B.94} and further investigated by Murty \cite{B.82}~\S2.

\definition{2.14. Definition}
 Let $A$ be a complex abelian variety, let $W= H_1(A,\QQ)$, and let $[E]$
be a polarization of $A$ represented by the Riemann form~$E$.  The
\dfn{Lefschetz group} of $A$ is the connected component of the identity in
the centralizer of $\EndoA$ in $\Sp(W,E)$,
$$
 \Lf(A) := \{g\in\Sp(W,E) : g\circ \phi = \phi\circ g \enspace \text{for
all } \phi\in\EndoA\}^\circ .
$$
\enddefinition

 In this definition $\Lf(A)$ appears to depend on the choice of
polarization, but if $[E']$ is another polarization,
then there is an element $\psi \in \EndoA$ and positive $m\in \ZZ$ such
that $mE' = E\psi$.  To see this, we take the point of view that $E$ and
$E'$ define isogenies, say $\phi$ and $\phi'$ respectively, from $A$ to
its dual~$A\dual$.  Then we can take $\psi = \phi\dual\circ\phi'$ and
$m=\deg \phi$.  Thus $\Lf(A)$ does not in fact depend on the choice of
polarization. Furthermore, it is clear that $\Lf(A)$ is an algebraic group
defined over $\QQ$, and
$$
 \Hg(A) \subseteq \Lf(A) .
$$
 The Lefschetz groups also has the following nice multiplicative property,
that the Hodge and Mumford-Tate groups in general do not.

\proclaim{2.15. Lemma {\rm (\cite{B.82}~Lem.2.1)}}
 If $A$ is isogenous to a product $B_1^{n_1} \times \dots \times
B_r^{n_r}$, with the $B_i$ simple and non-isogenous, then
$$
 \Lf(A) \simeq \Lf(B_1) \times \dots \times \Lf(B_r) .
$$
\endproclaim

\demo{Proof}
 First let $A_i = B_i^{n_i}$, for $1\le i \le r$, and choose polarizations
$[E_i]$ of $A_i$.  Then $[E_1 \oplus \dots \oplus E_r]$ is a polarization
of $A$, since $W = H_1(A,\QQ) \simeq \bigoplus_{i=1}^r H_1(A_i,\QQ)$.
Further, since the $B_i$ are non-isogenous, $\Hom(A_i, A_j) =0$ for $i\ne
j$, whence $\EndoA = \prod_{i=1}^r \Endo(A_i)$.  Therefore any
automorphism of $W$ that commutes with the action of $\EndoA$ must
preserve each $ H_1(A_i,\QQ)$, and thus $\Lf(A) \simeq \Lf(A_1)\times
\dots
\times \Lf(A_r)$.

 Now fix $i$, and let $B=B_i$ and $A=B^n$.  Then if $[E]$ is a
polarization of $B$, then $[E\oplus \dots \oplus E]$ is a polarization of
$A$.  Further, $\EndoA \simeq M_n(\Endo(B))$, so the centralizer of
$\EndoA$ in $\Sp(H_1(A,\QQ),E\oplus \dots \oplus E)$ can be identified
with the centralizer of $\Endo(B)$ in $\Sp(H_1(B,\QQ),E)$.  Therefore
$\Lf(A)\simeq \Lf(B)$.
\Qed
\enddemo

 For later reference we state some variants of ``Goursat's Lemma'' that
turn out to be useful, especially when extending results from simple
abelian varieties to products of abelian varieties.  The formulations
given below come from \cite{B.90}, \cite{B.83} and
\cite{B.75}.

\proclaim{2.16. Proposition {\rm(Goursat's Lemma)}}
\roster
\item
 Let $G$ and $G'$ be groups and suppose $H$ is a subgroup of $G\times G'$
for which the projections $p:H\to G$ and $p': H\to G'$ are surjective.
Let $N$ be the kernel of $p'$ and let $N'$ be the kernel of $p$.  Then $N$
is a normal subgroup of $G$ and $N'$ is a normal subgroup of $G'$, and the
image of $H$ in $G/N \times G'/N'$ is the graph of an isomorphism $G/N
\simeq G'/N'$.
\item
 Let $V_1$ and $V_2$ be two finite-dimen\-sion\-al complex vector spaces.
Let $\frak s_1, \frak s_2$ be simple complex Lie subalgebras of
$\frak{gl}(V_1)$, $\frak{gl}(V_2)$ respectively, of type A, B or C.  Let
$\frak s$ be a Lie subalgebra of $\frak s_1 \times \frak s_2$ whose
projection  to each factor is surjective.  Then either $\frak s = \frak
s_1 \times \frak s_2$ or $\frak s$ is the graph of an isomorphism $\frak
s_1 \simeq \frak s_2$ induced by an $\frak s$-module isomorphism $V_2
\simeq V_1$ or $V_2 \simeq V_1\dual$.
\item
 Let $\frak s_1, \dots , \frak s_d$ be simple finite-dimen\-sion\-al Lie
algebras and let $\frak g$ be a subalgebra of the product $\frak s_1
\times \dots \times \frak s_d$.  Assume that for $1\le i\le d$ the
projection $\frak g \to \frak s_i$ is surjective, and that whenever $1\le
i < j\le d$ the projection of $\frak g$ onto $\frak s_i \times \frak s_j$
is surjective.  Then $\frak g = \frak s_1 \times \dots \times \frak s_d$.
\item
 Let $I$ be a finite set and for each $\sigma\in I$, let $\frak s_\sigma$
be a finite-dimensional complex simple Lie algebra.  Let $\frak g, \frak
h$ be two algebras such that
\itemitem"{(a)}" $\frak g \subseteq \frak h$.
\itemitem"{(b)}" $\frak h$ is a subalgebra of $\prod_{\sigma\i I} \frak
s_\sigma$ such that the projection to each $\frak s_\sigma$ is surjective.
\itemitem"{(c)}" $\frak g, \frak h$ have equal images on $\frak s_\sigma
\times \frak s_\tau$ for all pairs $(\sigma, \tau) \in I\times I$, $\sigma
\ne \tau$.
\item""  Then $\frak g = \frak h = \prod_{\sigma \in J} \frak s_j$ for
some subset $J\subseteq I$.
\item
 Let $V_1, \dots , V_n$ be finite-dimensional vector spaces over an
algebraically closed field of characteristic zero, and let $\frak g$ be a
semisimple Lie subalgebra of $\End(V_1) \times \dots \times \End(V_n)$.
For $1\le i \le n$ let $\frak g_i\subseteq \End(V_i)$ be the projection of
$\frak g$ onto the $i$-th factor.  Assume that $\frak g_i$ is nonzero and
simple for all ~$i$.  Then for any simple Lie algebra $\frak h$ let
$I(\frak h) \subset \{1,\dots, n\}$ be the set of indices for which $\frak
g_i \simeq \frak h$.  Assume that for any $\frak h$ with $\#I(\frak h) >
1$ the following conditions are satisfied:
\itemitem"(a)" All automorphisms of $\frak h$ are inner.
\itemitem"(b)" For $i\in I(\frak h)$ the representations $V_i$ are all
isomorphic.
\itemitem"(c)" $\End_{\frak g}(\bigoplus_{i\in I(\frak h)} V_i) =
\prod_{i\in I(\frak h)} \End_{\frak g_i}(V_i)$.
\item""  Then $\frak g \simto \frak g_1 \times \dots \times \frak g_n$.
\endroster
\endproclaim

\head 3. Products of Elliptic curves
\endhead

 Tate seems to be the first to have checked the (usual) Hodge conjecture
for powers $E^n$ of an elliptic curve, see \cite{B.130},
\cite{B.43}~\S3, but he never published his proof.  In
\cite{B.80} Murasaki showed the $\Hdg^p(E^n) = \Div^p(E^n)$ for
all $p$ by exhibiting explicit differential forms  that give a basis for
$\Hdg^1(E^n)$ and then carrying out explicit computations with them.  In a
different direction, Imai \cite{B.58} showed that when $E_1, \dots,
E_n$ are pairwise non-isogenous elliptic curves, then $\Hg(E_1\times\dots
\times E_n) \simeq \Hg(E_1)\times \dots \times \Hg(E_n)$.  A unified
approach to computing the Hodge and Mumford-Tate groups, and verifying the
Hodge conjecture, for arbitrary products of elliptic curves can be found
in \cite{B.84}.  Since Murty's approach provides a nice example of
how the Hodge and Mumford-Tate groups can be used to verify the usual
Hodge conjecture, we summarize his exposition here.

\proclaim{Theorem}
 Let $A= E_1^{n_1} \times \dots \times E_r^{n_r}$, where the $E_i$ are
pairwise non-isogenous elliptic curves.  Then
\roster
\item $\Hg(A) = \Hg(E_1) \times \dots \times \Hg(E_r)$.
\item $\Hdg(A) = \Hdg(E_1^{n_1}) \tensor \dots \tensor \Hdg(E_r^{n_r}) =
\Div(A)$.
\endroster
\endproclaim

\demo{Proof {\rm (after \cite{B.84})}}
 Let $E$ be an elliptic curve.  The cases where $E$ has or does not have
complex multiplication have to be handled separately.  If $E$ has complex
multiplication, then $\Endo(E) =:K$ is an imaginary quadratic field, and
(the proof of) Proposition~2.12 shows that $\MT(E,\QQ)\subseteq K^\times$,
as algebraic groups over $\QQ$.  Since the two-dimensional $\SS(\RR)
\subset \MT(E,\RR)$, we see that $\MT(E) = \Res_{K/\QQ}(\Gm_{/K})$.  When
$E$ does not have complex multiplication, then $\hg(E)$ is a simple
subalgebra of $\frak{sl}(V)$ which is already simple, so $\hg(E) =
\frak{sl}(V)$ and thus $\Hg(E) =\SL(V)$ and $\MT(E) =\GL(V)$.

 Next consider $A= E^n$.  Then
$$ \align
 \Hdg(A) &= H^*(E^n,\QQ)^{\Hg(A)} \\
 &= \twedge^*(H^1(E,\QQ) \oplus \dots \oplus H^1(E,\QQ))^{\Hg(A)} \\
 &= \bigoplus (H^1(E,\QQ) \tensor \dots \tensor H^1(E,\QQ))^{\Hg(E)} .
\endalign
$$
 Now if $E$ has complex multiplication then $\alpha \in K^\times = \MT(E)$
acts on $H^1(E,\QQ)\tensor_\QQ\CC \simeq \CC \oplus \CC$ by $\alpha(z,w) =
(\alpha z, \bar\alpha w)$.  Let $K^\times_1$ denote the elements of
$K^\times$ of norm~$1$.  Then
$$
(H^1(E,\QQ) \tensor \dots \tensor H^1(E,\QQ))^{\Hg(E)} \tensor_\QQ \CC =
\big( (H^1(E,\QQ)\tensor \CC) \tensor \dots \tensor (H^1(E,\QQ)\tensor
\CC) \big)^{K^\times_1} ,
$$
 in which any invariant class arises as a combination of products of
elements of
$$
 \big( (H^1(E,\QQ)\tensor \CC) \tensor ((H^1(E,\QQ)\tensor \CC))
\big)^{K^\times_1} \subseteq (H^2(E\times E,\QQ)\tensor\CC)^{K^\times_1} .
$$
 Therefore the invariants are generated by those in $H^2(A,\QQ)$, which
means that $\Hdg(A) =\Div(A)$ in this case.

 Next suppose that $E$ does not have complex multiplication.  Then $\Hg(E)
= \SL(2)$ and acts on $H^1(E,\QQ)$ by the standard representation.  Now we
invoke the well-known fact that the tensor invariants of $\SL(2)$ are
generated by the determinant; see \cite{B.137} or \cite{B.10}~App.1 for
this.  Since the determinant is a
representation of degree~$2$ lying in $H^1(E,\QQ) \tensor H^1(E,\QQ)
\subset H^2(A,\QQ)$, again we find that all Hodge cycles of $A= E^n$ are
generated by divisors.

 Finally let $A= E_1^{n_1} \times \dots \times E_r^{n_r}$, where the $E_i$
are pairwise non-isogenous elliptic curves.  First suppose all the $E_i$
have complex multiplication by an imaginary quadratic field $K_i
=\Endo(E_i)$.  Since the field $\Endo(E_i)$ determines the isogeny class
of $E_i$, all the $K_i$ are distinct.  Then
$$
 \Hg(A) \subseteq K^\times_{1,1} \times \dots \times K^\times_{r,1},
$$
 and moreover from the definition, $\Hg(A)$ surjects onto each factor.
Therefore there is a surjection of character groups
$$
 \lambda: M\to \X(\Hg(A)) ,
$$
 where
$$
 M:= \X(K^\times_{1,1}) \oplus \dots \oplus \X(K^\times_{r,1}) .
$$
 Now to see that $\lambda$ is an isomorphism and prove the theorem in the
case where all the $E_i$ have complex multiplication, we observe that for
each $i$ the composition
$$
 \X(K^\times_{i,1}) \hra M \to \X(\Hg(A))
$$
 of $m\mapsto (0,\dots,0,m,0,\dots,0)$ with $\lambda$ is injective.  In
addition, all of these character groups are $\script G =
\Gal(\QQ^{\text{ab}}/\QQ)$-modules and the maps $\script G$-module maps. 
Then since the fields are distinct there is some $\sigma\in\script G$ that
acts as
$+1$ on $\X(K^\times_{1,1})$ and $-1$ on the other components.  Thus if $m
= (m_1,\dots ,m_r)$ is in the kernel of $\lambda$ then $\sigma m + m =
(2m_1, 0,\dots,0)$ must be as well.  Then the injectivity of the
composition above forces $m_1 =0$, and by induction the kernel of
$\lambda$ is zero.

 Next suppose none of the $E_i$ has complex multiplication.  Then we have
$$
 \hg(A) \subseteq \hg(E_1) \times \dots \times \hg(E_r) ,
$$
 and mapping surjectively onto each factor.  Then by Proposition~2.16.4,
if
it also maps surjectively onto each pair of factors, it is the entire
product.  But by Proposition~2.16.2, if $\hg(A)$ does not project onto
$\hg(E_i) \times \hg(E_j)$ for all pairs $i\ne j$, then it projects to the
graph of an isomorphism between them, which in turn could be used to
produce an isogeny between $E_i$ and $E_j$, contrary to assumption.

 Finally it remains to see that if $A$ is an abelian variety isogenous to
a product $B\times C$ with $\Hg(B)$ a torus and $\Hg(C)$ semisimple, then
$\Hg(A) = \Hg(B) \times \Hg(C)$.  However, this is a consequence of
Proposition~2.16.1.  This completes the proof of the theorem.
\Qed
\enddemo



\head 4. Abelian varieties of Weil or Fermat type
\endhead
  We have already defined abelian varieties of Weil type
(1.13.2), and an abelian variety of Fermat type is one which is isogenous
to a product of certain factors of the Jacobian variety of a Fermat curve
$x^m + y^m + z^m =0$ \cite{B.116}.  The important thing about these
examples, insofar as the Hodge $(p,p)$ conjecture goes, is that they
contain the only known examples where the conjecture has been verified for
abelian varieties $A$ for which $\Hdg(A) \ne \Div(A)$.  However, both
types also provide explicit examples of Hodge cycles that are not known to
be algebraic.  Indeed, as is well-known, Weil has suggested that a place
to look for a counterexample to the Hodge $(p,p)$ conjecture might be
among
what are now called abelian varieties of Weil type \cite{B.135}.

  We will begin by summarizing Shioda's results on abelian varieties of
Fermat type.  Then, since there is a nice presentation of the issues
concerning abelian varieties in \cite{B.35}, we will just
summarize the saliant points here.  Finally we will recall the work of
Schoen \cite{B.104} and van Geemen \cite{B.36} verifying
the Hodge conjecture for special four-dimensional abelian varieties.

\subhead Shioda's results on abelian varieties of Fermat type and
Jacobians of hyperelliptic curves {\rm (\cite{B.116})}
\endsubhead
  We will attempt to give a careful statement of the results, and refer
the reader to the original for the proofs.

\remark{4.1. Notation}
 Fix an integer $m>1$, and for $a\in\ZZ$ not congruent to zero modulo~$m$,
let $1\le \bar a \le m-1$ be the unique integer such that $\bar a \cong a
\pmod m$.  Let
$$\align
 \frak A_m^n &:= \{ \alpha =(a_0,\dots , a_{n+1}) : 1\le a_i \le m-1, \
\sum_{i=0}^{n+1} a_i \cong 0 \pmod m \} \\
 \frak B_m^n &:= \{\alpha \in \frak A_m^n : |t\cdot \alpha| = (n/2)+1
\quad \text{for all } t\in (\ZZ/m\ZZ)^\times\} ,
\endalign
$$
 where in the latter case $n$ must be even, and where for
$t\in(\ZZ/m\ZZ)^\times$ and $\alpha \in \frak A_m^n$,
$$
 t\cdot \alpha := (\ol{ta}_0, \dots , \ol{ta}_{n+1}), \qquad |\alpha| :=
\frac 1 m \sum_{i=0}^{n+1} a_i .
$$
 For $\alpha = (a_0,\dots , a_{r+1}) \in \frak A_m^r$ and $\beta =
(b_0,\dots,b_{s+1}) \in \frak A_m^s$ let
$$
 \alpha * \beta := (a_0,\dots , a_{r+1}, b_0,\dots,b_{s+1}) \in \frak
A_m^{r+s+2} .
$$

 Further, let
$$\alignat2
 M_m &:= \{ \xi=(x_1, \dots,x_{m-1};y) :
 \sum_{\nu=1}^{m-1} \ol{t\nu}\, x_\nu = my \enspace
  &&\text{for all }t\in(\ZZ/m\ZZ)^\times,
 \\  \vspace{-2\jot}
  &&& x_\nu,y\in\ZZ,\ x_\nu\ge 0,\ y>0\}
 \\   \vspace{2\jot}
 M_m(d) &:= \{ (x_1, \dots,x_{m-1};d) \in M_m\} .
\endalignat
$$
\endremark

\definition{4.2. Definition}
 An element $\xi\in M_m$ is said to be \dfn{indecomposable} if $\xi \ne
\xi' +\xi''$ for any $\xi',\xi''\in M_m$.  An element $\xi\in M_m$ is
called \dfn{quasi-decomposable} if there exists $\eta\in M_m(1)$ such that
$\xi +\eta = \xi' +\xi''$ for some $\xi',\xi'' \in M_m$ with $\xi', \xi''
\ne \xi$.
\enddefinition

    It is easy to see that the set $M_m$ is an additive semigroup with
only a finite number of indecomposable elements.

\definition{4.3. Definition}
 Let $X_m : x^m +y^m +z^m=0$ denote the Fermat curve of degree~$m$, let
$J(X_m)$ be its Jacobian, and let $\frak S_m = (\ZZ/m\ZZ)^\times \bs \frak
A_m^1$ be the orbit space.  Then there is an isogeny
$$
 \pi :J(X_m) \to \prod_{S\in \frak S_m} A_S
$$
 where
\roster
\item"(i)" $A_S$ is an abelian variety of dimension $\phi(m')/2$ admitting
complex multiplication by $\QQ(\zeta_{m'})$, where $\zeta_{m'} = e^{2\pi
i/{m'}}$  and $m' = m/ \operatorname{gcd}(a,b,c)$ for $(a,b,c)\in \frak
A_m^1$ belonging to the orbit~$S$.
\item"(ii)" $H^1(A_S,\CC)$ has the eigenspace decomposition
$$
 H^1(A_S,\CC) = \bigoplus_{\alpha\in S} W(\alpha)
$$
 for the complex multiplication of~(i), where $\dim W(\alpha) =1$ and such
that, if $\pi_S : J(X_m) \to A_S$ is the composition of $\pi$ and the
projection to the $S$ factor, then $\pi_S^* W(\alpha) = U(\alpha)$ with
$\alpha \in S$.  Here $U(\alpha)$ is defined by
$$
 H^{1,0}(J(X_m)) = \bigoplus \Sb \alpha\in \frak A_m^1 \\ |\alpha|=1
\endSb U(\alpha) , \qquad  H^{0,1}(J(X_m)) = \bigoplus \Sb \alpha\in \frak
A_m^1 \\ |\alpha|=2 \endSb U(\alpha)
$$
 and $U(\alpha)$ is one-dimensional \cite{B.41} \cite{B.59}.
\endroster
 Then an abelian variety will be said to be \dfn{of Fermat type of
degree~$m$} if it is isogenous to a product of a finite number of factors
$A_S$ satisfying~(i) and~(ii) as above.
\enddefinition

 Thus an abelian variety of Fermat type of degree~$m$ is given by
$$
 A= \prod_{i=1}^k A_{S_i}
\tag 4.3
$$
 with $S_1,\dots,S_k \in \frak S_m$ not necessarily distinct.  With this
notation we can now state Shioda's results on the Hodge $(p,p)$ conjecture
for these abelian varieties.

\proclaim{4.4. Theorem {\rm (\cite{B.116}~Thm.4.3)}}
 Let $A$ be an abelian variety of Fermat type of degree~$m$.  Assume that
for any set $\{\alpha_1,\dots,\alpha_{2d}\}$ of distinct elements in the
disjoint union of the $S_i$, $i=1,\dots,k$, such that $\alpha_1 * \dots *
\alpha_{2d} = \gamma \in \frak B_m^{6d-2}$ there exists some
$\beta_1,\dots , \beta_l \in \frak B_m^0 \cup \frak B_m^2 \cup (\frak
B_m^4\cap \frak A_m^1 *\frak A_m^1)$ such that $\beta_1 * \dots * \beta_l$
coincides with $\gamma$ up to permutation.  Then the Hodge $(p,p)$
conjecture is true for $A$ in codimension~$d$.
\endproclaim

\proclaim{4.5. Theorem {\rm (\cite{B.116}~Thm.4.4)}}
 If every decomposable element of $M_m(y)$ with $3\le y \le 3d$, if any,
is quasi-decomposable, then the Hodge $(p,p)$ conjecture is true in
codimension ~$d$ for all abelian varieties of Fermat type of degree~$m$.
In particular, if $m$ is a prime or $m\le 20$, then the Hodge conjecture
is
true in any codimension for all abelian varieties of Fermat type of
degree~$m$.
\endproclaim

\proclaim{4.6. Theorem {\rm (\cite{B.116}~Thm.5.6)}}
 For any given $d\ge 2$ there exists some abelian variety of Fermat type
$A$ such that the Hodge ring $\Hdg^*(A)$ is not generated by
$\sum_{r=1}^{d-1} \Hdg^r(A)$.
\endproclaim

\proclaim{4.7. Theorem {\rm (\cite{B.116}~Thm.5.3)}}
 Let $C_m : y^2 = x^m -1$ be the hyperelliptic curve of genus
$g=[(m-1)/2]$, and let $J(C_m)$ be its Jacobian.  If $m>2$ is a prime
number, then the Hodge ring $\Hdg^*(J(C_m))$ is generated by
$\Hdg^1(J(C_m))$, i.e., $\Hdg^*(J(C_m)) = \Div^*(J(C_m))$ and the Hodge
conjecture is true for $J(C_m)$.  The same result also holds for arbitrary
powers of $J(C_m)$.
\endproclaim

\proclaim{4.8. Theorem {\rm (\cite{B.116}~Thm.5.4)}}
 For any odd $m\ge 3$ the Hodge conjecture is true for $J(C_m)$ in
codimension~$2$.
\endproclaim

\subhead Abelian varieties of Weil type
\endsubhead
  From 1.13.2, a Weil abelian variety of dimension $g=2n$ is an abelian
variety $A$ together with an imaginary quadratic field $K$ embedded in
$\EndoA$ such that the action of $\alpha \in K$ has the eigenvalues
$\alpha$ and $\bar \alpha$ with equal multiplicity~$n$.  Here we briefly
review some of the important points about Weil abelian varieties, mainly
following the exposition of \cite{B.35}, to which we refer the
reader for more detail; the original source is \cite{B.135}.

 If $(A,K)$ is an abelian variety of Weil type, and $K=\QQ(\sqrt{-d})$,
then it is possible to choose a polarization $[E]$ on $A$ normalized so
that $(\sqrt{-d})^* E = d E$, where here we are viewing $\sqrt{-d}
\in \scO_K \hra \End(A)$ as an endomorphism of $A$, so
$(\sqrt{-d})^*$ and $(\sqrt{-d})_*$ denote the induced pullback and
push-out maps, respectively.   Hereafter when we speak of a polarized
abelian variety of Weil type we will assume that its polarization is
normalized this way.  Let $W=H_1(A,\QQ)$.  Then the $K$-valued Hermitian
form $H: W\times W \to K$ associated to $E$ is given by
$$
 H(u,v) := E(u,(\sqrt{-d})_* v) + (\sqrt{-d})E(u,v) .
$$

\proclaim{4.9. Theorem {\rm(\cite{B.135})}}
 The Hodge group of a general polarized abelian variety $(A,K,E)$ of Weil
type is $\SU(W,H)$ (as an algebraic group over~$\QQ$).
\endproclaim

   It can be shown that polarized abelian varieties of Weil type of
dimension~$2n$ are parameterized by an $n^2$-dimensional space which can
be described as the bounded symmetric domain assoicated to
$\SU(W,H;\RR)\simeq \SU(n,n)$. Thus the word ``general'' in the statement
of the theorem refers to a general point in this parameter space,
analogously to the usage in~1.13.8.

\definition{4.10. Definition}
 Let $(A,K)$ be an abelian variety of Weil type and dimension~$2n$.
Then $V = H^1(A,\QQ)$ has the structure of a vector space over~$K$.  The
space of \dfn{Weil-Hodge cycles} on $A$ is the subset of $H^{2n}(A,\QQ)$
$$
 \Weil(A) := \twedge_K^{2n} H^1(A,\QQ)  ,
$$
 where $\twedge_K^{2n}$ signifies the $2n^{\text{th}}$-exterior power of
$H^1(A,\QQ)$ as a $K$@-vector space.
\enddefinition

\proclaim{Lemma}
 $\dim_\QQ \Weil(A) =2$, and $\Weil(A) \subset \Hdg^n(A)$.
\endproclaim

\proclaim{4.11. Theorem {\rm(\cite{B.135})}}
 Let $(A,K)$ be an abelian variety of Weil type of dimension $g=2n$.  If
the Hodge groups of $A$ is $\Hg(A) = \SU(W,H)$, then
$$
 \dim \Hdg^p(A) = \cases 1, & p\ne n, \\ 3, &p=n \endcases
$$
 and $\Hdg^n(A) = \Div^n(A) \oplus \Weil(A)$.
\endproclaim

  Thus abelian varieties of Weil type provide examples of Hodge cycles
which do not arise from products of those in codimension one.

  In the following the determinant of the Hermitian form $H$ is
well-defined as an element of $\QQ^\times$ modulo the subgroup of norms
from $K^\times$.

\proclaim{4.12. Theorem {\rm(\cite{B.104})}}
 The Hodge $(2,2)$ conjecture is true for a general abelian variety of
Weil type $(A,K)$ of dimension~$4$ with $K=\QQ(\sqrt{-3})$ or $K=\QQ(i)$,
when $\det H =1$.
\endproclaim

A different proof for the case where $K=\QQ(i)$ and $\det H =1$ can be
found in \cite{B.36}.  In \cite{B.35}~7.3 it is
pointed out that Schoen's methods together with a result of \cite{B.32}
imply the Hodge $(3,3)$ conjecture for a general $6$-dimensional
abelian variety of Weil type with $K=\QQ(\sqrt{-3})$ and $\det H =1$.

\bigpagebreak

  Recently Moonen and Zarhin \cite{B.76} have considered
the extent to which Weil's construction of exceptional Hodge classes can
be generalized.  For $K\hra \EndoA$ a subfield and $r= 2\dim A/[K:\QQ]$,
let
$$
 \Weil_K(A) := \twedge_K^r H^1(A,\QQ) ,
$$
 and call this the space of \dfn{Weil classes with respect to~$K$.}  Let
$V = H^1(A,\QQ)$.  Then $K$ acts on $V$ and there is a decomposition
$$
 V\tensor\CC = \bigoplus_{\sigma\in\Hom(F,\CC)} V_{\CC,\sigma} =
\bigoplus_{\sigma\in\Hom(F,\CC)} (V_{\CC,\sigma}^{1,0} \oplus
V_{\CC,\sigma}^{0,1}).
$$

\proclaim{4.13. Proposition {\rm(\cite{B.76})}}
  Let $A$ be a complex abelian variety, $K\hra \EndoA$ a subfield, and $r=
2\dim A/[K:\QQ]$.  With the notation above,
\roster
\item If $\dim V_{\CC,\sigma}^{1,0} = V_{\CC,\ol\sigma}^{1,0}$ for all
$\sigma \in \Hom(F,\CC)$, where $\ol\sigma$ denotes the complex conjugate
of $\sigma$, then $\Weil_K(A)$ consists entirely of Hodge classes, i.e.,
$\Weil_K(A) \subset \Hdg^r(A)$; if $\dim V_{\CC,\sigma}^{1,0} \neq
V_{\CC,\ol\sigma}^{1,0}$ for some $\sigma \in \Hom(F,\CC)$, then the zero
class is the only Hodge class in~$\Weil_K(A)$.
\item Suppose $A$ is isogenous to a power $B^m$ of a simple abelian
variety~$B$, and suppose $K\hra \EndoA$ is a subfield such that
$\Weil_K(A)$ consists of Hodge classes.  Let $F$ be the center of
$\Endo(B)$, and $F_0$ the maximal totally real subfield of~$F$.  Then
either $\Weil_K(A) \subset \Div(A)$, or all nonzero classes in
$\Weil_K(A)$ are exceptional; this last possibility occurs precisely in
the following cases:
\itemitem"(a)" $B$ is of type~III, $m=1$ and $K\subsetneqq F$,
\itemitem"(b)" $B$ is of type~III, $m\ge 2$ and $2m[F:\QQ]/[K:\QQ]$ is
odd,
\itemitem"(c)" $B$ is of type~IV, $(\dim_F(\Endo(B))^{1/2} =1$, $m=1$ and
$K\subsetneqq F_0$,
\itemitem"(d)" $B$ is of type~IV with  $(\dim_F(\Endo(B))^{1/2} \ge 2$ or
$m \ge 2$ and the map
$$
 \Lie(\U_F(1)) \hra \End_K(V) @> \operatorname{Tr} >> K
$$
\item"" is nonzero.  Here $\U_F(1) = \Res_{F/\QQ}\Gm_F \cap \U(1)$
(see~1.5.2).
\endroster
\endproclaim



\head 5. Simple abelian fourfolds
\endhead
  Since the Hodge $(p,p)$ conjecture is true for any smooth projective
complex threefold, it might seem that fourfolds would be the next case to
attack.  However, as we've just seen, simple abelian fourfolds provide the
first examples of abelian varieties of Weil type, for which the Hodge
conjecture is mostly not known, and they also provide the first examples
of abelian varieties of type~(III) in the Albert classification, as well
as the first examples of abelian varieties that are not characterized by
their endomorphism rings \cite{B.78}~\S4.

  Recently Moonen and Zarhin \cite{B.75} have analyzed
the Hodge structures of simple abelian fourfolds, and their work makes an
instructive example.  The main result is the following, where a subalgebra
of $\EndoA$ is said to be stable under all Rosati involutions if for every
polarization it is stable under the associated Rosati involution.  Recall
that an exceptional Hodge class is one which is not accounted for by
linear combinations of intersections of divisors.

\proclaim{5.1. Theorem {\rm (\cite{B.75}~Thm.2.4)}}
 When $A$ is a simple abelian fourfold, then $A$ supports exceptional
Hodge classes if and only if $\EndoA$ contains an imaginary quadratic
field $K$ which is stable under all Rosati involutions and such that with
the induced action of $K \subseteq \EndoA$ the complex Lie algebra
$\Lie(A,\CC)$ of $A$ becomes a free $K\tensor\CC$-module.
\endproclaim

 In the situation of the theorem, $\Lie(A)$ becomes a free $K\tensor
\CC$-module if and only if $\alpha \in K$ acts as $\alpha$ and as
$\bar\alpha$ with equal multiplicity~$2$.  Thus a less precise but more
simply stated corollary would be the following.

\proclaim{Corollary}
 When $A$ is a simple abelian fourfold, if $\Hdg^2(A) \ne \Div^2(A)$ then
$A$ must be an abelian variety of Weil type.
\endproclaim

 For the abelian fourfolds $A$ for which  $\Hdg^2(A) \ne \Div^2(A)$ Moonen
and Zarhin prove the following theorem.

\proclaim{5.2. Theorem {\rm (\cite{B.75}~Thm.2.12)}}
 Let $A$ be a simple abelian fourfold, and let
$$
 \V(A) := \sum_K \left( \twedge^4_K H^1(A,\QQ) \right) ,
$$
 where the sum runs over all imaginary quadratic subfields $K\subset
\EndoA$ that act on $A$ with multiplicities $(2,2)$.  Then $\Hdg^2(A) =
\Div^2(A) + \V(A)$.
\endproclaim

 This theorem should be compared with Theorem~4.11:  Theorem~5.1 applies
to fourfolds, whereas in Theorem~4.11 it was assumed that $\Hg(A) =
\SU(W,H)$, which is a sort of a generality assumption.

\remark{5.3. Remark: Earlier work of Tankeev}
   In 1978 and 1979 Tankeev published two papers \cite{B.122}
\cite{B.123} containing results about the Hodge structure and Hodge
conjecture for abelian fourfolds.  In particular, in \cite{B.123}~Thm.3.2
he proved that when $A$ is a simple abelian fourfold of type~(I)
or type~(II), then $\Hdg(A) = \Div(A)$.  First he showed that if the
center of $\EndoA$ is a product of totally real fields, then $\Hg(A)$ is
semisimple \cite{B.123}~Lemma~1.4, and then he proved the Hodge
conjecture for simple abelian fourfolds of type~(I) or~(II) by considering
the possible symplectic representations of the complexified Lie algebra
$\hg_\CC(A)$ and showing that in each case that its invariants in
$H^*(A,\CC)$ are generated by those of degree~$2$.  The earlier paper
\cite{B.122} considered the possible pairs $(\frak g, \rho)$, where
$\frak g$ is the semisimple part of $\hg_\CC(A)$ and $\rho:\frak g \to
\End_\CC(W_\CC)$ denotes its action on $W_\CC=H_1(A,\CC)$, under the
assumption that that there does not exist an abelian variety $A_0$ defined
over $\Qbar$ such that $A_0 \tensor_\Qbar \CC \simto A$; however, the
proof there contains some gaps.  Then he derived, under the same
assumption that the abelian fourfold $A$ cannot be defined over $\Qbar$,
that when $\EndoA$ is neither an imaginary quadratic field nor a definite
quaternion algebra (type~(III)) then $\Hdg(A) = \Div(A)$.
\endremark

\subhead \nofrills \endsubhead
 We now proceed to sketch the outline of the proof of Theorem~5.1.  One
direction is covered by the following more general result.

\proclaim{5.4. Theorem {\rm (\cite{B.75}~Thm.3.1)}}
 Let $A$ be a simple abelian variety, and assume that either
\roster
\item"(a)" $A$ is of type~\rom{(III)}, or
\item"(b)" $\EndoA$ is a CM-field $K$ which contains a CM-subfield $F$
such that the multiplicity with which $\alpha \in F$ acts as
$\sigma(\alpha)$ is the same for all $\sigma\in \Hom(F,\CC)$
\endroster
 Then $A$ supports exceptional Hodge classes
\endproclaim

\remark\nofrills
 The result for abelian varieties of type~(III) is due to Murty
\cite{B.82}, see 8.6 below, although Moonen and Zarhin give a
different proof.  For case ~(b) what they show is that except for the zero
element, $\twedge_F^m H^1(A,\QQ)$ consists entirely of exceptional Hodge
classes, where $m= 2\dim A /[F:\QQ]$.
\endremark

  The other direction of Theorem~5.2 is proved case by case, running
through the different possible endomorphism algebras (see~1.13.4).  It
turns out that except when $\EndoA = \QQ$, knowing $\EndoA$ together with
its action on $\Lie(A)$ suffices to determine $\hg(A)$, which in turn is
enough to identify the absence or presence of exceptional cycles.  We will
run through the results, giving only some comments on the ingredients of
the proofs.  As usual, $W := H_1(A,\QQ)$.

\proclaim{5.5. Type I(1)}
 Let $A$ be a simple abelian fourfold with $\EndoA =\QQ$.  Then the Lie
algebra $\hg(A)$ together with its representation on $W$ is isomorphic
over $\Qbar$ to one of
\roster
\item"(i)" $\frak{sp}_4$ with the standard representation, or
\item"(ii)" $\frak{sl}_2 \times \frak{sl}_2 \times \frak{sl}_2$ with the
tensor product of the standard representation of each of the three
factors.
\endroster
 Both possibilities occur, and in both cases $\Hdg(A) = \Div(A)$.
\endproclaim

\remark\nofrills
 On the assumption that $\hg$ is simple, Theorem~2.11 can be used to show
that case~(i) is the only possibility.  When $\hg = \frak h_1 \times \dots
\times \frak h_t$ is not simple, then $W= W_1\tensor \dots \tensor W_t$
and at least one $W_i$ must be $2$-dimensional.  Then $\frak h_i =
\frak{sl}_2$, and since the representation is symplectic, the complement
must be $\frak{so}_4 \simeq \frak{sl}_2\times \frak{sl}_2$.  In both cases
$(\twedge^4 W)^{\mt}$ (the subspace of $\mt$-invariants in $\twedge^4 W$)
is computed to be $1$-dimensional.
\endremark

\remark{5.6. Notation}
 Before treating the remaining type~(I) cases, suppose in general that
$\EndoA$ contains a totally real field $F$, and suppose a polarization
$[E]$ on $A$ is given.  Then there is a unique $F$-bilinear alternating
form $\psi: W\times W \to F$ whose trace
$\operatorname{Tr}_{F/\QQ}(\psi(u,v)) = E(u,v)$.  Then from the uniqueness
of $\psi$ it follows that $\hg(A)$ is contained in
$$
 \frak{sp}_F(W,\psi) := \{ \phi\in \End_F(W) : \psi(\phi(u),v) +
\psi(u,\phi(v)) =0 \enspace\text{for all } u,v \in W \} ,
$$
 regarded as a Lie algebra over $\QQ$.
\endremark

\proclaim{5.7. Type I(2)}
  Let $A$ be a simple abelian fourfold with $\EndoA = F$ a real quadratic
field.  Then in the notation above,  $\hg(A) \simeq \frak{sp}_F(W,\psi)$.
In particular, $\Hdg(A^n) = \Div(A^n)$ for all ~$n$.
\endproclaim

\demo{Sketch of proof}
 The representation $\hg_\CC$  on $W_\CC$ splits as a direct sum $W_\CC =
W_1 \oplus W_2$ with $\dim W_1 = \dim W_2 =4$.  Further, the restriction
of $\psi$ to $W_i$ is a nondegenerate skew-symmetric bilinear form
$\psi_i: W_i\times W_i \to \CC$, and $\hg_\CC \subseteq
\frak{sp}(W_1,\psi_1) \times \frak{sp}(W_2,\psi_2)$.  Then the projection
of $\hg_\CC$ onto $\frak{sp}(W_i,\psi_i)$ acting on $W_i$ must be on the
list of Theorem~2.11, and since it is irreducible, symplectic,
$4$-dimensional, it must be equal $\frak{sp}(W_i,\psi_i)$. Then since all
automorphisms of $\frak{sp}_{4\,\CC}$ are inner, Proposition~2.16.5
implies
that
$\hg_\CC = \frak{sp}(W_1,\psi_1) \times \frak{sp}(W_2,\psi_2)$, and thus
$\hg = \frak{sp}_F(W,\psi)$.
\enddemo

\proclaim{5.8. Type I(4)}
  Let $A$ be a simple abelian fourfold with $\EndoA = F$ a totally real
field with $[F:\QQ] =4$.  Then in the notation of~5.6,  $\hg(A) \simeq
\frak{sp}_F(W,\psi) \simeq \frak{sl}_{2\, F}$.  In particular, $\Hdg(A^n)
= \Div(A^n)$ for all ~$n$.
\endproclaim

\remark\nofrills
 The method of proof is similar to and easier than the previous case.
Over $\Qbar$ or $\CC$ the representation of $\hg$ on $W$ splits into a sum
of four mutually nonisomorphic, irreducible, symplectic, $2$-dimensional
representations, and Proposition~2.16.3 applies.
\endremark

\proclaim{5.9. Type II}
  Let $A$ be a simple abelian fourfold of type~\rom{(II)}, i.e., $\EndoA$
is an indefinite quaternion algebra $D$ over a totally real field $F$ of
degree $e\in\{1,2\}$ over~$\QQ$.  Then $\hg(A)$ is the centralizer of $D$
in $\frak{sp}(W,E)$.  In particular, $\Hdg(A^n) = \Div(A^n)$ for all ~$n$.
\endproclaim

 For both $e=1$ and $e=2$ this is a special case of \cite{B.21}~Thm.4.10
and \cite{B.22}~Thm.7.4.  Compare also with the method of
Ribet \cite{B.94}, 6.3 below.

\proclaim{5.10. Type III}
   Let $A$ be a simple abelian fourfold of type~\rom{(III)}, i.e.,
$\EndoA$ is a definite quaternion algebra $D$ over $\QQ$.  Then $\hg(A)$
is the centralizer of $D$ in $\frak{sp}(W,E)$, which is a $\QQ$-form of
$\frak{so}_4$.  Moreover, $\dim \Hdg^2(A) =6$, and $\dim \Div^2(A) = 1$,
and $\Hdg^2(A) = \Div^2(A) + \V(A)$, where $\V(A)$ is as in~5.2.
\endproclaim

\proclaim{5.11. Type IV(1,1)}
 Let $A$ be a simple abelian fourfold such that $\EndoA =K$ is an
imaginary quadratic field.
\roster
\item"(i)" If $K$ acts with multiplicities $\{1,3\}$ then $\hg(A) = \frak
u(W/K)$, and $\Hdg(A^n) = \Div(A^n)$ for all ~$n$.
\item"(ii)" If $K$ acts with multiplicities $(2,2)$.  In this case $\hg(A)
= \frak{su}(W/K)$, and $\dim \Hdg^2(A) =3$, and $\dim \Div^2(A) = 1$, and
$\Hdg^2(A) = \Div^2(A) \oplus \V(A)$, where $\V(A)$ is as in~5.2.
\endroster
\endproclaim

\remark\nofrills
 In case~(i) $A$ is of Ribet type (1.13.2), see 6.3 below.  In case~(ii)
the equality $\hg(A) = \frak{su}(W/K)$ must be proved, and then
Theorem~4.11 applies.
\endremark

\proclaim{5.12. Type IV(2,1)}
   Let $A$ be a simple abelian fourfold such that $\EndoA=K$ is a CM-field
of degree~$4$ over $\QQ$.  Then $\hg(A) \simeq \frak u_K(W,\psi)$.  In
particular, $\Hdg(A^n) = \Div(A^n)$ for all ~$n$.
\endproclaim

 Similarly as in 5.6, here $\psi : W\times W \to K$ is the unique
$K$-Hermitian form such that $\operatorname{Tr}_{K/\QQ}(\alpha\cdot
\psi(u,v))$ is a Riemann form for $A$, for $\alpha \in K$ such that $\bar
\alpha = -\alpha$ (The uniqueness is proved in \cite{B.27}).  In this
case $\hg(A)$ is contained in
$$
\frak u_K(W,\psi) := \{ \phi \in \End_K(W) : \psi(\phi(u),v) +
\psi(u,\phi(v)) =0 \enspace\text{for all } u,v \in W \} ,
$$
 regarded as a Lie algebra over $\QQ$.

\proclaim{5.13. Type IV(4,1)}
 Let $A$ be a simple abelian fourfold such that $\EndoA =K$ is a CM-field
of degree~$8$ over $\QQ$.
\roster
\item"(i)" If $K$ does not contain an imaginary quadratic field $F$ acting
on $A$ with multiplicities $(2,2)$, then $\hg(A) = \frak u_K$, which is a
commutative Lie algebra of rank~$4$, and $\Hdg(A^n) = \Div(A^n)$ for all
{}~$n$.
\item"(ii)" If $K$ does contain an imaginary quadratic field $F$ acting on
$A$ with multiplicities $(2,2)$, then $\hg(A) = \frak{su}_{K/F}$.  In this
case $\dim \Hdg^2(A) =8$, and $\dim \Div^2(A) = 6$, and $\Hdg^2(A) =
\Div^2(A) + \V(A)$.
\endroster
\endproclaim

 We discuss abelian varieties with complex multiplication further below,
see section nine.



\head 6. Simple abelian varieties with conditions on dimension or
endomorphism algebra
\endhead
\rightheadtext{simple abelian varieties with conditions}
  After Tankeev's early work on simple abelian fourfolds \cite{B.122}
\cite{B.123}, the next progress on the Hodge $(p,p)$
conjecture was the work of Tankeev \cite{B.124} \cite{B.125}
\cite{B.126} and Ribet \cite{B.93} \cite{B.94}
on simple abelian varieties of types~(I), (II) or~(IV) whose dimension and
endomorphism algebras satisfy various conditions.  More precisely, Tankeev
proved the following.

\proclaim{6.1. Theorem {\rm (\cite{B.125} \cite{B.126})}}
 Let $A$ be a simple abelian variety of dimension~$d$.  Then if
\roster
\item $A$ is of nondegenerate CM-type \rom{(2.13)}, or
\item $\EndoA$ is a totally real field of degree~$e$ over $\QQ$, and $d/e$
is odd, or
\item $\EndoA$ is a totally indefinite division quaternion algebra over a
totally real field $K$ of degree~$e$ over $\QQ$, and $d/2e$ is odd, or
\item $d$ is a prime,
\endroster
 then $\Hdg(A) = \Div(A)$.
\endproclaim

 However, soon afterward Ribet extended some of those results by first
observing the following basic criterion, which he says was used implicitly
in Tankeev's work, and then identifying some instances where it is
satisfied.

\proclaim{6.2. Theorem {\rm (\cite{B.94}~Theorem~0)}}
 Let $A$ be an abelian variety, and suppose
\roster
\item"(a)" $\EndoA$ is a commutative field, and
\item"(b)" $\Hg(A) = \Lf(A)$ (the Lefschetz group, see~2.14).
\endroster
 Then $\Hdg(A^n) =\Div(A^n)$ for $n\ge 1$.
\endproclaim

\proclaim{6.3. Theorem {\rm (\cite{B.94}~Theorems~1--3)}}
 Let $A$ be an abelian variety of dimension~$d$, and suppose
\roster
\item $\EndoA$ is a totally real field of degree~$e$ over $\QQ$, and $d/e$
is odd, or
\item $d$ is prime and $A$ is of CM-type, or
\item $\EndoA$ is an imaginary quadratic field $K$, and the multiplicities
$n'$ and $n''$ with which $\alpha \in K$ acts as $\alpha$ and $\bar
\alpha$ respectively are relatively prime.
\endroster
 Then $\Hg(A) = \Lf(A)$ and thus $\Hdg(A^n) =\Div(A^n)$ for $n\ge 1$.
\endproclaim

\proclaim{Corollary}
 When $A$ is a simple abelian variety of prime dimension, then $\Hdg(A^n)
=\Div(A^n)$ for $n\ge 1$.
\endproclaim

\demo\nofrills
 For if $A$ is simple and of prime dimension, then one of the conditions
of Theorem~6.3 must be satisfied, see 1.13.3.
\enddemo

\remark{Remark}
 In \cite{B.139} Yanai showed that a prime-dimensional abelian
variety of simple CM-type is nondegenerate (2.13).
\endremark

 In a similar spirit as 6.3.2 above, Hazama proved the following.

\proclaim{6.4. Theorem {\rm (\cite{B.45})}}
 Let $A$ be a simple abelian variety of CM-type.  Then $\Hdg(A^n)
=\Div(A^n)$ for all $n$ if and only if $\dim \Hg(A) =\dim A$.
\endproclaim

\demo{Sketch of proof of Theorem~6.2 {\rm (after \cite{B.94})}}
  In order to give some flavor of the techniques involved, consider the
proof of Theorem~6.2.
  Since it is always the case that $\Hg(A) \subseteq \Lf(A)$, see 2.14,
the
condition that these two groups are equal should be thought of as saying
that $\Hg(A)$ is as large as possible, whence has as few invariants as
possible.  Then the proof is separated into two cases, according as
$\EndoA =K$ is a totally real or a CM-field.

 Consider first the case where $K$ is totally real.  Then similarly as in
5.6 there is a unique $K$-bilinear alternating form $\psi: W\times W \to
K$ whose trace $\operatorname{Tr}_{K/\QQ}(\psi(u,v))$ is a Riemann form
$E(u,v)$ representing the chosen polarization on $A$, where $W=H_1(A,\QQ)$
\cite{B.27}~4.7.  The uniqueness of $\psi$ implies that a
$K$-automorphism of $W$ preserves $\psi$ if and only if it preserves $E$.
Then from the definition 2.14 we get that $\Hg(A) = \Lf(A)$ is the
symplectic group of $\psi$ acting on $W$ as a $K$-vector space, i.e.,
$$
 \Hg(A) = \Res_{K/\QQ}\Sp(W_0, \psi),
$$
 where $W_0$ is $W$ as a $K$-vector space.

  Now what we need to prove is that $(\twedge_\QQ^*((W\dual)^n))^{\Hg(A)}$
is generated by its elements of degree~$2$, and for this question it
suffices to extend scalars to $\CC$.  Then $W\tensor \CC$ is a free
$K\tensor_\QQ \CC$-module of rank $r= 2\dim A /[K:\QQ]$, and thus
$$
 W\tensor \CC \simeq \bigoplus_{\sigma\in\Hom(K,\CC)} U_\sigma,
$$
 where $U_\sigma$ is a complex vector space of dimension~$r$.  Therefore
$$
 \Hg(A,\CC) \simeq \prod_{\sigma\in\Hom(K,\CC)} \Sp(U_\sigma, \psi_\sigma)
,
$$
  from which we see that
$$
 (\twedge_\QQ^*((W\dual)^n))^{\Hg(A)} \tensor\CC \simeq \bigotimes_\sigma
\ (\twedge_\CC^* (U_\sigma\dual)^n)^{\Sp(U_\sigma, \psi_\sigma)} .
$$
  So it suffices to know that this algebra of invariants is generated by
its elements of degree~$2$, which is the case; see \cite{B.94}, or
derive this fact using the methods of \cite{B.10},
\cite{B.137} or \cite{B.34}.

  The proof for the case where $\EndoA =K$ is a CM-field, either of degree
$2\dim A$ or of degree $2$ over $\QQ$ follows a very similar pattern.  As
in 5.12 there is an element $\alpha \in K$ such that $\bar \alpha =
-\alpha$ and a unique Hermitian form $\psi : W\times W \to K$ such that a
Riemann form representing a polarization on $A$ is given by $E(u,v) =
\operatorname{Tr}_{K/\QQ}(\alpha\psi(u,v))$.  Then the centralizer
$\Lf(A)$
of $K$ in $\Sp(W,E)$, which by hypothesis coincides with $\Hg(A)$, is
$\Res_{K_0/\QQ}\U(W_0, \psi)$, where $K_0$ is the maximal totally real
subfield of $K$ and $W_0$ is $W$ as a $K$-vector space.

 Now when we extend scalars to $\RR$,
$$
 W\tensor_\QQ \RR \simeq \bigoplus_{\sigma\in \Hom(K_0,\RR)} U_\sigma.
$$
 Moreover, $\psi$ induces a nondegenerate Hermitian form $\psi_\sigma$ on
each $U_\sigma$, from which
$$
 \Hg(A,\RR) \simeq \prod_{\sigma\in \Hom(K_0,\RR)}
\U(U_\sigma,\psi_\sigma) .
$$
 Thus in this case we need to know that for each $\sigma$ and for all~$n$
the algebra of invariants
$(\twedge_\RR^*(U\dual_\sigma)^n)^{\U(U_\sigma,\psi_\sigma)}$ is generated
by elements of degree~$2$, which is the case.  The first step towards
proving this is to extend scalars to $\CC$, so that the unitary group
$\U(U_\sigma,\psi_\sigma)$ becomes a general linear group; we omit the
invariant theory arguments here, see {\sl op.~cit.}\space for more
details.
\Qed
\enddemo

\demo{Sketch of proof of Theorem~6.3.3}
 We are assuming that $\EndoA =K$ is an imaginary quadratic field, and the
multiplicities $n',~n''$, with which $\alpha\in K$ acts as $\alpha$ and
acts as $\bar\alpha$ are relatively prime, and we want to show that
$\Hg(A) =\Lf(A)$.  In fact, it turns out to be more convenient to show
that
$\MT(A) = \Gm \cdot\Hg(A)$ coincides with $G = \Gm\cdot \Lf(A)$.  This
group
$G$ may also be described as the largest connected subgroup of the
symplectic similitude group $\GSp(W,E)$ that commutes with $K$, and in the
present case $G=\Gm\cdot\Res_{K/\QQ}\U(W,\psi)$.  It is clear that $\MT(A)
\subseteq G$.

 Now, if $d= \dim A$, then $W$ is free of rank~$d$ as a vector space
over~$K$, and thus $W\tensor_\QQ\CC$ is free of rank~$d$ over
$K\tensor_\QQ \CC \simeq \CC\oplus \CC$, where the two copies are
naturally indexed by the embeddings of $K$ into $\CC$.  Therefore we may
write
$$
 W\tensor_\QQ \CC = W' \oplus W'',
$$
 and this decomposition is compatible with the Hodge decomposition of and
the action of $\MT(A)$ on $W\tensor_\QQ \CC$, since it is induced by
endomorphisms of $A$.  In particular,
$$
 W' = (W'\cap H^{-1,0}(A)) \oplus (W'\cap H^{0,-1}(A)) .
$$

  Now the action of $\MT(A,\CC)$ on $W\tensor_\QQ \CC$ induces an action
of $\MT(A,\CC)$ on $W'$, and the next step of the proof is to see that the
induced map $\MT(A,\CC) \to \GL(W')$ is surjective.  However, this follows
from \cite{B.106}~Prop.5; it is here that the relative primality of
$n'$ and $n''$ is a required hypothesis.  It follows that the commutator
subgroup of $\MT(A,\CC)$ maps onto $\SL(W')$.  What this means is that
when we write $G$ as the product of its center~$C$ and its semisimple part
$G_{\text{ss}}$, then $\MT(A) \supset G_{\text{ss}}$.  Thus it remains to
show that $C\subset \MT(A)$ as well.  And since $\dim C =2$, it will
suffice to show that the dimension of the center of $\MT(A)$ is at
least~$2$, which in turn would follow from showing that $\MT(A)$ maps onto
a $2$-dimensional torus.

 Since $\MT(A) \subseteq \Res_{K/\QQ}\GL(W_0)$, where $W_0$ is $W$
considered as a vector space over $K$, we may consider the determinant map
$\Theta: \MT(A) \to \Res_{K/\QQ} \Gm_{/K} =: T$.  Then from the fact that
$\MT(A)$ contains $\Gm_{/\QQ}$ acting as homotheties on $W$, the image of
$\Theta$ contains $\Gm_{/\QQ} \subset T$.  If we extend scalars to $\CC$,
then $T_\CC \simeq \Gm_{/\CC} \times \Gm_{/\CC}$, and the image of
$\Gm_{/\QQ}$ in $T$ becomes the diagonal.  On the other hand, $\MT(A)_\CC$
also contains $h(\SS_\CC) \simeq \Gm_{/\CC} \times \Gm_{/\CC}$.  Then
$\Theta(h(z,1)) = (z^{(n')}, z^{(n'')})$, and because $n' \ne n''$, this
generates a torus distinct from the diagonal.  Therefore $\Theta$ is
surjective, which completes the proof.
\Qed
\enddemo



\head 7. More abelian varieties with conditions on dimension or
endomorphism algebra
\endhead
\rightheadtext{more abelian varieties with conditions}
  During the 1980's Hazama, Murty and others continued to generate results
related to the Hodge conjecture by examining the interactions among the
dimension, endomorphism algebra, and Hodge or Mumford-Tate group, in a
spirit akin to the work of Tankeev and Ribet described in section~6.  In
particular, Hazama and Murty, working at about the same time but using
different methods, produced a number of overlapping results about the
Hodge conjecture for not-necessarily-simple abelian varieties extending
the results of Tankeev and Ribet.

\subhead Abelian varieties with generalized real multiplication
\endsubhead
  A first set of results can be loosely grouped together as dealing with
abelian varieties with generalized real multiplication.

\definition{7.1. Definition}
 Several different definitions of what is meant by \dfn{real
multiplication} appear in the literature.  The narrowest would be that the
abelian varietiey $A$ is of type~(I), i.e., every simple factor $A_s$ of
$A$ has $\Endo(A_s)$ equal to a totally real field.  A slightly broader
definition would be to require that $\EndoA$ contains a product $R$ of
totally real fields such that $[R:\QQ] = \dim A$ \cite{B.46}.
Zarhin \cite{B.142} calls a $g$-dimensional abelian variety \dfn{of
RM-type} if it contains a commutative semisimple $\QQ$-algebra of
degree~$g$ over $\QQ$, and notes that this means that any abelian variety
of CM-type is also automatically of RM-type.  Murty variously considers
the cases where a commutative semisimple $\QQ$-algebra $R\subseteq \EndoA$
is its own centralizer in $\EndoA$ \cite{B.81}, or $R$ is maximal
among commutative semisimple subalgebras of $\EndoA$ and is a product of
totally real fields \cite{B.83}.  Altogether the most useful general
definition might be to say that an abelian variety $A$ has generalized
real multiplication if $\EndoA$ contains a commutative semisimple
subalgebra $R$ with $[R:\QQ] =\dim A$, and $A$ is not of CM-type, i.e.,
$\EndoA$ does not contain a commutative semisimple subalgebra of
degree~$2\dim A$ over $\QQ$.  To avoid ambiguity we will try to give
precise statements of results without using this terminology.
\enddefinition

  The following theorem tries to summarize the main results concerning the
Hodge $(p,p)$ conjecture for abelian varieties with some generalized real
multiplication.

\proclaim{7.2. Theorem {\rm (\cite{B.46} \cite{B.81}
\cite{B.83})}}
 Let $A$ be an abelian variety.
\roster
\item Suppose $\EndoA$ contains a product $R$ of totally real fields with
the property that $[R:\QQ] = \dim A$, and no simple component of $A$ of
CM-type has dimension greater than~$1$.  Then $\Hdg(A) =\Div(A)$.
\item Suppose $\EndoA$ contains a commutative semisimple subalgebra $R$
that is its own centralizer in $\EndoA$, and $H^0(A,\Omega^1)$ is free of
rank~$1$ over $R\tensor_\QQ \CC$.  Then $\Hdg(A) =\Div(A)$.
\item  Suppose that a maximal commutative semisimple subalgebra $R$ of
$\EndoA$ is a product of totally real fields, and that $W= H_1(A,\QQ)$ is
free over $R$ of rank $2m$, where $m$ is odd.  Then $\Hg(A) = \Lf(A)$ and
thus $\Hdg(A^k) = \Div(A^k)$ for all $k\ge 1$.
\endroster
\endproclaim

\remark{7.3.1. Remarks on Theorem 7.2.1 {\rm (\cite{B.46})}}
  The first observation about an abelian variety $A$ satisfying the
conditions of 7.2.1 is that its simple isogeny factors must be of
type~(I), or type~(II), or elliptic curves with complex mutliplication, as
elliptic curves without complex multiplication are included in type~(I)
\cite{B.37}.  In particular, the condition is stable under taking
products or abelian subvarieties.

 Without going into the proof at too great a length, some of the main
ingredients include firstly a lemma, due to Tankeev
\cite{B.123}~Lemma~1,4,
that if the center of $\EndoA$ is a product of totally real
fields, then $\Hg(A)$ is semisimple.  Then Hazama proceeds to work out the
Hodge Lie algebra $\hg(B)$, where $B$ is the isogeny factor of $A$
containing the simple factors of type~(I) or~(II).  After complexifying
and applying Goursat's Lemma he finds that $\hg(B,\CC) \simeq \frak{sl}_2
\times \dots  \times \frak{sl}_2$.  Finally, he observes, similarly as we
did in section~3 on elliptic curves, that when $B$ is a abelian variety
whose Hodge group is semisimple and $C$ is an abelian variety of CM-type,
then $\Hg(B\times C) \simeq \Hg(B)\times \Hg(C)$.  Then Theorem~7.2.1
follows from the invariant theory of $\frak{sl}_2$ and the known results
for elliptic curves.
\endremark

\remark{7.3.2. Remarks on Theorem 7.2.2 {\rm (\cite{B.81})}}
  Murty calls a pair $(A,R)$ consisting of an abelian variety and a
commutative semisimple subalgebra $R\subset \EndoA$ \dfn{of type~{\rm
(H)}} when the hypotheses of Theorem~7.2.2 are satisfied.  He observes
that the product of two abelian varieties of type~(H) is again of
type~(H), and further proves that in general when a commutative semisimple
subalgebra of $\EndoA$ of degree $\dim A$ over $\QQ$ is a product $R$ of
totally real fields, then $H^0(A,\Omega^1)$ is free of rank~$1$ over
$R\tensor \CC$.  Thus some examples of abelian varieties of type~(H)
include $(E,\QQ)$, where $E$ is an elliptic curve without complex
multiplication, and $(A,F)$ where $A$ is an abelian surface with
quaternionic multiplication by a quaternion algebra $B$ over $\QQ$, as in
1.13.7, and $F$ is any real quadratic subfield of $B$ which splits~$B$.

 Although the analysis is somewhat different, some of the main ideas of
the proof of Theorem~7.2.2 are similar to some of the main points of the
proof of Theorem~7.2.1.  In particular, under the hypothesis of type~(H)
(and making use of Goursat's Lemma), Murty shows that
$$
 \hg(A) = \{ m \in \End_{\EndoA}(W) : \operatorname{tr}_R m =0\} ,
$$
 and that not only is this semisimple, but over $\CC$ it is a product of
$\frak{sl}_2$'s.  Thus it is possible to deduce that all invariants are
generated by those of degree~$2$.
\endremark

\example{7.3.3. Example: Jacobians of elliptic modular curves}
  Among the motivating examples for both Hazama and Murty were Jacobians
of elliptic modular curves.  To recall briefly, for $N\ge 3$ let
$$
 \Gamma_1(N) := \left\{ \pmatrix a&b\\ c&d\endpmatrix \in \SL(2,\ZZ) :
c\cong 0 \text{ and } a\cong d \cong 1 \pmod N \right\} ,
$$
  and let $\frak H = \{z\in \CC : \Im z >0 \}$ denote the upper
half-plane.  Then the quotient
$$
 X_1(N)(\CC) :=  \Gamma_1(N) \big\bs (\frak H \cup \QQ \cup \{ i\infty\})
$$
 can be identified with the complex points of a nonsingular projective
algebraic curve which can, in fact, be defined over~$\QQ$.  Moreover,
Shimura has shown that in the Jacobian $J_1(N) :=
\operatorname{Jac}(X_1(N))$ all the isogeny factors with complex
multiplication are products of elliptic curves, and all the isogeny
factors without complex multiplication are of real multiplication type in
the sense that the endomorphism algebras of the simple factors contain a
totally real number field whose degree over $\QQ$ is the dimension of that
factor \cite{B.112} \cite{B.113} \cite{B.114} (see
also \cite{B.89}).
\endexample

\remark{7.3.4. Remark on Theorem 7.2.3}
 Theorem 7.2.3 is an artifact of Murty's study \cite{B.83} of the
semisimple parts of the Hodge groups of abelian varieties, and their
relationship with the Lefschetz groups.  We will return to this again
briefly below.
\endremark

\subhead  Stably nondegenerate abelian varieties
\endsubhead
  A next group of results concerns conditions under which $\Hdg(A^k)
=\Div(A^k)$ for all $k\ge 1$.  Again we combine closely relatied results
of Murty and Hazama; but first we need a definition.

\definition{7.4. Definition {\rm (\cite{B.47})}}
 When $A$ is a simple abelian variety, the \dfn{reduced dimension} of $A$
is defined by
$$
 \rdim A := \cases \dim A, &\text{for $A$ of type~(I) or of type~(III),}
\\
   (\dim A)/2, &\text{for $A$ of type~(II),} \\
  (\dim A)/d, &\text{for $A$ of type~(IV), and } [\EndoA: C(\EndoA)] =
d^2 .
\endcases
$$
 By $C(\EndoA)$ here we mean the center of $\EndoA$.
 When $A$ is isogenous to $\prod _i A_i^{m_i}$ with the $A_i$ simple and
nonisogenous, then the \dfn{reduced dimension} of $A$ is
$$
 \rdim A := \sum_i \rdim A_i .
$$
\enddefinition

\proclaim{7.5. Theorem {\rm (\cite{B.82}, \cite{B.47})}}
 For an abelian variety $A$, the following are equivalent.
\roster
\item $\Hdg(A^k) =\Div(A^k)$ for all $k\ge 1$.
\item $A$ has no factor of type~\rom{(III)}, and $\Hg(A) = \Lf(A)$.
\item $\rank \Hg(A)_\CC = \rdim A$.
\endroster
\endproclaim

\definition{7.6. Definition}
 An abelian variety satisfying the conditions of Theorem~7.5 may be called
\dfn{stably nondegenerate.}
\enddefinition

\remark{7.6.1. Remarks}
 Hazama makes the following elementary observations about stable
nondegeneracy \cite{B.47}:
\roster
\item If $A$ is stably nondegenerate, and $B$ is an abelian subvariety of
$A$, then $B$ is stably nondegenerate.  For up to isogeny $A\simeq B\times
B'$, and thus if stable nondegeneracy (in the sense of 7.5.1) failed for
$B$ it would fail for $A$.
\item For any $k\ge 1$, $A$ is stably nondegenerate if and only if $A^k$
is stably nondegenerate.  This follows from the definition 7.5.1 and the
previous observation.
\item For abelian varieties $A_i$ and integers $k_i$, the product $\prod_i
A_i^{k_i}$ is stably nondegenerate if and only if $\prod_i A_i$ is stably
nondegenerate.  Observe that
$$
 \prod_i A_i^{k_i} \subset \big( \prod_i A_i \big)^{\max k_i} .
$$
\endroster
\endremark

 He also proves the following.

\proclaim{7.6.2. Theorem {\rm (\cite{B.49)}}}
 If $A$ and $B$ are stably nondegenerate abelian varieties and contain no
factors of type~(IV), then $A\times B$ is also stably nondegenerate.
\endproclaim

 The difficulty with type~(IV) arises in taking products of, or with,
abelian varieties of CM-type, see section nine below.  Or from a different
point of view, there is the theorem of Tankeev
that if all the simple factors of
an abelian variety $A$ are of types~(I), (II) or (III), then $\Hg(A)$ is
semisimple \cite{B.123}, and this may fail for type~(IV).   What
can be said is that if $A$ is stably nondegenerate and has no factors of
types~(IV), and $B$ is stably nondegenerate and of CM-type, then $A\times
B$ is stably nondegenerate \cite{B.49}.  On the other hand, since
by 7.5.2 no abelian variety with a factor of type~(III) can be stably
nondegenerate, Theorem~7.6.2 applies when all simple factors of $A$ and
$B$ are of type~(I) or~(II).

\remark{7.7. Some remarks on Theorem~7.5}
 In \cite{B.82} Murty proves (1) if and only if~(2).  Much of the
paper is devoted to a careful analysis of the structure of $\Lf(A)$.
Given
the multiplicativity of $\Lf(A)$, Lemma~2.15, we may assume $A$ is simple.
Then fix a maximal commutative subfield $F \subset \EndoA$ which is
totally real for type~(I) and a CM-field in the other three cases, and let
$F_0$ be the maximal totally real subfield of $F$.  Now extending scalars
to $\RR$, there is a decomposition of $\Lf(A)$ into factors indexed by the
embeddings $F_0 \hra \RR$.  Then these factors are of the form: for
type~(I), a symplectic group; for type~(II), the intersection of a unitary
group and a symplectic group; for type~(III), the intersection of a
unitary group and a special orthogonal group; for type~(IV), a unitary
group \cite{B.82}~Lemma~2.3.  Moreover, after complexifying, these
act on the corresponding components of $W\tensor \RR$ as: for type~(I), as
a standard symplectic repesentation; for type~(II), two copies of the
standard representation of the complex symplectic group; for type~(III),
two copies of the standard representation of the complex special
orthogonal group; for type~(IV), the sum of a standard representation of
the complex general linear group and its contragredient.  Using this
structural analysis, Murty is then able to prove the following.

\proclaim{7.7.1. Proposition {\rm (\cite{B.82})}}
 If $A$ contains no simple factors of type~\rom{(III)}, then for all
$k\ge 1$
$$
 H^*(A^k, \QQ)^{\Lf(A)} = \Div(A^k) .
$$
\endproclaim

  To complete the proof that 7.5.1 is equivalent to 7.5.2, Murty shows
that a simple abelian variety of type~(III) supports an exceptional Hodge
class, see~8.6 below for more discussion of this.

\smallpagebreak

 Hazama's proof in \cite{B.47} that 7.5.1 is equivalent to 7.5.3 is
based on a careful type by type analysis of the Lie algebra $\hg(A)_\CC$
and its action on $W_\CC$ using that all the possibilities are as listed
in Theorem~2.11.  For example, for a simple abelian variety of type~(I),
the action of a simple component $\frak g_i$ of $\hg(A)_\CC$ on the
corresponding component $W_i$ of $W_\CC$ is a symplectic representation,
and indeed $\frak g_i \simeq \frak{sp}(W_i,\CC)$.  A similar result holds
for type~(II), whereas for type~(IV), the simple components of
$\hg(A)_\CC$ are of the form $\frak{sl}_{d_i}$.  In all these cases,
careful invariant theory arguments using \cite{B.137} show that the
invariants of $H^*(A^k,\CC)$ are generated by those of degree~$2$ if and
only if the rank is as claimed.  On the other hand, for type~(III) Hazama
finds that
$$
 \rank \Hg(A)_\CC \le (\dim A)/2 < \rdim A ,
$$
 i.e., equality never holds, and an abelian variety with a factor of
type~(III) fails to be stably nondegenerate.  It comes about in the proof
that in general
$$
 \rank \Hg(A) \le \rdim A .
$$
  So both criteria 7.5.2 and 7.5.3 can be understood philosophically as
saying that $A$ is stably nondegenerate when $\Hg(A)$ is as large as
possible.
\endremark

\subhead Further work on Hodge and Mumford-Tate groups
\endsubhead
  We conclude this section with certain additional results derived from
close study of Hodge and Mumford-Tate groups.

\subhead  \nofrills \endsubhead
  In \cite{B.83} Murty examines the semisimple part of the Hodge
group of an abelian variety, and finds the following.  As usual, $W=
H_1(A,\QQ)$.

\proclaim{7.8. Theorem {\rm (\cite{B.83})}}
 If a maximal commutative subalgebra $R$ of $\EndoA$ is a product of
CM-fields, and $W$ is free over $R$ of odd rank, and if $\Hdg(A)
=\Div(A)$, then $\Hg(A)_{\text{ss}} = \Lf(A)_{\text{ss}}$.
\endproclaim

 This together with 7.2.3 implies the following.

\proclaim{Corollary {\rm (\cite{B.83})}}
 When $A$ is simple and of odd dimension, then $\Hdg(A) = \Div(A)$ implies
that $\Hg(A)_{\text{ss}} = \Lf(A)_{\text{ss}}$.
\endproclaim

\subhead \nofrills \endsubhead
  In \cite{B.57} Ichikawa studies groups of Mumford-Tate type,
and extending \cite{B.94} \cite{B.126} and his own earlier
work \cite{B.56}, he obtains the following result.  First we need
some notation.

\definition{7.9. Definition}
  Let $A$ be a simple abelian variety of dimension~$g$, let $K$ the center
of $\EndoA$, let $e=[K:\QQ]$ and let $d^2=[\EndoA:K]$.  Then the
\dfn{relative dimension} of $A$ is defined by
$$
 \operatorname{rel\,dim}(A) := \cases g /e, &\text{if $A$ is of
type~(I),}\\
 g /2e, &\text{if $A$ is of type~(II) or type~(III),}\\
 2g/de, &\text{if $A$ is of type~(IV).}
\endcases
$$
\enddefinition

\proclaim{7.10. Theorem {\rm (\cite{B.57})}}
 Let $A$ be an abelian variety all of whose simple factors have odd
relative dimension.
\roster
\item When $A$ is isogenous to $A'\times A''$, where each simple factor of
$A'$ is of type ~\rom{(I)}, \rom{(II)} or~\rom{(III)} and each simple
factor of $A''$ is of type~\rom{(IV)}, then all Hodge cycles on $A$ are
generated by the Hodge cycles on $A'$ and $A''$.
\item When $A$ is isogenous to $\prod_j A_j^{m_j}$, where the $A_j$ are
simple and mutually non-isogenous, then all Hodge cycles on $A$ are
generated by the Hodge cycles on the $A_j$.
\endroster
\endproclaim



\head 8.  Exceptional Hodge cycles
\endhead
  Thus far we have looked mainly at examples and conditions under which
the
Hodge $(p,p)$ conjecture is true.  Now we consider the known examples of
Hodge cycles that are not known to be algebraic, and thus might be
considered potential counterexamples to the conjecture.

\definition{8.1. Definition}
 By an \dfn{exceptional Hodge cycle} on $A$ we mean an element of
$\Hdg^p(A) = H^{2p}(A,\QQ) \cap H^{p,p}(A)$, for some $p$, which is not in
$\Div^p(A)$, that is to say, which cannot be written as a $\QQ$-linear
combination of classes of $p$-fold intersections of divisors.
\enddefinition

\subhead 8.2. Mumford's CM fourfold
\endsubhead
  Perhaps the first example of an abelian variety where $\Hdg(A) \ne
\Div(A)$ was Mumford's example of the abelian fourfold with complex
multiplication corresponding to a particular  CM-type (see~1.13.6) for the
splitting field of $(3X^4 - 6X^2 + X +1)(X^2 +1)$ \cite{B.88}.
This example is described in Lecture~7, 7.23--7.28.

\subhead 8.3. Abelian varieties of Weil type
\endsubhead
  It was Weil's observation, however, that the crucial feature of
Mumford's example was not that it was of CM-type, but rather that there
was an imaginary quadratic field $F$ acting on $A$ in such a way that
$\Lie(A)$ becomes a free $K\tensor \CC$-module, or equivalently, such that
$\alpha \in K$ acts as $\alpha$ and as $\bar \alpha$ with equal
multiplicity \cite{B.135}.  Moreover, as we saw in Theorem~4.11, the
general such abelian variety, what we now refer to as an abelian variety
of Weil type, has a $2$-dimensional space of exceptional Weil-Hodge cycles
in $\Hdg^n(A)$, where $\dim A =2n$.  In Theorem~4.12 we recalled Schoen's
examples of general abelian fourfolds with $K= \QQ(i)$ or
$K=\QQ(\sqrt{-3})$ where he showed that the Weil-Hodge cycles are
algebraic
\cite{B.104}, and little else is known.

\subhead 8.4. Abelian varieties of Fermat type
\endsubhead
  Shioda's work on abelian varieties of Fermat type, see 4.1--4.8 above,
provides examples of abelian varieties $A$ where, at least for some $p$,
the space of Hodge cycles $\Hdg^p(A) \supsetneqq \Div^p(A)$ but is
nonetheless generated by classes of algebraic cycles \cite{B.116},
see Theorems~4.4 and~4.5.  In the same work he also shows the existence,
for any $d\ge 2$, of an abelian variety $A$ of Fermat type whose Hodge
ring $\Hdg^*(A)$ is not generated by $\sum_{r=1}^{d-1} \Hdg^r(A)$, let
alone by $\Hdg^1(A)$.

\subhead 8.5. Exceptional cycles in codimension~$2$
\endsubhead
  In \cite{B.124} Tankeev produced a family of abelian varieties
of dimension~$4^m$ with exceptional cycles in codimension~$2$ when $m\ge
2$.

\proclaim{Theorem {\rm (\cite{B.124}~Thm.5.6)}}
 For any $m\ge 1$ there exist abelian varieties $A$ such that
\roster
\item $\dim A = 4^m$, and
\item $\EndoA = \QQ$, and
\item $\hg_\CC(A) \simeq (\frak{sl}_{2\, \CC})^{2m+1}$, acting on
$H^1(A,\CC) \simeq (\CC^2)^{\tensor(2m+1)}$ as the tensor product of a
standard representation of each factor.
\endroster
 Moreover, for any abelian variety satisfying these conditions, $\dim_\QQ
\Hdg^2(A)
\mathbreak
 = (4^m -1)/3$.  In particular, if $m\ge 2$ then $\Hdg^2(A)$ is
not generated by classes of intersections of divisors.
\endproclaim

  The existence part of this theorem is obtained by generalizing Mumford's
example in \cite{B.78} (not the example in \cite{B.88}
mentioned in 8.2 above) of an abelian fourfold $A$ with $\EndoA = \QQ$ and
thus not characterized by its endomorphism ring.  The computation of
$\dim_\QQ \Hdg^2(A)$ is proved by induction, and a computation with the
roots of $\hg_\CC(A)$ and the character of its representation.

\subhead 8.6. Abelian varieties of type~(III)
\endsubhead
  In the same paper where he proved that $\Hdg(A^k) =\Div(A^k)$ for all
$k\ge 1$ if and only if $A$ has no factor of type~(III) and $\Hg(A) =
\Lf(A)$ \cite{B.82}, Murty also proved the existence of an exceptional
Hodge cycle on abelian varieties of type~(III).

\proclaim{Theorem {\rm (\cite{B.82})}}
 If an abelian variety $A$ has a factor of type~\rom{(III)}, then it
supports an exceptional Hodge class $\omega$ with the property that
$\pi_1^*(\omega) \tensor \pi_2^*(\omega) \in \Div(A^2)$, where $\pi_1$ and
$\pi_2$ are the projections from $A^2 = A\times A$ to its first and second
factors respectively.
\endproclaim

\remark{Remark}
  In \cite{B.135}, where he presented abelian varieties of Weil type
as a place to look for counterexamples of the Hodge conjecture, Weil also
asked whether a weaker statement might be true, that is, whether the
presence of a Hodge cycle on $A$ might imply the presence of an algebraic
cycle on some power of $A$.  This result of Murty is the first example
where the Hodge conjecture itself is not known to be true, but Weil's
question is answered affirmatively.
\endremark

\demo\nofrills
 To get a flavor of the proof, suppose $A$ is simple and of type~(III),
and let $F$ be the center of $\EndoA$, let $m= \dim_{\EndoA} W$, and let
$d=(\dim A)/[F:\QQ]$.  Then $d=2m$, and by \cite{B.109}~Prop.15,
$m\ge 2$.  Then as a consequence of his analysis of $\Lf(A)$, Murty shows
that
$$
 \big( \twedge^* (W\dual)\big)^{\Lf(A)} \tensor_\QQ \RR =
\bigotimes_{\sigma
\in \Hom(F,\RR)} \big( \twedge^* X\dual_\sigma\big)^{\Lf(A)_\sigma}
$$
 with $\dim_\RR X_\sigma =4m$.  Then $X_\sigma \tensor \CC$ becomes
isomorphic to two copies of a standard representation of
$\operatorname{SO}(V,\psi)$ for a suitable $V$ and $\psi$.  Then by
\cite{B.137}~p.53 the covariant tensors of
$\operatorname{SO}(V,\psi)$ are generated by $\psi$ and the determinant,
say~$\Delta$.  Then $\Delta$ cannot be written as a polynomial in the
degree~$2$ invariant $\psi$, but $\Delta^2$ can.  Take $\omega$ to be the
class corresponding to $\Delta$.
\enddemo

\subhead 8.7. Determinant cycles
\endsubhead
  In \cite{B.57} Ichikawa uses the idea of Murty's construction
to develop a certain extension of the work of Tankeev and Ribet on simple
abelian varieties \cite{B.124} \cite{B.126} \cite{B.94},
see section~6.  Firstly he observes that on any abelian variety of
type~(I), (II) or~(III) there exist Hodge cycles that are $\CC$-linear
combinations of the determinant forms on the spaces $V$ as in the last
paragraph.  He calls these \dfn{determinant Hodge cycles.}  In this
language, Murty's result above is that on an abelian variety of type~(III)
no determinant cycle is generated by classes of divisors.  Then Ichikawa
proves the following result.  Recall the definition of relative dimension
from~7.9.

\proclaim{Theorem {\rm (\cite{B.57})}}
  Let $A$ be an abelian variety whose simple factors are all of odd
relative dimension, and suppose $A$ is isogenous to $A'\times A''$ where
each simple factor of $A'$ is of type ~\rom{(I)}, \rom{(II)}
or~\rom{(III)} and each simple factor of $A''$ is of type~\rom{(IV)}.
Further, assume that the relative dimension of any simple factor of $A$ of
type~\rom{(III)} is not equal to $\frac 1 2 \binom{2k}{k}$ for any power
$k$ of~$2$, and that $A''$ is a power of a simple abelian variety of odd
prime dimension.  Then any Hodge cycle on $A$ is generated by classes of
divisors and determinant Hodge cycles.
\endproclaim

\definition{8.8. Definition {\rm (\cite{B.47})}}
 Recall (definition~7.6) that a stably nondegenerate abelian variety is
one which satisfies the conditions of Theorem~7.5, in particular,
$\Hdg(A^k) =\Div(A^k)$ for all $k\ge 1$.  Then a \dfn{stably degenerate}
abelian variety $A$ is one which is not stably nondegenerate, that is,
$\Hdg^p(A^n) \supsetneqq \Div^p(A^n)$ for some $p,n$.  Then the least $n$
for which this occurs is called the \dfn{index of degeneracy,} which we
will denote by $\ind(A)$.
\enddefinition

\subhead 8.9. Stably degenerate abelian varieties
\endsubhead
  Hazama has given two examples of stably degenerate abelian varieties of
type~(I) having index of degeneracy~$2$.

\proclaim{8.9.1. Theorem {\rm (\cite{B.47})}}
 There exists a simple abelian variety $A$ of dimension~$4$ with the
following properties:
\roster
\item"(a)" $A$ is of type~\rom{(I)},
\item"(b)" $\Hdg(A) = \Div(A)$,
\item"(c)" $\Hdg(A^2) \supsetneqq \Div(A^2)$.
\endroster
\endproclaim

\demo\nofrills
 To give a rough idea of the construction, let $K$ be a totally real
number field of degree~$3$, let $B$ be a quaternion algebra over $K$ such
that $B\tensor_\QQ \RR \simeq M_2(\RR) \times \Bbb H \times \Bbb H$, and
let $G = \Res_{K/\QQ}\SL(1,B)$.  Then $G(\RR) \simeq \SL(2,\RR) \times
\SU(2) \times \SU(2)$, and there exists an $8$-dimensional $\QQ$-rational
symplectic representation $\rho : G\to \Sp(8)$ satisfying the necessary
analyticity conditions so that the induced map $\tau: X \to \frak H_4$ of
Hermitian symmetric domains pulls back the universal family $\script A \to
\frak H_4$ of polarized abelian fourfolds over the Siegel upper half-space
to an analytic family of abelian fourfolds $A\to X$ \cite{B.64}.
Moreover, if $A_0$ denotes a generic member of the family $A\to X$, then
$\Hdg^p(A_0^k) = H^{2p}(A_0,\QQ)^G$ for all $k\ge 1$.  Then computations
with the complexified Lie algebra $\Lie(G,\CC) \simeq \frak{sl}_2 \times
\frak{sl}_2 \times \frak{sl}_2$ show firstly that $\dim \Hdg^1(A_0) =1$,
from which follows that $A_0$ is simple, of type~(I), and $\Hdg(A)
=\Div(A)$, and secondly that $\Hdg^2(A^2)$ is not generated by the
elements of $\Hdg^1(A^2)$.
\enddemo

\proclaim{8.9.2. Theorem {\rm (\cite{B.48})}}
 There exists a simple abelian variety $A$ of dimension~$10$ with the
following properties:
\roster
\item"(a)" $\hg(A,\CC) \simeq \frak{sl}(6,\CC)$,
\item"(b)" the representation $\hg(A,\CC) \to \End(H^1(A,\CC))$ is
equivalent to the representation $\frak{sl}(6,\CC) \to \End(\twedge^3
\CC^6)$ induced by the natural action of $\frak{sl}(6,\CC)$ on~$\twedge^3
\CC^6$,
\item"(c)" $\Hdg(A) = \Div(A)$,
\item"(d)" $\Hdg(A^2) \supsetneqq \Div(A^2)$.
\endroster
\endproclaim

\demo\nofrills
  The existence is worked out similarly as in the previous case, except
that here $G$ is a $\QQ$-form of $\SU(5,1)$.  Again $\dim \Hdg^1(A) = 1$,
and $\EndoA =\QQ$.  However, the actual computations are based on using
Young diagrams and branching rules, see \cite{B.48} for the
details.
\enddemo

\remark{Remark}
 In \cite{B.48} Hazama constructs a family $A_n$ of simple abelian
varieties of dimension $\frac 1 2 \binom{4n+2}{2n+1}$, with $\Endo(A_n)
=\QQ$, and $\hg(A_n,\CC) \simeq \frak{sl}(4n+2,\CC)$ and $H^1(A_n,\CC)
\simeq \twedge^{2n+1}\CC^{(4n+2)}$ as a representation of $\hg(A_n,\CC)$.
The abelian variety of Theorem~8.9.2 is the $A_1$ in this family.  Then he
also shows that the index of degeneracy $\ind(A_n) \le 2$ for $n\ge 2$,
where the theorem shows that $\ind(A_1) =2$.
\endremark

\subhead 8.10.  Invariants of partially indefinite quaternion algebras
\endsubhead
  In the 1960's Kuga asked which semisimple algebraic groups $G$ defined
over $\QQ$ together with which of their symplectic representations $\rho:
G\to \Sp(W,\beta)$ satisfy the necessary and sufficient analyticity
conditions to allow the construction of an algebraic family of polarized
abelian varieties parameterized by $\Gamma\bs X$, where $\Gamma$ is a
discrete subgroup of $G$ and $X$ is the Hermitian symmetric domain
associated to~$G$ \cite{B.61} \cite{B.62} \cite{B.63}.
Shortly thereafter Satake answered Kuga's question under the assumption
that for each $\QQ$-simple factor $\rho$ comes from an absolutely
irreducible representation of an absolutely simple factor of~$G$
\cite{B.99} \cite{B.100} \cite{B.102} \cite{B.103}.
It turned out that the list was quite small, and nearly all cases
had been considered by Shimura in his analysis of families of abelian
varieties characterized by polarization, endomorphism ring and level
structure \cite{B.109} \cite{B.110} \cite{B.112};
one more case was treated in \cite{B.101}.  Some time later
Addington considered Kuga's question without Satake's assumption, and for
the groups corresponding to units of norm~$1$ in a partially indefinite
quaternion algebras $B$ over a totally real field~$F$, i.e.,
$$
 B\tensor_\QQ \RR \simeq M_2(\RR)^n \oplus \Bbb H^m  \qquad \text{and}
\qquad G = \Res_{F/\QQ}\SL(1,B) ,
$$
 she developed a combinatorial scheme (called ``chemistry'') to describe
which symplectic representations give rise to an algebraic family of
abelian varieties \cite{B.5}.  Then Tjiok \cite{B.131} and
Abdulali \cite{B.1} \cite{B.2} \cite{B.3}
showed that under certain reasonable hypotheses (``rigidity'' or
``condition~(H$_2$)'') the space of Hodge cycles in a generic fiber $A_0$
of the family is the space of $G$-invariants, $H^{2r}(A_0,\QQ)^G =
(\twedge^{2r} W\dual)^G$.   Thus, for the purposes of this appendix, where
the issue is Hodge cycles on abelian varieties, statements about Hodge
cycles on the generic fiber of such a family can be understood as
statements about Hodge cycles on an abelian variety $A_0$ with specified
semisimple Hodge group~$G$.  Then the problem is to describe the
invariants of $G$ in the exterior algebra $\twedge^*W\dual$.

  This is the problem taken up by Kuga in \cite{B.64}, \cite{B.65}
and the series of papers \cite{B.66} \cite{B.67}, and Lee in
\cite{B.71}.  The results are rather involved to
state precisely.  In \cite{B.64} Kuga looks at some simple examples
of the situation just described, and finds conditions (``totally
disconnected triangular polymer'') where all the Hodge cycles in the
abelian variety $A_0$ are generated by those of degree two or where only
the Hodge cycles in codimension~$2$ or~$4$ are generated by those of
degree two (``triangular polymer without double bond, short cycle or
Hexatram''), which is to say that there are exceptional cycles in higher
codimensions.  At the end of \cite{B.65} is an example of a
$16$-dimensional abelian variety $A$ for which the dimensions of the
spaces $\Hdg^r(A)$ of Hodge cycles are determined, where $\dim \Hdg^2(A)=
82$ and $\dim \Div^2(A) = 10$.  In \cite{B.66} and \cite{B.67}
the focus is more on the complicated invariant theory in
the exterior algebra for the groups and representations under
consideration, in particular the latter papers examine the asymptotic
behavior of $\dim \twedge^{2r}(\mu W\dual)^G$ as the multiplicity~$\mu$
grows.  In \cite{B.71} Lee computes the dimensions of the spaces
$\Hdg^r(A_0)$ for the $8$-dimensional abelian variety constructed from a
particular form of $16$-dimensional representation of $G$.



\head 9. The problem of complex multiplication
\endhead
  In this section we look at what is known about the Hodge $(p,p)$
conjecture for abelian varieties with complex multiplication.  We have
already seen that the Hodge group of such an abelian variety is an
algebraic torus, necessarily contained in $\Res_{K/\QQ}\Gm_{/K}$, where
$K$ is the field of complex mutliplication (Proposition~2.12).  For the
general theory of complex multiplication, see \cite{B.115}
\cite{B.70} and parts of \cite{B.112}.

\definition{9.1. Definition}
 Recall (1.13.6) that a CM-type $(K,S)$ consists of a CM-field $K$ and a
subset $S\subset \Hom(K,\CC)$ containing exactly one from each pair of
conjugate embeddings.  Moreover, given a CM-type $(K,S)$, there is a
natural construction of an abelian variety $A$ with that CM-type, that is,
with $\EndoA = K$ and $K$ acting on $H^{1,0}(A)$ as $\bigoplus_{\phi\in S}
\phi$.  Then we define the \dfn{rank} of the CM-type $(K,S)$ by
$$
 \rank (K,S) := \dim \MT(A) .
$$
\enddefinition

\remark{Remark}
 The rank of a CM-type seems to have originally been defined by Kubota
\cite{B.60}, who defined it as $\dim_\QQ\{\sum_{\phi\in S} \phi(x)
: x\in K\}$.  The equality of this with the dimension of the Mumford-Tate
group follows from the methods in \cite{B.91}, see also
\cite{B.26} \cite{B.27}.
\endremark

\subhead 9.2. Pohlmann's criterion
\endsubhead
  One of the first results about Hodge cycles on abelian varieties with
complex multiplication is a theorem of Pohlmann \cite{B.88} that
describes $\dim \Hdg(A)$ in terms of the Galois theory of~$K$.  To state
the theorem we need some notation.  Let $A$ be an abelian variety with
CM-type $(K,S)$, let $S = \{\phi_1, \dots , \phi_g\}$ and let $\Sbar = \{
\phibar_1, \dots ,\phibar_g\}$, so $\Hom_\QQ(K,\CC) = S\cup \Sbar$.  Then
$\alpha \mapsto (\phi_1(\alpha), \dots , \phi_g(\alpha))$, for $\alpha \in
K$, induces an isomorphism of $K$ onto $H_1(A,\QQ)$ (in 1.13.6 we mapped
$\script O_K$ onto $H_1(A,\ZZ)$), via which $H^1(A,\CC)$ can be identified
with $\Hom_\QQ(K,\CC)$.  Further, without loss of generality we may assume
$K\subset \CC$, and let $L$ be the Galois closure of $K$ in $\CC$ and $G =
\Gal(L/\QQ)$.  Then we let $\sigma \in \Aut(\CC/\QQ)$ act on $f\in
H^r(A,\CC)$, with $f:\twedge^r K \to \CC$, by $(\sigma f)(\lambda) =
\sigma(f(\lambda))$ for $\lambda \in \twedge^r K$.  Finally, for an
ordered subset $\Delta \subset S$, let $|\Delta|$ denote the cardinality
of $\Delta$ and let $\angled \Delta := \twedge_{\phi\in S} \phi$.

\proclaim{Theorem {\rm (\cite{B.88}~Thm.1)}}
 When $A$ is an abelian variety with CM-type $(K,S)$, then $\Hdg^p(A)
\tensor \CC$ has a basis consisting of those $\angled \Delta \in
H^{2p}(A,\CC)$ such that
$$
 |\tau \Delta \cap S| = | \tau \Delta \cap \Sbar|
\tag 9.2.1
$$
 for every $\tau \in G$.  Thus $\dim \Hdg^p(A)$ is the number of ordered
subsets $\Delta \subset S$, with $|\Delta| = 2p$, that satisfy the
condition~\rom{(9.2.1)}.
\endproclaim

\demo{Proof}
 If $f= \sum_i c_i\angled{\Delta_i} \in \Hdg^p(A)$ with $c_i\in\CC$, then
$\sum_i \sigma(c_i)\angled{\sigma\Delta_i} = \sigma f = f$ is in
$H^{p,p}(A)$ for all $\sigma\in \Aut(\CC/\QQ)$, hence $\Delta_i$
satisfies~(9.2.1).  Then every element of $\Hdg^p(A)$ is a linear
combination of $\Delta_i$ satisfying~(9.2.1).  Conversely, let $\Delta$ be
such taht $|\Delta | = 2p$ and (9.2.1) is satisfied.  Let
$\{u_1,\dots,u_s\}$ be a basis for $L$ over $\QQ$, and let $f_i =
\sum_{\tau\in G} \tau(u_i) \angled{\tau\Delta}$ for $1\le i \le s$.  Then
$\sigma f_i = f_i$ for all $\sigma\in\Aut(\CC/\QQ)$, and by (9.2.1) $f_i
\in H^{p,p}(A)$, so $f_i\in \Hdg^p(A)$. Further, since
$\det(\tau(u_i))_{\tau, i} \ne 0$, we can solve the system of linear
equations $f_i = \sum_{\tau \in G} \tau(u_i) \angled{\tau \Delta}$ and
find that $\angled {\tau\Delta} \in \Hdg^p(A) \tensor \CC$ for $\tau \in
G$.  Thus $\angled \Delta \in \Hdg^p(A) \tensor \CC$, as required.
\Qed
\enddemo

  Pohlmann's theorem give a criterion for exceptional cycles, also see 9.3
below.

\proclaim{9.2.2. Corollary {\rm (\cite{B.138})}}
 $\dim \Hdg^p(A) - \dim\Div^p(A)$ is the number of subsets $\Delta \subset
\Hom(K,\CC)$ such that
\roster
\item"(a)" $\Delta - \ol\Delta \ne \varemptyset$,
\item"(b)" $|\Delta \cap g S| = p$ for all $g\in G$.
\endroster
\endproclaim

\subhead 9.3. Sporadic cycles
\endsubhead
  In \cite{B.138} White observes that Pohlmann's criterion shows that
when a CM abelian variety $A$ is nondegenerate, as defined in~2.13, then
$\Hdg(A) =\Div(A)$.  He then recounts that between 1977 and 1978 Ribet
asked if these two conditions were equivalent, and that Lenstra was
quickly able to show that they are, for a simple abelian variety of simple
CM-type $(K,S)$, under the additional hypothesis that the CM-field $K$ is
abelian over $\QQ$ (see \cite{B.138}~Thm.3).  Then later Hazama
showed that a simple abelian variety is nondegenerate if and only if
$\Hdg(A^k) = \Div(A^k)$ for all $k\ge 1$, see Theorem~6.4 \cite{B.45} and
Theorem~7.5 \cite{B.47}.  Only recently, however, White
showed the following.

\proclaim{Theorem {\rm (\cite{B.138}~Thm.1)}}
 There exists an abelian variety of CM-type with $\Hdg(A) =\Div(A)$ and
$\dim \Hg(A) \lneqq \dim A$.
\endproclaim

 The argument involves a rather technical analysis of the $\QQ$-group ring
of a non-abelian group.  Eventually, however, the counterexample is a
CM-type for a CM-field whose Galois group is
$$
 \ZZ/2\ZZ \times \ZZ/2\ZZ \times \ZZ/5\ZZ \times D_5,
$$
 where the last factor is the dihedral group and the first factor
corresponds to complex conjugation.  The splitting field of
$$
 X^5 - 10X^4 -70X^3 -25X^2 +190X +12
$$
 is a totally real field with Galois group $D_5$, and it is easy to make
disjoint totally real quadratic and quintic extensions, and then a totally
imaginary quadratic extension.  For the abelian variety $A$ with the
requisite CM-type for this field, $\dim\Hg(A) = 84 < 100 =\dim A$.

\subhead 9.4. Degenerate CM types
\endsubhead
   Recall from definitions 2.13 and 9.1 that a CM-type $(K,S)$ is said to
be nondegenerate if $\rank(K,S) = \dim A +1$, and is called degenerate
otherwise, if $\rank (K,S) \le \dim A$, where $A$ an abelian variety
CM-type $(K,S)$.  A number of examples of degenerate CM-types and lower
bounds for the rank as a function of $\dim A$ have been given by Ribet
\cite{B.92}, Dodson \cite{B.28} \cite{B.29}
\cite{B.30}, Mai \cite{B.72} and Yanai \cite{B.140}.

\smallpagebreak

  The following proposition of Kubota can be a useful way of measuring the
rank of a CM-type.  Let $c$ denote complex conjugation.

\proclaim{9.4.1. Proposition {\rm (\cite{B.60})}}
$$
\rank (K,S) = 1 + \#\big\{\chi : \Gal(K/\QQ) \to \CC : \chi(c) =-1\enspace
\&\enspace \sum_{s\in S} \chi(s) \ne 0 \big\},
$$
 where only irreducible $\chi$ are included.
\endproclaim

\example{9.4.2. Examples {\rm (\cite{B.92})}}
  First, let $p \ge 5$ be a prime, let $K = \QQ(\zeta_p)$ be the field of
$p^{\text{th}}$ roots of unity, and identify $G = \Gal(K/\QQ) \simeq
(\ZZ/p\ZZ)^\times$.  For $g\in G$ let $\angled g \cong g \pmod p$ with
$1\le \angled g \le p-1$.  Then for $1\le a \le p-2$ with $a^3\not \cong 1
\pmod p$ the set
$$
 S_a = \{g\in G : \angled g + \angled{ag} <p\}
$$
 is a simple CM-type.  It is nondegenerate when $a=1$, but is degenerate
for $p= 67$ and $a = 10$, $19$, $47$, $56$, $60$ \cite{B.40}.
Lenstra and Stark also noticed that for $p\cong 7 \pmod {12}$ and
sufficiently large there always exists a number~$a$ for which $S_a$ is
degenerate, {\sl loc..~cit.}

 Next let $K=\QQ(\zeta_{32})$ and, identifying $\Gal(K/\QQ)$ with
$(\ZZ/32\ZZ)^\times$, let
$$
 S = \{ 1,\ 7,\ 13,\ 21,\ 23,\ 27,\ 29\}, \qquad S'=\{1,\ 7,\ 9,\ 11,\
13,\ 15,\ 27,\ 29\} .
$$
 Then $(K,S)$ and $(K,S')$ are both degenerate (simple) CM-types.  This
example is due to Lenstra.

  Let $K= \QQ(\zeta_{19})$ and, identifying $\Gal(K/\QQ)$ with
$(\ZZ/19\ZZ)^\times$, let
$$
 S = \{ 1,\ 3,\ 4,\ 5,\ 6,\ 7,\ 8,\ 10,\ 17\}.
$$
 Then again $(K,S)$ is a degenerate CM-type.  This example is due to
Serre.

  Finally, let $p$, $q$, $r$ be odd primes, let $G = \ZZ/2pqr\ZZ$ as
cyclic group, and let $K$ be an extension of $\QQ$ with $\Gal(K/\QQ)\simeq
G$.  Then let $S$ be the subset of elements having order $1$, $pqr$, $2p$,
$2q$, $2r$, $2pq$, $2pr$ or $2qr$.  Then $(K,S)$ is a simple CM-type and
$$
 \rank(K,S) = 1 + pqr - (p-1)(q-1)(r-1) .
$$
 This example is due to Lenstra.

  All of these examples are verified in \cite{B.92} using
Proposition~9.4.1 and exhibiting odd characters $\chi$ such that
$\sum_{s\in S} \chi(s) = 0$.
\endexample

\example{9.4.3.  Degenerate CM-types in composite dimension}
 In Theorem~6.3.2 we saw that an abelian variety of CM-type with prime
dimension is always nondegenerate.  Dodson has proved a number of theorems
exhibiting the existence of degenerate abelian varieties in composite,
i.e., non-prime, dimensions \cite{B.28} \cite{B.29}
\cite{B.30}, and more recently Yanai has come up with a method for
generating degenerate CM-types that encompasses some of the previous
examples
\cite{B.140}.

\proclaim{Theorem{\rm (\cite{B.28}~Thm.3.2.1)}}
 When $n$ is composite, $n>4$, then there exist abelian varieties of
CM-type with $\dim A =n$ and rank $n-l+2$ for a divisor $l\ge 2$ of $n$
such that $n/l >2$.
\endproclaim

\proclaim{Theorem{\rm (\cite{B.28}~Thm.3.2.2)}}
 Let $p$ be a prime and suppose the ideal $(2)$ decomposes in the
cyclotomic field $\QQ(\zeta_p)$ into $g>1$ factors of degree~$f$.  Note
that $p\cong \pm1 \pmod 8$ is sufficient but not necessary to insure
$g>1$.  Then there exist abelian varieties of CM-type having dimension
$2^{fs}p$ and rank~$p-1$, for $0\le s\le g-1$.  In particular, for $s>0$
these abelian varieties are simple and degenerate.
\endproclaim

  In the following $K_0^{\Gal}$ denotes the Galois closure of the field
$K_0$, and $[(\ZZ/2\ZZ)^m]^+$ is the even subgroup of $(\ZZ/2\ZZ)^m$

\proclaim{Theorem{\rm (\cite{B.28}~Thm.3.3.1)}}
 \rom{(A)} Suppose there exists a totally real field $K_0^{\Gal}$ with
Galois group isomorphic to the wreath product $(\ZZ/k\ZZ) \wr (\ZZ/l\ZZ)$,
with $k>2$ and $l\ne 1$.  Then there exist simple degenerate abelian
varieties with complex multiplication by a non-Galois CM-field such that
the varieties have dimension~$kl$ and rank~$n-l+2$.
\par
\rom{(B)} Under the same hypotheses except that $k=2$ is allowed, there
exist simple abelian varieties of CM-type with dimension~$n' =k^2l$ having
rank $\le kl+1$.
\par
\rom{(C)} Further, in the even case, the existence of totally real fields
with Galois groups $[(\ZZ/2\ZZ)^m]^+\rtimes \ZZ/m\ZZ$ with $m\ge 3$,
respectively $D_m$ with $m>5$ and odd, supplies simple abelian varieties
of CM-type with dimension $n=4m$, respectively $n=2m$, and rank $\le \frac
n 2 +1$
\endproclaim

\proclaim{Theorem{\rm (\cite{B.29}~Thm.3.1)}}
 Let $d$ be a composite number.  Then there exist simple abelian varieties
of
CM-type with dimension ~$d$ and rank $\le (d/2)+1$ whenever $d$ is
\roster
\item even;
\item divisible by a square;
\item of the form $d = \binom p m$ with $0\le m \le (p-1)/2$; or
\item of the form $d= pt$ with $t\mid (p-1)$ and $t< p-1$.
\endroster
 Further, the rank $p+1$ occurs in dimension~$d$ at least in the following
cases:
\roster
\item when $d$ and $p$ are as in ~3 or ~4 above;
\item when $d= 2^{(p-1)/2}pt$ with $t\mid (p-1)/2$;
\item when $d=p q^2$ where $q$ is an odd prime and $p= q^2+ q +1$.
\endroster
\endproclaim

 The following is a nondegeneracy result

\proclaim{Theorem{\rm (\cite{B.30}~Thm.2.1)}}
 Let $n$ be odd and suppose that $K$ is a CM-field of degree $2n$ such
that the maximal totally real subfield $K_0$ has $\Gal(K_0^{\Gal}/\QQ)$
isomorphic to the symmetric groups or the alternating group on
$n$~letters.  Then every primitive CM-type $(K,S)$ is nondegenerate.
\endproclaim

 Recently Yanai has developed a method for generating degenerate CM-types
in higher dimension starting with degenerate CM-types in lower dimension.

\proclaim{Theorem {\rm (\cite{B.140})}}
 Let $K$ be a CM-field with $[K:\QQ] = 2d$ and let $K_1$ be a proper
subfield of $K$ with $[K_1:\QQ] =2d_1$.  Further, let $\pi: X_K \to
X_{K_1}$ be the canonical surjection from the character group of
$\Res_{K/\QQ}\Gm_{/K}$ to the character group of
$\Res_{K_1/\QQ}\Gm_{/K_1}$.  Suppose that the CM-types $(K,S)$ and
$K_1,S_1)$ satisfy the condition
$$
 \pi\big(\sum_{\sigma\in S}\sigma\big) = a \sum_{\sigma\in S_1}\sigma +
b\sum_{\sigma\in S_1}\bar\sigma
$$
 with some nonnegative integers $a$ and $b$ such that $a+b =[K:K_1]$.
Then
$$
 d+1 -\rank S \ge d_1 +1 -\rank S_1 .
$$
 Moreover, if $a=b$ then
$$
 d+1 - \rank S \ge d_1.
$$
 In particular, if the CM-type $(K_1,S_1)$ is degenerate or if $a=b$ then
the CM-type $(K,S)$ is degenerate.
\endproclaim
\endexample

\example{9.4.4. Lower bounds for CM-types}
 Ribet \cite{B.92} and Mai \cite{B.72} have given some lower
bounds for the rank of a CM-type, and particularly the latter discusses
how sharp these might be.

\proclaim{Proposition {\rm (\cite{B.92})}}
  $\rank(K,S) \ge 2 + \log_2 (\dim A)$
\endproclaim

  Mai considers the case where the CM-filed $K$ is Galois over $\QQ$.  For
the next proposition, note that when $V =\ZZ[\Gal(K/\QQ)]\tensor\CC$ is
considered as a $\Gal(K/\QQ)$-module, there is a decompositon $V=
\bigoplus_\pi d_\pi V_\pi$, where $\pi$ ranges over the irreducible
representations of $\Gal(K/\QQ)$ and $d_\pi = \dim V_\pi$.  A
representation $\pi$ is called \dfn{odd} if the value of its character at
complex conjugation is~$-1$.

\proclaim{Proposition {\rm(\cite{B.72}~Prop.1)}}
 When $K/\QQ$ is a Galois extension and $(K,S)$ is a simple CM-type, then
$$
 \rank (K,S) \ge 1 + \sum d_\pi ,
$$
 where the sum ranges only over those odd irreducible representaions $\pi$
such that $\pi(\sum_{s\in S}s) \ne 0$.
\endproclaim

\proclaim{Proposition {\rm(\cite{B.72}~Prop.2)}}
 When $K/\QQ$ is a Galois extension and $(K,S)$ is a simple CM-type, then
$$
 \rank(K,S) \le \max \Big\{ \frac{(p-1)^2\alpha}{p} : p \text{ an odd
prime, and } p^\alpha \,\big\|\, ([K:\QQ]/2) \Big\} .
$$
\endproclaim

 The notation $p^\alpha \,\big\|\, N$ means $p^\alpha$ exactly
divides~$N$, i.e., $p^\alpha \,\big|\,N$ and $p^{(\alpha+1)}
\not{\big|}\,N$.

\smallpagebreak

 In the following $S_a$ is the same CM-type that occured the first
paragraph of 9.4.2.  Such CM-types occur among factors of the Jacobians of
Fermat curves.

\proclaim{Proposition {\rm(\cite{B.72}~Prop.3)}}
 Let $K= \QQ(\zeta_p)$ and identify $\Gal(K/\QQ)$ with
$(\ZZ/p\ZZ)^\times$.  For $1\le a \le p-2$ let $S_a$ be the CM-type
defined by
$$
 S_a = \{g\in \Gal(K/\QQ) : 1 \le \angled g + \angled{ag} < p \} ,
$$
 where $1\le \angled g \le p-1$ and $\angled g \cong g \pmod p$.  Then
$$
 \rank(K,S) \ge 1 + \frac{19}{21} d , \qquad  d= \frac{p-1}2 .
$$
\endproclaim
\endexample

\subhead 9.5. Andr\'e's description of CM-Hodge cycles as Weil cycles
\endsubhead
  Finally we turn to a recent result of Andr\'e \cite{B.9} that
every Hodge cycle on an abelian variety $A$ of CM-type is a linear
combination of inverse images under morphisms $A\to B_J$ of Weil-Hodge
cycles on various abelian varieties $B_J$ of CM-type.  The following
definition should be compared with 1.13.6, 4.10, and the discussion
preceding Theorem~4.9.

\definition{9.5.1. Definition}
 Let $A$ be an abelian variety, let $F$ be a CM-field contained in
$\EndoA$, and let $V= H^1(A,\QQ)$.  Then $A$ or $V$ is said to be \dfn{of
Weil type relative to $F$} if there exists an $F$-Hermitian form $\psi$ on
$V$ admitting a totally isotropic subspace whose dimension over $F$ is
$\frac 1 2 \dim_F V$, and there exists a purely imaginary element $\alpha
\in F$ such that $\operatorname{Tr}_{F/\QQ}(\alpha\cdot\psi(u,v))$ defines
a polarization on~$A$.  Then the elements of $\twedge_F^{2p}V$ are called
\dfn{Weil-Hodge cycles relative to~$F$.}
\enddefinition

\proclaim{9.5.2. Theorem {\rm(\cite{B.9})}}
  Let $A$ be an abelian variety of CM-type, and $p$ a positive integer.
Then there exists a CM-field $F$, a finite number abelian varieties $A_J$
with complex multiplication of Weil type relative to $F$, and morphisms
$A\to A_J$, such that every Hodge cycle $\xi \in \Hdg^p(A)$ is a sum of
inverse images of Weil-Hodge cycles $\xi_J \in \Hdg^p(A_J)$.
\endproclaim

\demo{Proof}
 We sketch Andr\'e's proof.  Up to isogeny write $A = \prod_i A_i$ as a
product of simple CM-abelian varieties, where $A_i$ is of CM-type
$(K_i,S_i)$.  Let $V = H^1(A,\QQ)$, let $V_{S_i} = H^1(A_i,\QQ)$, and let
$F$ be the Galois closure of the compositum of all the $\Endo(A_i)$.  Then
$V = \bigoplus_{i\in I} V_{S_i}$ and $V_{S_i} \tensor F =
\bigoplus_{\sigma \in \Hom(F,\CC)} V_{S_i,\sigma}$.  Then
$$ \spreadlines{1\jot} \align
 \big( \twedge_\QQ^{2p} \big) \tensor F & \simeq \sum_{\sum d_i = 2p}
\big( \bigotimes_{i\in I} \twedge^{d_i} V_{S_i}\big) \tensor F \\
  & \simeq \sum\Sb \sum d_{i,\sigma} =2p \\ d_{i,\sigma} \in \{0,1\}
\endSb \bigotimes_{(i,\sigma)\in I\times \Hom(F,\CC)}
V_{S_i,\sigma}^{\tensor d_{i,\sigma}} .
\endalign
$$
  Let $T_F = \Res_{F/\QQ}\Gm_{/F}$.  Then the action of $(T_F)^I$ on
$\twedge^{2p}_\QQ V$ commutes with the action of $\Hg(V)$.  Further, the
action of $T_F$ can be extended by $F$-linearity to an action on
$\bigotimes V_{S_i,\sigma}^{\tensor d_{i,\sigma}}$.  It follows that every
Hodge cycle $\xi \in \twedge^{2p}_\QQ$ can be written as $\xi = \sum
\lambda_J \theta_J$, where $J$ indexes the set of sequences
$(d_{i,\sigma})_{(i,\sigma)\in I\times \Hom(F,\CC)}$ with $d_{i,\sigma}
\in \{0,1\}$ and $\sum d_{i,\sigma} =2p$, and where $\lambda_J \in F$, and
$\theta_J \in \bigotimes_{(d_{i,\sigma})\in J} V_{S_i,\sigma}^{\tensor
d_{i,\sigma}}$, and the restriction of $\Hg(V)$ acting on $F$
fixes~$\theta_J$.

  Now observe that each $\tau \in \Aut(F)$ induces an isomorphism of
rational Hodge structures $V_{S_i,\sigma} \to V_{\tau S_i, \sigma\tau}$,
although this isomorphism does not respect the action of $T_F$.  Therefore
we may write $V_J = \sum_{(d_{i,\sigma})\in J}
V^{d_{i,\sigma}}_{S_i,\sigma}$, where $S_{i,\sigma} = \sigma ^{-1} S$, in
such a way that via the isomorphisms induced by $\Aut(F)$ we get a
morphism of rational Hodge structures of CM-type $V_J \to V$ and
$\theta_J$ comes from an element $\zeta_J \in \bigotimes_{j\in J}
V_{S_{i,\sigma},\id}^{\tensor d_{i,\sigma}} \subset \big( \twedge_F^{2p}
V_J\big)\tensor_\QQ F$ which is invariant under the action of $\Hg(V_J)$
on~$F$.

  Note also that there is a natural basis $\chi_{(j,\sigma)}$ of
characters of $T_{F^J}$, where $(j,\sigma)$ runs over $J\times
\Hom(F,\CC)$.  Let $\gamma_{(j,\sigma)}$ denote the dual basis of
cocharacters.  Then if $\zeta_J \ne 0$, it generates the character
$\sum_{j\in J} j \chi_{(j,\sigma)}$ of~$T_{F^J}$.

  On the other hand, the Hodge structure of $V_J$ is determined by the
cocharacter $h:\U(1) \to (T_{F^J})_{/\RR}$, whose complexification may be
written out as
$$
 h_\CC = \sum_{(j,\sigma) \in J\times \Hom(F,\CC)} j(2S_j(\sigma) -1)
\gamma_{(j,\sigma)},
$$
 where $S_j(\sigma)$ is $1$ or $0$ according as $\sigma\in S_j$ or not.

  Then the fact that $\zeta_J$ is $\Hg(V_J)$-invariant implies that
$\sum_{j\in J} j \angled{\tau h, \chi_{(j,\id)}} =0$ for all $\tau \in
\Gal(F/\QQ)$.  Then expanding this expression, we find
$$
 \sum_{j\in J} j(2S_j(\sigma) -1) =0 .
$$
 And since $\sum_{j\in J} j = 2p$, it follows that
$$
 \sum_{j\in J} j S_j({\ssize{\bullet}}) = p ,
$$
 which implies that $V_J$ is a rational Hodge structure of Weil type.
Moreover,
$$
 \xi_J := [F:\QQ]^{-1} \sum_{\tau \in \Gal(F/\QQ)} \lambda^\tau
\zeta_J^\tau \in \twedge_F^{2p} V_J
$$
 is a Weil-Hodge cycle, and $\xi$ is a sum of images of $\xi_J$ under the
maps $V_J \to V$, since $\xi = [F:\QQ]^{-1} \sum_{\tau \in \Gal(F/\QQ)}
\lambda^\tau \theta^\tau_J$.
\Qed
\enddemo



\head 10. The general Hodge conjecture
\endhead
  In this section, when we speak of the general Hodge conjecture we always
mean the Grothendieck amended version as in~(7.12) of the text, that the
$r^{\text{th}}$ step of the arithmetic filtration $F_a^rH^i(A,\QQ)$ is the
largest rational Hodge structure contained in $F^rH^i(Z,\CC) \cap
H^i(A,\QQ)$, where $F^rH^i(A,\CC)$ is the Hodge filtration.  In those
cases where the stronger statement that $F_a^rH^i(A,\QQ) = F^rH^i(Z,\CC)
\cap H^i(A,\QQ)$ we will speak of Hodge's original conjecture, or the
strong form of the general Hodge conjecture.  Based on the results
assembled below, it would seem that when it is true, this stronger version
is more amenable to being proved.

\example{10.1. The general abelian variety}
  The earliest results about the general Hodge conjecture for abelian
varieties are those of Comessatti \cite{B.23} and Mattuck
\cite{B.73}, which show that Hodge's original conjecture is true
for the general abelian variety described in~1.13.8.  Mattuck's proof
proceeds by induction and explicit computation with period matrices.
Since the general $g$-dimensional abelian variety $A$ has
$$
 \Hg(A) = \Sp(H^1(A,\QQ),E) = \Lf(A),
$$
 (see~2.14) and also $\EndoA = \QQ$, it satisfies the hypotheses of
Theorem~6.2.  Thus $\Hdg(A^k) = \Div(A^k)$ for all $k\ge 1$, which is to
say that the general abelian variety is stably nondegenerate, and
Theorem~10.9 below applies.
\endexample

\subhead Abelian varieties of low dimension
\endsubhead
\nopagebreak
\example{10.2. Various abelian threefolds}
 The first interesting case of the general Hodge conjecture is
$\GHC(1,3,X)$ for a threefold $X$, and indeed, it was a special abelian
threefold, the product of three copies of an elliptic curve whose period
satisfies a cubic relation, that Grothendieck exhibited to show that
Hodge's original conjecture needed to be modified, see (7.5) in the text
or \cite{B.43}.  In \cite{B.12} Bardelli took up the
question of whether Grothendieck's counterexample to Hodge's original
conjecture satisfies Grothendieck's amended version, and at the same time
he considered a number of other abelian threefolds.

\proclaim{10.2.1. Theorem {\rm (\cite{B.12}~Prop.3.8)}}
  The Grothendieck generalized Hodge conjecture holds for:
\roster
\item The generic abelian threefold;
\item The generic member of the family of Jacobians of smooth genus three
curves admitting a morphism onto some elliptic curve;
\item The generic member of the family of Jacobians of smooth genus three
curves admitting a morphism onto some genus two curve;
\item The generic product of an abelian surface and an elliptic curve;
\item The generic product of three elliptic curves;
\item All products of three copies of the same elliptic curve, in
particular when the period $\tau$ is quadratic over $\QQ$ (the CM case) or
cubic over $\QQ$ (as in Grothendieck's counterexample).
\endroster
\endproclaim

  Bardelli's arguments are very geometric in nature.   To give a hint of
their flavor, let $J_H(A)$ denote the maximal subtorus of the intermediate
Jacobian of an abelian threefold $A$ that is orthogonal, with respect to
cup product, to $H^{3,0}(A)$.  Then the key lemma has the following form.

\proclaim{10.2.2. Lemma {\rm (\cite{B.12}~Lem.2.2)}}
  Let $T$ be an irreducible analytic subvariety of the Siegel upper
half-space of genus three and let $p:\script A \to T$ be the restriction
of the universal family of principally polarized abelian varieties over
the Siegel half-space.  Let $t_0 \in T$ be a generic point of $T$ at which
$T$ is smooth.  Then
$$
 \dim_\CC J_H(A_{t_0}) \le 9 - \dim T.
$$
\endproclaim

  Recall (1.13..8) that the Siegel half-space of genus~$g$ consists of
complex symmetric $g\times g$ matrices with positive-definite imaginary
part, and thus represents the possible complex structures for an abelian
variety of dimension~$g$.

\medpagebreak

  Another example of abelian threefold for which the general Hodge
conjecture is true comes up in \cite{B.86}.  In the course of
constructing a counterexample to a conjecture of Xiao, Pirola finds the
following.

\proclaim{10.2.3. Proposition}
 The general Hodge conjecture is true for the generic member of the family
of abelian threefolds of the form $W=\operatorname{Jac}(C)/f^*E$, where
$C$ is a smooth genus~$4$ curve, $E$ is an elliptic curve, and $f:C\to E$
is a $3$~to~$1$ cyclic Galois covering.
\endproclaim
\endexample

\example{10.3. An abelian fourfold of Weil type}
  In \cite{B.105} Schoen proves the general Hodge conjecture for
abelian fourfolds $A$ of Weil type with multiplication by $\EndoA = K
=\QQ(i)$ (and determinant~$1$, associated to a Hermitian form of
signature~$(3,1)$).  Previously, in \cite{B.104}, see Theorem~4.12,
he had proved the Hodge $(2,2)$ conjecture for these fourfolds.  Thus in
\cite{B.105} it only remained to verify $\GHC(1,4,A)$.  Therefore
the focus is on the rational Hodge substructure $U'\subset \twedge^4_K
H^1(A,\QQ)$, where $U'$ is the unique $\Res_{K/\QQ}\Gm_{/K}$
subrepresentation of $H^4(A,\QQ)$ which after tensoring with $\CC$ becomes
isomorphic to the sum of weight spaces $\alpha^4 \oplus \bar\alpha^4$.
The Hodge type of $U'$ is $\{(3,1),\,(1,3)\}$.  The highly geometric
arguments are lengthy and intricate, so we will not go into them here,
except to say that one of the main points is that the generic $A$ as above
is a generalized Prym variety associated to a cyclic $4$-fold covering
$\pi: C\to X$ of curves.  See \cite{B.105} for the details.
\endexample

\example{10.4. Powers of elliptic curves or abelian surfaces with
quaternionic multiplication}
  The simplest elliptic curves to deal with are those which entirely avoid
the kind of problem in Grothendieck's counterexample to Hodge's original
conjecture, namely where the period $\tau$ is either quadratic over $\QQ$,
in which case the elliptic curve has complex mutliplication, or a general
elliptic curve, whose period $\tau$ is transcendental over $\QQ$.

\proclaim{10.4.1. Proposition {\rm (\cite{B.117})}}
 If $E$ is an elliptic curve with complex multiplication, then Hodge's
original conjecture is true for $E^k$, for all $k\ge 1$.
\endproclaim

\proclaim{10.4.2. Proposition {\rm (\cite{B.38})}}
 If $E$ is a general elliptic curve, then Hodge's original conjecture is
true for $E^k$, for all $k\ge 1$.
\endproclaim

\proclaim{10.4.3. Proposition {\rm (\cite{B.39})}}
 If $A$ is a general abelian surface with quaternionic multiplication,
Hodge's original conjecture is true for $A^k$, for all $k\ge 1$.
\endproclaim

\proclaim{10.4.4. Proposition {\rm(\cite{B.4}~Thm.6.1)}}
 When $A$ is a product of elliptic curves, then the general Hodge
conjecture is true for $A$.
\endproclaim

  The first of these results is discussed in the text (7.18)--(7.20), and
the next two are supersceded by Theorem~10.9 below.  For the last, the
multiplicativity of the Lefschetz group (Lemma 2.15) reduces the problem
to powers of a single elliptic curve, see section three.  Then the case
where the curve is of CM-type is covered by 10.4.1 above, whereas the case
where the curve is not of CM-type is a special case of Theorems~10.9
or~10.12 below.
\endexample

\subhead Abelian varieties with conditions on endomorphisms, dimension or
Hodge group
\endsubhead
  It is only quite recently that results about the general Hodge
conjecture for abelian varieties of comparable generality to what has been
proved for the usual Hodge conjecture have begun to appear.  Here we
collect together the main results before discussing some of what is
involved in proving them.

\proclaim{10.5. Theorem {\rm (\cite{B.127}~Thm.1)}}
  If $A$ be a simple abelian variety of type~\rom{(I)} such that $\dim A /
[\EndoA :\QQ]$ is odd, then Hodge's original conjecture holds for~$A$.
\endproclaim

\proclaim{10.6. Theorem {\rm (\cite{B.127}~Thm.2)}}
  Let $A$ be a simple abelian variety of CM-type, with $\EndoA = K$ and
$K_0$ the maximal totally real subfield of~$K$.  If $[K^{\Gal}:
K_0^{\Gal}] = 2^{\dim A}$, then Hodge's original conjecture holds for~$A$.
\endproclaim

\proclaim{10.7. Theorem {\rm (\cite{B.128}~Thm.1)}}
  Let $A$ be an abelian variety with $\EndoA =\QQ$.  If
\roster
\item $\dim A \ne 4^l$,
\item $\dim A \ne \frac 1 2 \binom{4l+2}{2l+1}^{2m-1}$,
\item $\dim A \ne 2^{8lm+4l-4m-3}$,
\item $\dim A \ne 4^l (m+1)^{2l+1}$,
\item $\dim A \ne 2^{8ln +2n-4l-2} (8l+4)^{m-1}$,
\endroster
  for any positive integers $l$, $m$, $n$, then the general Hodge
conjecture holds for~$A$.  Furthermore, $\Hdg(A) = \Div(A)$, and $\Hg(A) =
\Sp(H^1(A,\QQ),E)$.
\endproclaim

\proclaim{10.8. Theorem {\rm (\cite{B.129}~Thm.1.1)}}
  Let $A$ be a simple complex abelian variety of dimension~$g$ with Hodge
group $\Hg(A)$, let $\hg(A,\CC) = \Lie \Hg(A)\tensor \CC$, and let
$\hg(A,\CC)_{\text{ss}}$ be the semisimple part of the reductive Lie
algebra $\hg(A,\CC)$.  Consider the following sets of natural numbers:
$$
 \operatorname{Ex}(1) := \Big\{ 4^l,\, \frac
12\binom{4l+2}{2l+1}^{2m-1},\, 2^{8lm+4l-4m-3},\, 4^l(m+1)^{2l+1} :
l,m\in\ZZ_+ \Big\};
$$
$$\multline
 \operatorname{Ex}(3) := \Big\{ 46{l+1},\, 6^{l+1},\,
\binom{4m+4}{2m+2}^l,\, \binom{4m+2}{2m+1}^{2l},\, 2^{(4m-1)}l ,  \\
  4^l(m+2)^{2l},\, 2^{l+1}(m+4)^{l+1} : l,m\in\ZZ_+ \Big\} ;
\endmultline
$$
$$\multline
 \operatorname{Ex}(4) := \Big\{ \binom{l+2}{m} \text{ for }1<m <(l+2)/2 ,
\\
  \binom{l+2}{m}^{n+1} \text{ for }1\le m <(l+2)/2 : l,m,n\in\ZZ_+\Big\}.
\endmultline
$$
\roster
\item If $\End(A)\tensor \RR =\RR$ and $g\notin \operatorname{Ex}(1)$,
then $\hg(A,\CC) = \frak{sp}(2g)$ and the general Hodge conjecture is true
for $A^k$, for $k\ge 1$.
\item If $\End(A)\tensor \RR = M_2(\RR)$ and $g\notin
2\cdot\operatorname{Ex}(3)$, then $\hg(A,\CC) = \frak{sp}(g)$ and the
general Hodge conjecture is true for $A^k$, for $k\ge 1$.
\item If $\End(A)\tensor \RR = \Bbb H$, the Hamiltonian quaternion
algebra, and $g\notin \operatorname{Ex}(3)$, then $\hg(A,\CC) =
\frak{so}(g)$ and for $0\le r \le g$ we have
$$
 \dim_\QQ \Hdg^r = \cases 1&\text{if }r\neq g/2,\\ g+2&\text{if }r=g/2
\endcases
$$
 (in particular, if $r\neq g/2$, then $\Hdg^r = \Div^r$).
\item If $\End(A)\tensor \RR =\CC$ and $g\notin \operatorname{Ex}(4)$,
then $\hg(A,\CC)_{\text{ss}} = \frak{sl}(q)$ and for all integers $r\neq
g/2$ we have $\Hdg^r =\Div^r$; in the case $rg/2$ we have the relations
$$
 \dim \Hdg^r = \cases 1&\text{if $\hg(A,\CC)$ is not semisimple,}\\
   3&\text{if $\hg(A,\CC)$ is semisimple.}\endcases
$$
\endroster
\endproclaim

\proclaim{10.9. Theorem {\rm (\cite{B.50}~Thm.5.1)}}
  If $A$ is a stably nondegenerate abelian variety (see 7.5 and 7.6) all
of whose simple components are of type~\rom{(I)} or~\rom{(II)}, the
general Hodge conjecture holds for~$A$ and all powers $A^k$ of $A$, for
$k\ge 1$.
\endproclaim

\remark{10.10. Remarks}
  1.\enspace It follows from Theorem~6.3.1 that a simple abelian variety
of type~\rom{(I)} such that $\dim A / [\EndoA :\QQ]$ is odd is stably
nondegenerate, so Theorem~10.5 is a special case of Theorem~10.9.  As
remarked already, Propositions~10.1, 10.4.2 and 10.4.3 are included in
10.9, as well.  Other examples of stably nondegenerate abelian varieties
include arbitrary products of elliptic curves (section~3), simple abelian
varieties of prime dimension (Theorem~6.3), or simple abelian varieties of
odd dimension without complex multiplication (Theorems~6.3 and~7.5).  Also
an abelian variety $A$ is stably nondegenerate if and only if $A^k$ is
stably nondegenerate, for any $k\ge 1$, by~7.6.1.2.

 2.\enspace Proposition~10.4.1 is a special case of Theorem~10.6.

 3.\enspace The methods used in the proof of Theorem~10.6 are similar to
those in the proof of Theorem~9.2 \cite{B.88}.

 4.\enspace The proof of Theorem~10.7 uses the classification result
Theorem~2.11.  Since $\EndoA=\QQ$, Theorem~2.7 applies, and over $\Qbar$
(or~$\CC$) the universal cover of $\Hg(A)$ is isomorphic to some number of
copies of an almost simple $\Qbar$-group, say~$G_1$.  Then if $G_1$ is
any of the types in Theorem~2.11 other than $\frak{sp}(2d)$, where
$d=\dim A$, then $d$ is one of the forbidden dimensions.  Then the known
representation theory of $\frak{sp}(2d)$ can be used to control the level
of sub-Hodge structures, in a similar spirit though by a different
argument as in the next paragraph.  The proof of Theorem~10.8 applies
similar ideas to the semisimple part of the Hodge group.
\endremark

 The proof of Theorem~10.9 provides an illustrative example of how the
representation theory of the symplectic group comes into proving the Hodge
conjecture.

\demo{10.11. Sketch of proof of Theorem~10.9 {\rm (following
\cite{B.50})}}
  Recall that the \dfn{level} $l(W)$ of a rational Hodge structure~$W$, in
particular a sub-Hodge structure of $H^m(A,\QQ)$, is the maximum of
$|p-q|$ for which $W^{p,q}\ne 0$.  Then it will suffice to prove that for
any irreducible rational sub-Hodge structure $W$ of $H^m(A,\QQ)$ with
$l(W) = m-2p$ there exists a Zariski-closed subset $Z$ of codimension~$p$
in $A$ such that
$$
 W \subset \Ker\{ H^m(A,\QQ) \to H^m(A-Z,\QQ)\} .
$$

  Now, it is a basic fact from the representation theory of
$\frak{sp}(2n,\CC)$ that there is a one-to-one correspondence between its
irreducible (finite-dimensional) representations and $n$-tuples
$(\lambda_1, \dots , \lambda_n)$ of nonnegative integers with
$\lambda_1\ge \lambda_2 \ge \dots \ge \lambda_n$, see \cite{B.33},
\cite{B.34}.  Such an $n$-tuple s called a \dfn{Young diagram of
length~$n$}.  Then the crucial proposition, whose proof we omit here, is
the following.

\proclaim{10.11.1.  Proposition {\rm (\cite{B.50})}}
 Let $A$ be an abelian variety with $\hg(A,\CC) \simeq
\frak{sp}(2n,\CC)$, let $W$ be an irreducible rational sub-Hodge structure
of $H^m(A^k,\QQ)$, and let $(\lambda_1,\dots,\lambda_n)$ be the associated
Young diagram.  Then
$$
 l(W) = \sum_{i=1}^n \lambda_i .
$$
\endproclaim

  For a Young diagram $(\lambda_1,\dots,\lambda_n)$ the number $\sum_i
\lambda_i$ is often referred to as \dfn{the number of boxes,} for the
traditional representation of a Young diagram as $n$ rows with $\lambda_i$
boxes in the $i^{\text{th}}$ row.  As a matter of notation, it is
convenient to write $(1^a)$ for the Young diagram with $\lambda_1 = \dots
= \lambda_a =1$ and $\lambda_{a+1} = \dots = \lambda_n =0$, and similarly
$(2^c,1^d)$ for the diagram with $\lambda_1 = \dots = \lambda_c =2$ and
$\lambda_{c+1} = \dots = \lambda_{c+d} = 1$ and $\lambda_{c+d+1} = \dots =
\lambda_n =0$.  By convention $(1^0)$ is the Young diagram for the trivial
representation.  Then we recall the following facts from representation
theory.

\proclaim{10.11.2. Lemma}
\roster
\item Let $V = \CC^{2n}$ as a standard representation of
$\frak{sp}(2n,\CC)$.  Then for $1\le i \le n$
$$
 \twedge^i V \simeq (1^i) \oplus (1^{i-2}) \oplus \dots \oplus
(1^{i-2[i/2]})  ,
$$
 where the inclusion of $(1^a)$ into $\twedge^i V$ is defined by taking
the exterior product $(i-a)/2$ times with
$$
 \Omega = \sum_{j=1}^n e_j\wedge e_{n+j} ,
$$
 where $\{e_1,\ldots,e_{2n}\}$ is a standard symplectic basis.
\item For nonnegative integers $a$, $b$ with $a\ge b$,
$$ \align
 (1^a) \tensor (1^b) \simeq\ & \{(1^{a+b}) \oplus (2, 1^{a+b-2}) \oplus
\dots \oplus (2^b, 1^{a-b})\} \\
 & \quad \oplus \{ (1^{a+b-2}) \oplus (2, 1^{a+b-4}) \oplus \dots \oplus
(2^{b-1}, 1^{a-b}) \} \\
 & \quad \oplus \dots \oplus \{ (1^{a-b}) \} ,
\endalign
$$
 with the convention that Young diagrams on the right-hand side with more
than $n$ rows are omitted.
\endroster
\endproclaim

\demo\nofrills
  See \cite{B.18}~Ch.VIII,\S13 and \cite{B.33}.
\enddemo

  Now the proof of the theorem is divided into three steps.  First,
consider the case where $A= B^k$, where $\hg(B,\CC) \simeq
\frak{sp}(2n,\CC)$ acting on $V= H^1(B,\CC) \simeq \CC^{2n}$ as a standard
representation.  Then any irreducible rational sub-Hodge structure $W$ in
$H^m(A,\QQ)$, say of level $l(W) = m-2p$, corresponds over $\CC$ to an
irreducible $\hg(B,\CC)$ representation (see Proposition~2.4) occuring in
one of the terms on the right-hand side of
$$
 H^m(A,\CC) \simeq \bigoplus_{m_1 +\dots + m_r =m} (\twedge^{m_1} V
\tensor \dots \tensor \twedge^{m_r} V) .
$$
  By Proposition~10.11.1 the number of boxes in the Young diagram
associated to $W_\CC$ is $m-2p$.  Then Lemma~10.11.2 implies that the
contraction, i.e., the reduction in the number of boxes, comes about only
by taking the exterior product with $\Omega^p$.  However, in the
dictionary between the representation theory and the cohomology, $\Omega$
corresponds to a divisor, say $D$, and taking the exterior product with it
corresponds to intersecting with~$D$.  Thus $W$ is the cup product of a
rational sub-Hodge structure in $H^{m-2p}(A,\QQ)$ with $D^p$, which
verifies the general Hodge conjecture in this case.

 Secondly, consider the case where $A = B_1^{k_1} \times B_2^{k_2}$ with
$B_1$ not isogenous to $B_2$, and as in the first case, $\frak g_i :=
\hg(B_i,\CC) \simeq \frak{sp}(2n_i,\CC)$ acting via a standard
representation on $V_i = H^1(B_i,\CC)$.  Then $\hg(A) \simeq \frak g_1
\times \frak g_2$, and any irreducible rational sub-Hodge structure
$W_\QQ$ of $H^m(A,\QQ)$ must correspond to an irreducible
$\hg(A,\CC)$-representation of the form $W_1\tensor W_2$ for some
irreducible $\hg(B_i,\CC)$-representation $W_i \subset
H^{m_i}(B_i^{k_i},\CC)$, with $m_1 + m_2 =m$.  Moreover, as in
Proposition~10.11.1, $l(W) = l(W_1) + l(W_2)$.  However, by the previous
case, if $l(W_i) = m_i -2p_i$, then $W_i$ is supported on a Zariski-closed
subset $Z_i$ of codimension $p_i$ on $B_i^{k_i}$.  Then $W_1\tensor W_2$
is supported on $Z_1\times Z_2$ of codimension $p_1 + p_2$ on $A$, which
verifies the general Hodge conjecture in this case.

  Finally, if $A$ is an arbitrary abelian variety which satisfies the
hypotheses of the theorem, then from \cite{B.47} and Theorem~2.11
it follows that
$$
 \hg(A,\CC) \simeq \frak{sp}(2n_1,\CC) \times \dots \times
\frak{sp}(2n_r,\CC)
$$
 acting in the standard way on
$$
 V_1^{\oplus k_1} \oplus \dots \oplus V_r^{\oplus k_r}.
$$
   Now, for each $i$ the fundamental form $\Omega_i \in \twedge^2 V_i$ is
$\frak{sp}(2n_i,\CC)$-invariant, thus by Lefschetz's theorem corresponds
to a linear combination of divisor classes.  Then arguing similarly as in
the previous paragraph for $r=2$ shows that the general Hodge conjecture
holds for~$A$, as was to be shown.
\Qed
\enddemo

\example{10.12. Passing from the usual Hodge conjecture to the general
Hodge conjecture}
 The last examples of abelian varieties for which the general Hodge
conjecture has been proved come through a theorem of Abdulali
\cite{B.4}.  In the following statement recall that the derived
group $G^{\text{der}} = (G,G)$ of a group $G$ is the (normal) subgroup
generated by all elements of the form $ghg^{-1}h^{-1}$.

\proclaim{Theorem {\rm (\cite{B.4})}}
 Let $A$ be an abelian variety whose Hodge group is semisimple and equal
to the derived group of the Lefschetz group of $A$, i.e.,  $\Hg(A) =
\Lf(A)^{\text{der}}$.  Further suppose that for every simple factor $B$ of
$A$ of type~\rom{(III)} the dimension of $H^1(B,\QQ)$ as a vector space
over $\Endo(B)$ is odd.  Then if the usual Hodge conjecture is true for
$A^k$ for all $k\ge 1$, then the general Hodge conjecture is also true for
$A$, and all $A^k$ for $k\ge 1$.
\endproclaim

\remark{10.12.1. Remarks}
 In the presence of the assumption that $\Hg(A)$ is semisimple, the
hypothesis that $\Hg(A) = \Lf(A)^{\text{der}}$ can be alternately
formulated
as follows:  There is a natural embedding $\Hg(A) \hra \Sp(H^1(A,\QQ))$
which induces a holomorphic embedding of the symmetric domain $D$ of
$\Hg(A,\RR)$ into the symmetric domain $\frak H$ of $\Sp(H^1(A,\RR))$.
Then the pull-back to $D$ of the universal family of polarized abelian
varieties of dimension $\dim A$ natually lying over $\frak H$ determines a
family of abelian varieties of \dfn{Hodge type} in the sense of
\cite{B.78}.  Then the hypothesis that $\Hg(A) =
\Lf(A)^{\text{der}}$, or in the absence of the assumption that $\Hg(A)$ is
semisimple, an assumption that $\Hg(A)^{\text{der}} =
\Lf(A)^{\text{der}}$,
is equivalent to requiring that the family of Hodge type be a family of
abelian varieties of PEL-type, in the sense of \cite{B.109}
\cite{B.110} \cite{B.111}.  That is to say, a family of
abelian varieties that is determined by polarization, endomorphism algebra
and level structures.

 The essential use of this hypothesis in the proof of the theorem is in
the multiplicativity of $\Lf(A)$.  If $A$ is isogenous to $A_1^{k_1}
\times \dots \times A_r^{k_r}$, then $\Lf(A) = \Lf(A_1) \times \dots
\times
\Lf(A_r)$, see Lemma~2.15, and thus under the assumptions at hand, $\Hg(A)
=
\Hg(A_1) \times \dots \times \Hg(A_r)$.  This makes it possible to reduce
the proof to the case where $A$ is isogenous to $A_0^k$ for a simple
abelian variety $A_0$.  The proof then proceeds by cases, according to
whether a simple factor of $\Endo(A_0) \tensor \RR$ is $\RR$, or $\CC$ or
$\Bbb H$.
\endremark

\example{10.12.2. Applications}
  What abelian varieties satisfy the hypotheses of Theorem~10.12?  A
stably nondegenerate abelian variety with a semisimple Hodge group cannot
have factors of type~(III) or~(IV), so this is the same class as covered
by Theorem~10.9.  However, Abdulali observes that whenever the usual Hodge
conjecture is true for an abelian four-fold $A$ of Weil type, then it is
true for all powers $A^k$ of $A$ \cite{B.4}.  Thus there is the
following consequence of Theorem~4.12.
\endexample

\proclaim{10.12.3. Corollary}
  The general Hodge conjecture is true for all powers of a general abelian
fourfold of Weil type $(A,K)$ with $K=\QQ(\sqrt{-3})$ or $K=\QQ(i)$, when
the determinant of the associated Hermitian form is~$1$.
\endproclaim
\endexample



\head 11. Other approaches to the Hodge conjecture
\endhead
  In this section we look at three conditional results on the Hodge
conjecture.

\example{11.1. Higher Jacobians}
 In \cite{B.98}, Sampson outlines one possible approach to proving
the Hodge conjecture for arbitrary abelian varieties.  Given an abelian
vareity $A$ over $\CC$, let $J^p(A)$ denote its $p^{\text{th}}$ Weil
intermediate Jacobian, for odd $p$ with $1 < p \le \dim A$.  Then Sampson
gives an explicit but complicated construction of a surjective
homomorphism $\pi: J^p(A) \to A$ which induces an isomorphism $f:
\Hdg^p(A) \to \Div^1(J^p(A))$.  By the Poincar\'e Reducibility Theorem
1.11.4, $J^p(A)$ splits, up to isogeny, as $\Ker(\pi) \times A'$.  Thus a
cycle class $[Z]$ of codimension~$r$ on $Ker(\pi)$ determines a cycle
class $\operatorname{proj}_A((Z\times A')\cdot f(\phi))$ of
codimension~$(r+1)$ on $A$, with $\phi \in \Hdg^p(A)$.  Now if we fix $Z =
H^{p-1}$ to be the $(p-1)$-fold self-intersection of a fixed hyperplane
section, then
$$
 \phi \mapsto f^*(\phi) := \operatorname{proj}_A((H^{p-1}\times A')\cdot
f(\phi))
$$
 defines a homomorphism from $\Hdg^p(A)$ into the group of cohomology
classes of algebraic cycles on $A$ of codimension~$p$.  Then it is not
hard to show that if $f^*$ were injective, then the Hodge $(p,p)$
conjecture would follow.  However, the highly transcendental nature of the
construction makes the connection between $f^*(\phi)$ and $\phi$ rather
obscure, as well as apparently making it very difficult to determine
whether $f^*$ is injective.
\endexample

\example{11.2. The Tate conjecture}
 The references for this are \cite{B.88}
\cite{B.87} \cite{B.27}, and see
also \cite{B.15} \cite{B.16} for related results.
If $A$ is a
complex abelian variety, then there is a subfield $F\subset \CC$ finitely
generated over $\QQ$ and a model $A_0$ of $A$ over $F$, meaning that $A =
A_0 \tensor_F \CC$.  Then the $\ell$-adic \'etale cohomology of $A_0$ over
the algebraic closure $F^{\text{alg}}$ of $F$, that is,
$H_{\text{\'et}}^{2p}(A_0\tensor F^{\text{alg}}, \QQ_{\ell}(p))$, is
naturally a
$\Gal(F^{\text{alg}}/F)$-module.  In \cite{B.130} Tate conjectured that
the
elements of $H_{\text{\'et}}^{2p}(A_0\tensor F^{\text{alg}},
\QQ_{\ell}(p))$ fixed
by some open subgroup of $\Gal(F^{\text{alg}}/F)$, or equivalently by (the
Zariski-closure of) the $\ell$-adic Lie subalgebra $\frak g_\ell \subset
\End(H_{\text{\'et}}^{2p}(A_0\tensor F^{\text{alg}}, \QQ_{\ell}(p)))$
generated by
the image of $\Gal(F^{\text{alg}}/F)$, is precisely the $\QQ_\ell$-span of
the
classes of algebraic cycles.  In \cite{B.130} Tate himself observed
that this conjecture has an air of compatibility with the Hodge
conjecture, and already in \cite{B.77}, with the introduction of
the Hodge group, Mumford reported that Serre conjectured that
$$
 \frak g_\ell = \mt(A) \tensor \QQ_\ell .
$$
 For an excellent early introduction to this conjecture and the
relationships between the Tate and Hodge conjectures, see \cite{B.108};
the literature in the 20~years since then is extensive, it would
take another appendix at least the size of this one to survey it.

 The main result of \cite{B.88} is that for abelian varieties of
CM-type, the Hodge and Tate conjectures are equivalent.  Then that the
validity of the Tate conjecture for an abelian variety $A$ implies the
validity of the Hodge conjecture for $A$ has been proved by
Piatetskii-Shapiro \cite{B.87}, Deligne
(unpublished) and \cite{B.27}.  \cite{B.15}
extends the result of \cite{B.87}, and
\cite{B.16} contains a weaker version of the main theorem
of \cite{B.27}, from which Tate implies Hodge for abelian
varieties follows as a corollary.
\endexample

\example{11.3. Standard conjectures}
  In \cite{B.3} Abdulali shows that if one assumes
Grothendieck's invariant cycles conjecture \cite{B.42} for
families of abelian varieties of Hodge type in the sense of \cite{B.78},
then the Hodge conjecture for abelian varieties follows.  He also
formulates the $L_2$-cohomology analogue of Grothendieck's standard
conjecture~(A) that the Hodge $*$-operator is algebraic \cite{B.44},
and shows that for the families of abelian varieties being
considered that this conjecture implies the invariant cycles conjecture
and thus the Hodge conjecture for abelian varieties.

\proclaim{11.3.1. Conjecture {\rm (Invariant cycles conjecture
\cite{B.42})}}
 Let $f: A\to V$ be a smooth and proper morphism of smooth quasiprojective
varieties over~$\CC$.  Let $P\in V$ and let $\Gamma:= \pi_1(V,P)$.  Then
the space of $s\in H^0(V, R^bf_*\QQ) \simeq H^b(A_P, \QQ)^\Gamma$ that
represent algebraic cycles in $H^b(A_P, \QQ)^\Gamma$ is independent
of~$P$.
\endproclaim

\definition{11.3.2. Families of abelian varieties of Hodge type {\rm
(\cite{B.78})}}
 Let $A_0$ be a polarized abelian variety, let $W= H_1(A_0,\QQ)$, let $L =
H_1(A_0,\ZZ)$, and let $E$ be a Riemann form on $W$ representing the
polarization. Also, let $h:\U(1) \to \GL(\Wr)$ be the complex structure on
$\Wr$, let $K^+$ be the connected component of the centralizer of
$h(\U(1))$ in $\Hg(A_0,\RR)$, and let $D= \Hg(A_0,\RR)^+/K^+$ be the
bounded symmetric domain associated to $\Hg(A_0)$, as in~2.10.  Then to
each point $x\in D$ we can associate the polarized abelian variety $A_x =
(\Wr/L, ghg^{-1}, [E])$, where $x = g K^+$.  Further, if $\Gamma \subset
\Hg(A_0)$ is a torsion-free arithmetic subgroup that preserves $L$, then
$\gamma \in\Gamma$ induces an isomorphism between $A_x$ and $A_{\gamma
x}$.  Thus we get a family $\{A_x : x\in V\}$ of polarized abelian
varieties parameterized by $V=\Gamma \bs D$, which may be glued together
into an analytic space $A \to V$ fibered over~$V$.  Such a family of
abelian varieties is said to be of \dfn{Hodge type} \cite{B.78}.
Furthermore, $V$ has a canonical structure as a smooth quasiprojective
algebraic variety \cite{B.11}, and the analytic map $A\to
V$ is an algebraic morphism \cite{B.13}.
\enddefinition

\proclaim{11.3.3. Theorem {\rm (\cite{B.3}~Thm.6.1)}}
 If Conjecture~11.3.1 is true for all families of abelian varieties of
Hodge type, then the Hodge conjecture is true for all abelian varieties.
\endproclaim

\demo\nofrills
  An outline of the proof may be sketched as follows.  The first step is
to deduce from Conjecture~11.3.1 that all Weil-Hodge cycles are algebraic.
 To do this, Abdulali shows that any abelian variety $A_1$ of Weil type is
a member of a Hodge family whose general member $A_\eta$ has Hodge group
equal to the full symplectic group.  Then since $\Hdg(A_\eta) =
\Div(A_\eta)$, the invariant cycles conjecture implies that all Weil
cycles become algebraic in this family.  The next point is to observe that
Theorem~9.5.2 implies that if all Weil-Hodge cycles are algebraic, then
the Hodge conjecture is true for all abelian varieties of CM-type.
However, Mumford showed that every family of abelian varieties of Hodge
type contains members of CM-type \cite{B.78}.  Then the invariant
cycles conjecture can be used again to deduce that a Hodge cycle on any
member of the family is algebraic.  See \cite{B.3} for more
details.
\enddemo
\endexample

\newpage


\head Chronological listing of work on \\
 the Hodge conjecture for abelian varieties
\endhead
\rightheadtext{Chronological listing}

\medskip
{\baselineskip=13.5pt
\halign{\hfil#&\qquad#\hfil&\qquad#\hfil\cr
\smc Year&\smc Author &\smc Topic\cr
\noalign{\medskip}
1950&Hodge &Presented conjecture  \cr
1958&Mattuck &GHC for general abelian variety  \cr
1966&Mumford &Introduced Hodge group  \cr
1968&Polhmann &Hodge if and only if Tate for CM-type  \cr
1969&Mumford &Families of Hodge type  \cr
1969&Grothendieck &Amended general Hodge conjecture  \cr
1969&Murasaki &Elliptic curves  \cr
1971&Piatetskii-Shapiro &Tate implies Hodge  \cr
1974&Borovo\u\i &Tate implies Hodge  \cr
1976&Imai &Elliptic curves  \cr
1977&Serre &Connections between Hodge and Tate conj.  \cr
1977&Borovo\u\i &Absolute Hodge cycles  \cr
1977&Weil &Weil type  \cr
1978&Tankeev &$4$-dimensional abelian varieties  \cr
1979&Deligne &Classification of semisimple part of $\hg$  \cr
1979&Serre &Classification of semisimple part of $\hg$  \cr
1979&Tankeev &$4$- and $5$-dimensional  \cr
1981&Borovo\u\i &Simplicity of Hodge group  \cr
1981&Shioda &Fermat type  \cr
1981&Tankeev &Simple abelian varieties  \cr
1981&Tankeev &Simple ($5$-dimensional) abelian varieties  \cr
1982&Tankeev &Simple abelian varieties, prime dimension  \cr
1982&Deligne &Absolute Hodge cycles, Tate implies Hodge  \cr
1982&Kuga &Exceptional cycles  \cr
1982&Sampson &Alternate approach  \cr
1982&Ribet &Simple abelian varieties  \cr
1983&Shioda &Survey  \cr
1983&Ribet &Simple abelian varieties, Lefschetz group  \cr
1983&Hazama &Nondegenerage CM-type  \cr
1983&Murty &Non-simple abelian varieties  \cr
1983&Hazama &Non-simple abelian varieties  \cr
1984&Kuga &Exceptional cycles  \cr
1984&Hazama &Non-simple abelian varieties  \cr
1984&Murty &Non-simple abelian var., exceptional cycles  \cr
1984&Dodson &Degenerate CM-types  \cr
1985&Zarhin &Classification of $\hg$, survey  \cr
1985&Yanai &Nondegenerate CM-types  \cr
1986&Dodson &Degenerate CM-types  \cr
1987&Steenbrink &General Hodge conjecture, survey  \cr
1987&Dodson &Degenerate CM-types  \cr
1987&Bardelli &GHC, low dimension \cr
1987&Ichikawa &Non-simple abelian varieties MT groups  \cr
1988&Hazama &Stably degenerate  \cr
1988&Murty &Lefschetz group, semisimple part of $\Hg(A)$  \cr
1988&Schoen &Weil type  \cr
1988&Gordon &GHC for powers of general QM-surfaces  \cr
1989&Kuga &Exceptional cycles  \cr
1989&Hazama &Non-simple abelian varieties  \cr
1989&Schoen &GHC Weil type  \cr
1989&Mai &Degenerate CM-types  \cr
1990&Murty &Survey, Hodge group  \cr
1990&Kuga, Perry, Sah &Exceptional cycles  \cr
1991&Ichikawa &Non-simple abelian varieties, Hodge groups  \cr
1992&Pirola &GHC special threefolds  \cr
1992&Andr\'e &CM-type and Weil cycles  \cr
1993&White &Degenerate CM-type  \cr
1993&Gordon &GHC powers of general elliptic curve  \cr
1993&Tankeev &GHC  \cr
1994&Zarhin &Survey, connections with arithmetic  \cr
1994&van Geemen &Survey, Weil type  \cr
1994&Yanai &Degenerate CM-types  \cr
1994&Hazama &GHC  \cr
1994&Abdulali &Alternate approach  \cr
1994&Tankeev &GHC  \cr
1995&Moonen \& Zarhin &Abelian $4$-folds  \cr
1996&Lee &Exceptional cycles  \cr
1996&Abdulali &GHC  \cr
1996&Tankeev &GHC\cr
1996&Silverberg and Zarhin&Hodge group, connection to arithmetic\cr
1996&Moonen \& Zarhin &Exceptional Weil cycles  \cr
}}

\newpage



\enddocument